\documentclass[A4]{amsart}
\usepackage{amsmath, amssymb, amstext, amsfonts, textcomp, amsxtra, amsbsy, amsgen, amsopn, amscd, mathrsfs, amsthm, latexsym, array}
\usepackage[textwidth=16cm,textheight=22cm,centering]{geometry}

\usepackage[dvips]{graphicx}

\usepackage[all]{xy}

\newcommand{\N}{\mathbb{N}}
\newcommand{\Z}{\mathbb{Z}}
\newcommand{\Q}{\mathbb{Q}}
\newcommand{\R}{\mathbb{R}}
\newcommand{\C}{\mathbb{C}}
\newcommand{\A}{\mathbb{A}}

\newcommand{\vO}{\mathcal{O}}
\newcommand{\vp}{\mathfrak{p}}

\newtheorem{theorem}             {Theorem}  [section]
\newtheorem{definition} [theorem] {Definition}
\newtheorem{lemma}      [theorem]{Lemma}
\newtheorem{corollary}  [theorem]{Corollary}
\newtheorem{proposition}[theorem]{Proposition}
\newtheorem{remark} [theorem] {Remark}

\numberwithin{equation}{section} \everymath{\displaystyle}

\makeatletter
\newcommand{\rmnum}[1]{\romannumeral #1}
\newcommand{\Rmnum}[1]{\expandafter\@slowromancap\romannumeral #1@}
\makeatother

\title{Burgess-like subconvex bounds for $\text{GL}_2 \times \text{GL}_1$}
\author{Han Wu}
\thanks{Research partially supported by SNF-grant 200021-125291}

\begin{document}

\begin{abstract}
	Let $F$ be a number field, $\pi$ an irreducible cuspidal representation of ${\rm GL}_2(\A_F)$ with unitary central character, and $\chi$ a Hecke character of analytic conductor $Q$. Then $L(1/2, \pi \otimes \chi) \ll Q^{\frac{1}{2} - \frac{1}{8}(1-2\theta)+\epsilon}$, where $0 \leq \theta \leq 1/2$ is any exponent towards the Ramanujan-Petersson conjecture. The proof is based on an idea of unipotent translation originated from P.Sarnak then developped by Ph.Michel and A.Venkatesh, combined with a method of amplification.
\end{abstract}

\maketitle

\tableofcontents

\section{Introduction}

	\subsection{Statement of the Main Result}
	Let $\A$ be the adele ring of a number field $F$. Let $\pi$, $\pi_1$, $\pi_2$ be generic automorphic representations of $G(\A) = {\rm GL}_2(\A)$, where at least one of $\pi_1,\pi_2$ is cuspidal. Let $\chi$ be a Hecke character. Denote by $C(\pi)$ (resp. $C(\chi)$) the analytic conductor of $\pi$ (resp. $\chi$).
	
	Ph.Michel and A.Venkatesh \cite{MV} solved the subconvexity problem for ${\rm GL}_2$. In fact, the main result of that paper is the existence of some $\delta > 0$ such that
	$$ L(1/2, \pi_1 \times \pi_2) \ll_{F, \epsilon, \pi_1} C(\pi_2)^{1/4-\delta+\epsilon}, \forall \epsilon > 0. $$
	That is to say, if one fixes $\pi_1$, then we have subconvex bound for $L(1/2, \pi_1 \otimes \pi_2)$ as $C(\pi_2)$ tends to infinity. As a preliminary result, they also obtained the following subconvex bound
	$$ L(1/2, \pi \times \chi) \ll_{F, \epsilon, \pi} C(\chi)^{1/2-\delta+\epsilon}, \forall \epsilon > 0. $$

	The main result of this paper is to give an explicit value of $\delta$.
	\begin{theorem}
		Let $\theta$ be such that no complementary series with parameter $> \theta$ appear as a component of a cuspidal automorphic representation of $G(\A)$. For any cuspidal automorphic representation $\pi$  of $G(\A)$ and any Hecke character $\chi$ of analytic conductor $C(\chi) = Q$, we have
		$$ L(1/2, \pi \otimes \chi) \ll_{F, \epsilon, \pi} Q^{1/2-\delta+\epsilon}, \forall \epsilon > 0 $$
		with $$ \delta = \frac{1-2\theta}{8}. $$
		Note that under the Ramanujan-Petersson conjecture $(\theta = 0)$, $\delta = 1/8$.
	\label{Main}
	\end{theorem}
	\begin{remark}
		This bound, when $\theta = 0$, is called a Burgess bound. Burgess \cite{B03} first obtained such bounds for Dirichlet $L$-functions in the level aspect. The best known value $\theta = 7/64$ is due to Kim and Sarnak \cite{KS} over $\Q$, and to Blomer and Brumley \cite{BB} over an arbitrary number field.
	\end{remark}
	\begin{remark}
		In the hybrid aspect, this result is new even for $F = \Q$. The previous best known result is due to Munshi \cite{M}.
	\end{remark}
	\begin{remark}
		Blomer, Harcos and Michel \cite{BHM} first established such a Burgess-like bound in the level aspect for $F=\Q $. It was then generalized by Blomer and Harcos \cite{BH4} for any totally real number field $F$. The best bound for $F=\Q$, $\delta = 1/8$, in the level aspect was obtained in Theorem 2 of \cite{BH3} by Blomer and Harcos. In the case $F = \Q$, $\pi$ being of trivial central character and $\chi$ being quadratic, $\delta = 1/6$ was obtained by Conrey and Iwaniec as Corollary 1.2 of \cite{CI}.
	\end{remark}

	\subsection{Plan of the Paper}

	Section 2 is concerned with some technical but fundamental aspects of the proof of Theorem \ref{Main}:
	
	In Section 2.1 we provide notations and conventions. In Sections 2.2 to 2.4 we recall how Hecke's theory can be extended from $K$-finite vectors to smooth vectors. In Section 2.5 we discuss Whittaker models and their norms. In Sections 2.6 and 2.7, we discuss various forms of the spectral decomposition of automorphic functions. In Section 2.8 we use results from Section 2.5 to construct and study local test vectors to be used in the sequel. In Section 2.9 we discuss the decay of matrix coefficients of automorphic representations.
	
	In Section 3 we start the proof of Theorem \ref{Main}, setting up the amplification method. We split to two sorts of arguments: local ones and global ones. The intuition behind the formal calculations is explained in the beginning. It seems that the idea of translation by $n(T)$ originates from P.Sarnak \cite{Sa}. The whole idea is the combination of his idea together with the amplification method.
	
	In Section 4 we deal with the local arguments and prove the first part of Proposition \ref{Est}.

	In Section 5 we give the decay of matrix coefficients in the special case without $n(T)$ translation and concerning classical vectors. This complements Section 2.9 for our application.

	In Section 6 we conclude the proof by putting local estimations into the global arguments.
	
	The reader is strongly recommended to read the beginning of Section 3 before entering into the subsequent calculations. The difference in methods between this paper and \cite{MV} is explained in Remark \ref{DifferenceFromMV}.

	\subsection{Acknowledgement}
	This paper is part of my Ph.D thesis under the supervison of Prof. Philippe Michel at the Ecole Polytechnique F\'ed\'erale de Lausanne. It is a great pleasure to thank Professor Philippe Michel for his stimulating discussions and patient explanations with inspirational references. I am also very grateful for his continuous encouragement during the whole preparation of this paper. I am greatly in debt of Paul Nelson for his interest in my work and pointing out a mistake of my paper. Special thanks are due to the referee, whose careful reading and detailed remarks have significantly improved the paper.

\section{Some Preliminaries}

	\subsection{Notations and Conventions}
	From now on, $F$ is a number field of degree $r=[F:\Q]=r_1+2r_2$, where $r_1$ is the number of real places and $r_2$ is the number of pairs of conjugated complex places. $V_F$ is the set of all places of $F$. For any $v\in V_F$, $F_v$ is the completion of $F$ at the place $v$. $\A=\A_F$ is the adele ring of $F$. $\A^{\times}$ is the idele group. We fix once for all an isometric section $\R_+ \to \A^{\times}$ of the adelic norm map $|\cdot |: \A^{\times} \to \R_+$, thus identify $\A^{\times}$ with $\R_+ \times \A^{(1)}$, where $\A^{(1)}$ is the kernel of the adelic norm map. We will constantly identify $\R_+$ with its image under the section map. Let $F_{\infty} = \prod_{v | \infty} F_v$ and $F_{\infty}^{(1)} $ the subgroup of $F_{\infty}^{\times}$ of adelic norm $1$. $\A_{\text{f}}$ is the subring of finite adeles. $\A_{\text{f}}^{\times}$ is the unit group of $\A_{\text{f}}$.
	
	We denote by $\psi = \prod_v \psi_v$ the additive character $ \psi = \psi_{\Q} \circ {\rm Tr}_{F/\Q} $ of $\A_F$, where $\psi_{\Q}$ is the additive character of $\Q \backslash \A_{\Q}$ taking $x \mapsto e^{2\pi i x}$ on $\R$. At each place $v\in V_F$, $dx_v$ denotes a self-dual measure w.r.t. $\psi_v$. Note if $v<\infty$, then $dx_v$ is the measure which gives the ring of integers $\vO_v$ of $F_v$ the mass $q_v^{-d_v/2}$, where $q_v$ is the cardinality of the residue field of $F_v$, and $\prod_{v<\infty} q_v^{d_v}$ is the discriminant $\text{disc}(F)$ of $F$. We set $v(\psi)=-d_v$. Define $dx = \prod_{v\in V_F} dx_v$ on $\A$. The quotient measure on $F\backslash \A$ has total mass $1$ (c.f. Chapter $\Rmnum{14}$ Proposition 7 of \cite{L}).  Define for $s\in\C$, if $v$ is complex, $\zeta_v(s) = \Gamma_{\C}(s) = 2(2\pi)^{-s} \Gamma(s)$; if $v$ is real, $\zeta_v(s) = \Gamma_{\R}(s) = \pi^{-s/2} \Gamma(s/2)$; if $v<\infty$, $\zeta_v (s) = (1-q_v^{-s})^{-1}$. Take $d^{\times}x_v = q_v^{d_v/2} \zeta_v(1) \frac{dx_v}{|x_v|_v }$ as the Haar measure on the multiplicative group $F_v^{\times}$ if $v < \infty$, which gives $\vO_v^{\times}$ mass $1$; $d^{\times}x_v = \frac{dx_v}{|x_v|_v }$ if $v | \infty$. Define $d^{\times}x = \prod_v d^{\times} x_v$ as the Haar measure on the idele group $\A^{\times}$.

	Unless otherwise specified, $G={\rm GL}_2$ as an algebraic group defined over $F$. Hence $G_v = {\rm GL}_2(F_v)$. If $v$ is a complex place, then $K_v = {\rm SU}_2(\C)$; if $v$ is a real place, then $K_v = {\rm SO}_2(\R)$; if $v<\infty$ then $K_v = G(\vO_v)$. We also define
	\begin{align*}
		Z_v &= \left\{ z(u) = \begin{pmatrix} u & 0 \\ 0 & u \end{pmatrix} : u\in F_v^{\times} \right\}, \\
		N_v &= \left\{ n(x) = \begin{pmatrix} 1 & x \\ 0 & 1 \end{pmatrix} : x\in F_v \right\}, \\
		A_v &= \left\{ a(y) = \begin{pmatrix} y & 0 \\ 0 & 1 \end{pmatrix} : y\in F_v^{\times} \right\}.
	\end{align*}
	The probability Haar measure on $K_v$ is $dk_v$. $Z_v$(resp. $N_v$, resp. $A_v$) is equipped with the measure $d^{\times}u$ (resp. $dx$, resp. $d^{\times}y$). Consider the Iwasawa decomposition $G_v = Z_v N_v A_v K_v$, a Haar measure of $G_v$ is given by $dg_v =  d^{\times}u dx d^{\times}y/|y|_v dk_v$, which in fact gives $K_v \subset G_v$ the mass $q_v^{-d_v/2}$ for $v<\infty$. View $Z_v \backslash G_v$ as $N_v A_v K_v$, equipped with the measure $d\bar{g}_v =  dx d^{\times}y/|y|_v dk_v$. The center of $G(\A)$ is $Z = \prod_{v\in V_F} Z_v$. Denote $A = \prod_v A_v$. The quotient group $Z \backslash G(\A)$ is equipped with the product measure $d\bar{g} = \prod_{v\in V_F} d\bar{g}_v$ which gives $K_v$ the mass $q_v^{-d_v/2}$. The quotient measure on $X(F) = Z G(F) \backslash G(\A)$ is also denoted by $d\bar{g}$, with total mass $\text{Vol}(X(F))$. $K=\prod_{v\in V_F} K_v$ is equipped with the product measure $dk = \prod_v dk_v$. Write $K_{\infty} = \prod_{v| \infty}K_v$ and $K_{\text{f}} = \prod_{v<\infty } K_v$.

	Given a Hecke character $\omega$, $L^2(G(F) \backslash G(\A), \omega)$ is the space of Borel functions $\varphi$ satisfying
	$$ \forall \gamma \in G(F), \varphi(\gamma g) = \varphi(g);$$
	$$ \forall z\in Z, \varphi(z g) = \omega(z) \varphi(g); $$
	$$ \lVert \varphi \rVert_{X(F)}^2 = \int_{X(F)} |\varphi(\bar{g})|^2 d\bar{g} < \infty. $$
	Let $L_0^2(G(F) \backslash G(\A), \omega)$ be the (closed) subspace of cusp forms $\varphi \in L^2(G(F) \backslash G(\A), \omega)$ satisfying 
	$$ \int_{F\backslash \A} \varphi(n(x)g) dx = 0, {\rm a.e.} \  g\in G(\A). $$
	Denote by $R_{\omega}$, or simply $R$ if no confusion, the right regular representation of $G(\A)$ on $L^2(G(F) \backslash G(\A), \omega)$. Denote by $R_0$ its subrepresentation on $L_0^2(G(F) \backslash G(\A), \omega)$. We know that each irreducible component $\pi$ of $R$ decomposes into $\pi = \hat{\otimes}_v' \pi_v$, where $\pi_v$'s are irreducible unitary representations of $G_v$. $R = R_0 \oplus R_{\text{res}} \oplus R_{\text{c}}$ is the spectral decomposition. $R_0$ decomposes as a direct sum of irreducible $G(\A)$-representations, whose components are called cuspidal representations. $R_{\text{res}}$ is the sum of all one dimensional subrepresentations. $R_{\text{c}}$ is a direct integral of irreducible $G(\A)$-representations, expressed via Eisenstein series. Components of $R_0$ and $R_{\text{c}}$ are the generic automorphic representations. Recall that a principal series representation $\pi(\mu_1, \mu_2)=\text{Ind}_{B(F_v)}^{G(F_v)} (\mu_1, \mu_2)$ with $\mu_1 \mu_2^{-1} (t) = |t|_v^s, \forall t \in F_v $ is a complementary series if $s$ is a non-zero real number in the interval $(-1,1)$. $|s|/2$ is called its parameter. Let $\theta \in [0,1/2)$ be such that no complementary series representation with parameter $> \theta$ appears as a local component of a cuspidal representation.
	
	A compact open subgroup $K_{\text{f}}' \subset G(\A_{\text{f}})$ is said to be of (congruence) type 0, if for every finite place $v$, there is an integer $m_v$ such that the local component
	$$K_v' = K_v^0[m_v] := \left\{ \begin{pmatrix} a & b\\ c & d \end{pmatrix} \in G(\vO_v) \mid c \equiv 0 \mod \varpi_v^{m_v} \right\} ,$$
where $\varpi_v$ is a uniformiser of the local field $F_v$. Let $\varphi \in \pi$ be a pure tensor vector in an automorphic representation. Suppose for every $v<\infty$, $\varphi$ is invariant by $K_v^0[m_v]$ but not by $K_v^0[m_v-1]$, then we define $m_v = v(\varphi)$. Define $v(\pi) = v(\pi_v) = \min_{\varphi \in \pi_v} v(\varphi)$. The local conductor is $C(\pi_v) = \varpi_v^{v(\pi_v)}$. We similarly define the principal congruence subgroups 
	$$K_v[n] := \left\{ \begin{pmatrix} a & b\\ c & d \end{pmatrix} \in G(\vO_v) \mid a-1,b,c,d-1 \equiv 0 \mod \varpi_v^n \right\}.$$
We use the convention $K_v^0[0] = K_v[0] = K_v$.
	
	For any semisimple (real) Lie group $G$, denote by $\mathcal{C}_G$ the Casimir element. In our case, $G = {\rm GL}_2$. At each place $v \mid \infty$, $Z_v\backslash G_v$ is semisimple, and $\Delta_v = -\mathcal{C}_{Z_v\backslash G_v} - 2 \mathcal{C}_{K_v} $ is an elliptic operator on $Z_v\backslash G_v$. Note that here we calculate $\mathcal{C}_{K_v}$ by using the Killing form of ${\rm Lie}(Z_v\backslash G_v)$ instead of $K_v$'s Killing form.

	\subsection{$L$-function Theory for $K$-finite Vectors}

	The proof of the fact that the representation of $G(\A)$ on $L_0^2(G(F) \backslash G(\A), \omega)$ decomposes as a discrete direct sum of irreducible representations, as in Lemma 5.2 of \cite{G}, actually gives important information on $K$-finite vectors in an irreducible component $\pi$. They consequently have representatives in the space of smooth functions on the automorphic quotient, and are rapidly decreasing in any Siegel domain (Lemma 5.6 of \cite{G}). Let the superscript ``fin'' mean ``$K$-finite''. The rapid decay is important, because it adds to the description of $W_{\pi}^{{\rm fin}}$, the image of $\pi^{{\rm fin}} \subset \pi \subset L_0^2(G(F) \backslash G(\A), \omega)$ under the Whittaker intertwiner
	\begin{equation}
		\varphi \mapsto W_{\varphi}(g)=\int_{F \backslash \A} \varphi(n(x)g) \psi(-x) dx
	\label{WhittakerIntertwiner}
	\end{equation}
	the important growth property, which is essential for the uniqueness of Whittaker model at archimedean places (Section 2.8 and 4.4 of \cite{B} for local uniqueness, Section 3.5 of \cite{B} for global uniqueness). If $\varphi$ has a prescribed $K$-type and is a pure tensor, i.e. $W_{\varphi}(g) = \prod_v W_{\varphi,v}(g_v)$ splits, $W_{\varphi,v}(a(y)k)$ is forced to have rapid decay at $\infty$, thus has the nice behavior around $0$
	\begin{equation}
		|W_{\varphi,v}(a(y)k)| \ll |y|_v^{1/2-\theta-\epsilon}, \forall \epsilon > 0.
	\label{Whittakerat0}
	\end{equation}
	
	Now let $\chi$ be a character of $F^{\times} \backslash \A^{\times}$ and $s\in \C$. Jacquet-Langlands \cite{JL} defined a functional on $\pi^{{\rm fin}}$, called the (global) zeta-functional :
	$$ \zeta(s,\varphi,\chi) = \int_{F^{\times} \backslash \A^{\times}} \varphi(a(y)) \chi(y) |y|^{s-1/2} d^{\times}y, \forall \varphi \in \pi, a(y)=\begin{pmatrix} y & 0 \\ 0 & 1 \end{pmatrix}. $$
	Since $\varphi(a(y))$ is rapidly decreasing at $\infty$, it is also rapidly decreasing at $0$ because
	$$ \varphi(a(y)) = \varphi(wa(y)) =  \omega(y) \cdot w.\varphi(a(y^{-1})), w = \begin{pmatrix}  & -1 \\ 1 & \end{pmatrix}. $$
	Thus $\zeta(s,\varphi,\chi)$ is well defined for all $s$, and the following functional equation characterizes the left invariance by $w$ of $\varphi$:
	\begin{equation}
		\zeta(s,\varphi,\chi) = \zeta(1-s,w.\varphi,\omega^{-1}\chi^{-1}).
	\label{GlobalFE}
	\end{equation}
	If $\varphi$ is a pure tensor in $\pi^{{\rm fin}} \simeq \otimes_v' \pi_v^{{\rm fin}}$, i.e. $W_{\varphi}$ factorizes, then since
	\begin{equation}
		\varphi(g) = \sum_{t\in F^{\times}} W_{\varphi}(a(t)g),
	\label{WhittakerInversion}
	\end{equation}
	we get
	$$ \zeta(s,\varphi,\chi) = \prod_v \zeta(s,W_{\varphi,v},\chi_v,\psi_v), \Re (s) > 1+\theta $$
	with
	$$ \zeta(s,W_{\varphi,v},\chi_v,\psi_v) = \int_{F^{\times}_v} W_{\varphi,v}(a(y)) \chi(y) |y|^{s-1/2} d^{\times}y. $$
	The convergence is justified by the above local growth property of $W_{\varphi,v}$ and the fact that at an unramified finite place $v$, the local zeta-function equals
	$$ \zeta(s,W_{\varphi,v},\chi_v,\psi_v) = L(s, \pi_v \otimes \chi_v) = (1-\mu_v\chi_v(\varpi_v)q_v^{-s})^{-1}(1-\nu_v\chi_v(\varpi_v)q_v^{-s})^{-1}, $$
	where $\pi_v = \text{Ind}_{B(F_v)}^{G(F_v)} (\mu_v, \nu_v)$ determines $\mu_v, \nu_v$.
	
	The analysis of local zeta-functions shows that $\zeta(s,W_{\varphi,v},\chi_v,\psi_v)$, as $W_{\varphi,v}$ varies over $W_{\pi,v}^{{\rm fin}}$, has a ``common divisor'' $L(s,\pi_v \otimes \chi_v)$, which is a meromorphic function in $s$ such that $\frac{\zeta(s,W_{\varphi,v},\chi_v,\psi_v)}{L(s,\pi_v \otimes \chi_v)}$, originally defined for $\Re(s) > \theta $, can be analytically continued into an entire function on $s\in\C$. It equals $1$ at almost all places $v$. Furthermore, there is a functional equation
	\begin{equation}
		\frac{\zeta(s,W_{\varphi,v},\chi_v,\psi_v)}{L(s,\pi_v \otimes \chi_v)} \epsilon(s,\pi_v,\chi_v,\psi_v) = \frac{\zeta(1-s,wW_{\varphi,v},\omega_v^{-1}\chi_v^{-1},\psi_v)}{L(1-s,\pi_v \otimes \omega_v^{-1}\chi_v^{-1})}
	\label{LocalFE}
	\end{equation}
	where $\epsilon(s,\pi_v,\chi_v,\psi_v)$ is an entire function of exponential type. Define usual and complete $L$-functions as, for $\Re(s) > 1+\theta $,
	$$ L(s,\pi\otimes\chi) = \prod_{v < \infty} L(s,\pi_v \otimes \chi_v),$$
	$$ \Lambda(s,\pi\otimes\chi) = \prod_v L(s,\pi_v \otimes \chi_v),$$
	then the analytic continuations and functional equations of these $L$-functions follow from the well-definedness of $\zeta(s,\varphi,\chi)$ and (\ref{GlobalFE}), (\ref{LocalFE}). The identity
	$$ \zeta(s,\varphi,\chi) = L(s,\pi\otimes\chi) \prod_{v|\infty} \zeta(s,W_{\varphi,v},\chi_v,\psi_v) \prod_{v<\infty} \frac{\zeta(s,W_{\varphi,v},\chi_v,\psi_v)}{L(s,\pi_v \otimes \chi_v)} $$
	can be evaluated at $s = 1/2$ without analytic continuation of any integral. Thus
	\begin{equation}
		L(1/2,\pi\otimes\chi) = \prod_{v|\infty} \zeta(1/2,W_{\varphi,v},\chi_v,\psi_v)^{-1} \cdot \prod_{v<\infty} \frac{L(1/2,\pi_v \otimes \chi_v)}{\zeta(1/2,W_{\varphi,v},\chi_v,\psi_v)} \cdot \zeta(1/2,\varphi,\chi).
	\label{JL}
	\end{equation}
	\begin{remark}
		In fact, the above theory is valid for smooth (not necessarily $K$-finite) vectors as we shall explain in the following sections.
	\end{remark}

	\subsection{Smooth Vectors in Different Models}

	For any Lie group $G$ and a unitary representation $(\rho, V)$ of $G$, let $\rho^{\infty}$ be the subspace of smooth vectors in $V$. This is naturally a Fr\'echet space, defined by the semi-norms $\lVert X.v \rVert, X \in U(\mathfrak{g})$. If $V \subset L^2(M)$ is realized as a space of functions on a orientable real manifold $M$ equipped with a smooth (right) $G$-action, and with a $G$-invariant volume form, then we can talk about Sobolev functions for the action. Note that the action $\rho: G\to U(V)$ need not coincide with the regular representation on $L^2(M)$ induced by the action of $G$ on $M$. One may think about $\rho=\pi(\mu_1,\mu_2)$ in the principal unitary series of $G={\rm GL}_2(\R)$. 
	\begin{definition}
		With the above notations, a function $f$ on $M$ is called Sobolev (for the $G$-action), if it is smooth for the differential structure of $M$, and if its class $[f]$ in $V \subset L^2(M)$ is a smooth vector. We write $V^{\infty}$ or $\rho^{{\rm nam},\infty}$, if ${\rm nam}$ is the name of the model, or just $\rho^{\infty}$ if the underlying model is clear, for the space of Sobolev functions.
	\end{definition}
	We obviously have $[\rho^{{\rm nam},\infty}] \subset \rho^{\infty}$. Reciprocally,
	\begin{lemma}
		Assume that:\\
		1. For any $p \in M$, the map $s_p: G \to M, g \mapsto p.g $ is a submersion at the identity $e \in G$.\\
		2. The action of any element $X\in \mathfrak{g}$ on $V \cap C^{\infty}(M)$ corresponds to a smooth vector field $v(X)$ on M.\\
		Then every vector $v \in \rho^{\infty} \subset L^2(M)$ has a representative in $C^{\infty}(M)$.
	\end{lemma}
	\begin{definition}
		Fix a basis $\mathcal{B}$ of $\mathfrak{g}$, for any positive integer $d>0$, one can define a Sobolev norm on $\rho^{\infty}$ by
		$$ S^{\rho}_d(v) = \max_{X_i \in \mathcal{B}, l\leq d} \lVert X_1...X_l . v \rVert. $$
	\end{definition}
	\proof Since the condition and the conclusion are of local nature, one may interpret everything on the open set $C_p$ of some euclidean space, diffeomorphic to some open neighborhood $U_p$ of some point $p\in M$. The assumptions 1,2 ensures that the Sobolev norms $S_d^{\rho}$ are equivalent to the usual Sobolev norms on $C_p$ in the underlying euclidean space. One can apply the classical Sobolev embedding theorem. \endproof
	\begin{corollary}
		Under the assumptions of the above lemma, for any  $p \in M$, there is an integer $d$ such that $ \forall f \in L^2(M) \cap \rho^{\infty} $,
		$$ \sup_{q \in U_p}  |f(q)| \ll_{p,U_p} S^{\rho}_d([f]). $$
	\label{SobIneq}
	\end{corollary}
	The assumptions of the above lemma apply to the following situations:\\
	-- $\rho \subset R_{\omega}^{K_{\text{f}}'} $ is a $G(F_{\infty})$-subrepresentation of the right regular representation on automorphic quotient space. $M$ is therefore $G(F) \backslash G(\A)/K_{\text{f}}'$ where $K_{\text{f}}'$ is a compact open subgroup of $G(\A_{\text{f}})$. In such situation, we say that $\rho$ is realized in the automorphic model: ``aut''.\\
	-- $\rho = \pi(\chi_1, \chi_2)^{K_{\text{f}}'}$ is a principal unitary series representation with a compact open $K_{\text{f}}' \subset K_{\text{f}}$. $M$ is just $K / K_{\text{f}}'$. We say that $\rho$ is realized in the induced model: ``ind''.\\
	-- $\rho = W_{\pi}^{K_{\text{f}}'}$ is the Whittaker model of a generic automorphic representation $\pi$ with the same $K_{\text{f}}'$. $M$ is thus $AK/K_{\text{f}}'$. We say it is realized in the Whittaker model.\\
	-- $\rho = K_{\pi}^{K_{\text{f}}'}$ is the Kirillov model of a generic automorphic representation $\pi$ with the same $K_{\text{f}}'$. $M$ is thus $A/A \cap K_{\text{f}}'$. We say it is realized in the Kirillov model.
	\begin{definition}
		If $G$ is a totally disconnected group, acting on a totally disconnected space $M$, then a function $f$ on $M$ is said to be smooth, if it is locally constant on $M$ and $K$-finite for any maximal compact subgroup $K$ of $G$.
	\end{definition}

	\subsection{Smooth Vectors and Extended $L$-function Theory}

	We generalize the theory of $L$-function to smooth vectors. Using Corollary \ref{SobIneq} and compactness of $F\backslash \A $, one may easily see (Corollary I.1.5 \cite{CP-S}) that the Whittaker functional
	$$ l: R^{\infty} \to \C, \varphi \mapsto W_{\varphi}(1) $$
is in the continuous dual space of $R^{\infty}$ verifying
	$$l(R(n(x))\varphi ) = \psi(x) l(\varphi), $$
and is related to the Whittaker intertwiner (\ref{WhittakerIntertwiner}) by
	$$ W_{\varphi}(g) = l(R(g).\varphi). $$
	When we restrict to an irreducible component $\pi$ of $R$, or more precisely to $ \otimes_v' \pi_v^{\infty} \subset \pi^{\infty} $, it splits as
	$$ l = \otimes_v' l_v, $$
where $l_v$ is a local (continuous) Whittaker functional of $\pi_v^{\infty}$ verifying
	$$l_v(n(x) w) = \psi_v(x) l_v(w), w \in \pi_v^{\infty}. $$
	The study of $l_v, v< \infty$ is the same as in the $K_v$-finite case. So the uniqueness, the local functional equation (\ref{LocalFE}), the rapid decay and the controlled behavior at $0$ i.e. (\ref{Whittakerat0}) remain valid. At a $v | \infty$, the uniqueness of $l_v$ is established by Shalika \cite{S}. So one can define the smooth Whittaker model associated with a unitary irreducible representation $\pi_v$ by
	\begin{equation}
		W_{\pi_v}^{\infty} = \left\{ W_w(g) = l_v(\pi_v(g)w); w \in \pi_v^{\infty} \right\},
	\label{SmoothWhittaker}
	\end{equation}
as well as its smooth Kirillov model
	\begin{equation}
		K_{\pi_v}^{\infty} = \left\{ K_w(y) = W_w(a(y)); w \in \pi_v^{\infty} \right\}.
	\label{SmoothKirillov}
	\end{equation}
	The rapid decay at infinity of the local Whittaker functions $W_w(g)$ can be found in Lemma I.1.2 \cite{CP-S}. Note that here, the rapid decay property is derived from the continuity of $l_v$. In fact, much more information is obtained by Jacquet, as a special case in Proposition 3.6 \cite{C}, where the behavior of $W_w(g)$ is completely characterized, which implies rapid decay and (\ref{Whittakerat0}) in this situation. Consequently, the rapid decay of $\varphi \in \otimes_v' \pi_v^{\infty} \subset \pi^{\infty} \subset R_0^{\infty}$ follows by using (\ref{WhittakerInversion}). Furthermore, local functional equations (\ref{LocalFE}) are obtained by Jacquet \cite{J} with absolute convergence for $\Re (s) > \theta $ as in the $K_v$-finite case.
\begin{remark}
	For a proof that rapid decay at infinity and local functional equation imply the controlled behavior at $0$, see Proposition 3.2.3 \cite{MV}.
\label{RWhittakerat0}
\end{remark}

	\subsection{An Identification of Norms}

	A by-product of the above theory, already known in the $K$-finite case, is the identification of the norm on $\pi \subset R_0$ and the natural norm we put on global Whittaker models. We begin with the case of Eisenstein series for motivation.
	\begin{lemma}
		If $\pi = \pi(\chi_1,\chi_2)$ is unitary Eisenstein, and $\varphi(g) = E(0,f)(g)$ with $E(s,f)(g)$ defined as in (\ref{EisDef}), for some $f = \prod_v f_v \in\pi^{{\rm ind}, {\rm fin}}$ in the induced model, then one can define the Eisenstein norm of $\varphi$ by
		$$ \lVert \varphi \rVert_{{\rm Eis}}^2 = \int_{K} |f(k)|^2 dk. $$
		The following relation holds
		$$ \prod_{v\in V_F} \frac{\zeta_v(2)}{\zeta_v(1)^2 } \int_{F_v^{\times}}  |W_{\varphi,v}(a(y))|^2 d^{\times} y  = \lVert \varphi \rVert_{{\rm Eis}}^2, $$
		and the local data are defined as the analytic continuation in $(\chi_1, \chi_2)$ of
		$$ W_{\varphi,v}(g) = W_{f,v}(g) = \int_{F_v} f_v(wn(x) g) \psi_v(-x) dx. $$
	\label{NormIdenEis}
	\end{lemma}
	\proof One can interpret $W_{\varphi,v}(a(y)) \chi_{2,v}(y)^{-1} |y|^{-1/2}$ as the Fourier transform of $x \mapsto f(wn(x))$. The above norm identification is then a formal consequence of Plancherel formula as discussed in Section 3.1.6 of \cite{MV}. One can also verify it by using Theorem 4.6.5 of \cite{B}. \endproof
	\begin{remark}
		Note that $\frac{\zeta_v(2)}{\zeta_v(1)^2} = \frac{1-q_v^{-1}}{1+q_v^{-1}}$ is bounded both from above and below by some constants, uniformly in $v < \infty$.
	\label{LocalBoundRSLFactorEis}
	\end{remark}
	Let's turn to the cuspidal case.
	\begin{lemma}
		If $\pi = \hat{\otimes}_v' \pi_v \subset R_0$ and $\varphi \in \otimes_v' \pi_v^{\infty}$ is a pure tensor, then		
		$$ \lVert \varphi \rVert_{X(F)}^2 = 2 \Lambda(1, \pi, {\rm Ad}) \prod_{v\in V_F} \frac{ \zeta_v(2) \int_{F_v^{\times} \times K_v} |W_{\varphi,v}(a(y)k)|^2 d^{\times} y dk}{L(1,\pi_v \times \bar{\pi}_v)}, $$
where $\Lambda_F$ is the complete Dedekind zeta-function, and
		$$ \Lambda_F(s)\Lambda(s, \pi, {\rm Ad}) = \Lambda(s,\pi \times \bar{\pi}) = \prod_{v\in V_F} L(s,\pi_v \times \bar{\pi}_v) $$
is the complete $L$-function associated with $\pi \times \bar{\pi}$.
	\label{NormIden}
	\end{lemma}
	\begin{remark}
		By \cite{HL} for $F=\Q$ and Lemma 3 of \cite{BH4} in general, $C(\pi)^{-\epsilon} \ll L^*(1,\pi \times \bar{\pi}) \ll C(\pi)^{\epsilon}$. $L(s,\pi \times \bar{\pi}) = \prod_{v< \infty} L(s,\pi_v \times \bar{\pi}_v)$ is the incomplete Rankin-Selberg $L$-function and $L^*(1,\pi \times \bar{\pi})$ is its residue at $s=1$. Analogously to Remark \ref{LocalBoundRSLFactorEis}, we also note that there is some constant $C(\theta)$ depending only on $\theta$ such that
		$$ C(\theta)^{-1} \leq \left\lvert \frac{\zeta_v(2)}{L(1, \pi_v \times \bar{\pi}_v)} \right\rvert \leq C(\theta), $$
uniformly in $v < \infty$.
	\label{LocalBoundRSLFactorCusp}
	\end{remark}
	\proof It is a standard use of Rankin-Selberg's method (c.f. \cite{MV} 4.4.2) : Unfold, for $\Re{s} \gg 1$ 
	$$ \int_{ZG(F)\backslash G(\A)} \varphi(g) \bar{\varphi}(g) E(s,f)(g) d\bar{g} $$
	to get
	\begin{align*}
		\int_{\A^{\times} \times K} &|W_{\varphi}(a(y)k)|^2 f_s(a(y)k)|y|^{-1} d^{\times}y dk \\
		&=  \int_{\A^{\times} \times K} |W_{\varphi}(a(y)k)|^2 |y|^{s-1/2} d^{\times}y dk,
	\end{align*}
	where $f_s \in \pi(|\cdot|^{s}, |\cdot|^{-s})$ is a spherical flat section taking value $1$ on $K$, and
	\begin{equation}
		E(s,f)(g) = \sum_{\gamma \in B(F) \backslash G(F)} f_s(\gamma g).
	\label{EisDef}
	\end{equation}
	Then take residue at $s=1/2$.
	
	In fact, $E(s,f)$ converges absolutely for $\Re(s) > 1/2$, has a meromorphic continuation to all $s \in \C$ admitting a simple pole at $s=1/2$ with residue $\frac{1}{2} \frac{\Lambda_F^*(1)}{\Lambda_F(2)}$, and is of moderate growth for any given $s$ (see for example Section 3.7 of \cite{B}). Here $\Lambda_F^*(1)$ is the residue of $\Lambda_F(s)$ at $s=1$. At a place $v < \infty$, for which $W_{\varphi,v}$ is spherical, and if $\Re(s)>-1/2+2\theta$, one has
	\begin{equation}
		\frac{ \zeta_v(2s+1) \int_{F_v^{\times}\times K_v} |W_{\varphi,v}(a(y)k)|^2 |y|_v^{s-1/2} d^{\times} y dk}{L(s+1/2,\pi_v \times \bar{\pi}_v)} = |W_{\varphi,v}(1)|^2,
	\label{Unramified}
	\end{equation}
	which is $1$ for almost all $v$. The product $\prod_{v\in V_F} L(s+1/2,\pi_v \times \bar{\pi}_v)$ converges for $\Re(s)>1/2$. Thus for $\Re{s} > 1/2$,
	\begin{align*}
		\int_{ZG(F)\backslash G(\A)} &\varphi(g) \bar{\varphi}(g) E(s,f)(g) d\bar{g} \\
		&= \frac{\Lambda(s+1/2, \pi \times \overline{\pi})}{\Lambda_F(2s+1)} \prod_{v\in V_F} \frac{ \zeta_v(2s+1) \int_{F_v^{\times} \times K_v} |W_{\varphi,v}(a(y)k)|^2 |y|_v^{s-1/2}d^{\times} y dk}{L(s+1/2,\pi_v \times \bar{\pi}_v)}.
	\end{align*}
By the local behavior (\ref{Whittakerat0}), one can evaluate the integrals on the right side at $s=1/2$. Whence
	$$ \lVert \varphi \rVert_{X(F)}^2 \cdot \frac{1}{2} \frac{\Lambda_F^*(1)}{\Lambda_F(2)} = \frac{\Lambda_F^*(1) \Lambda(1, \pi, \text{Ad})}{\Lambda_F(2)} \prod_{v\in V_F} \frac{ \zeta_v(2) \int_{F_v^{\times} \times K_v} |W_{\varphi,v}(a(y)k)|^2 d^{\times} y dk}{L(1,\pi_v \times \bar{\pi}_v)}. $$
	\endproof
	We can simplify by taking into account the theory of Kirillov model.
	
	Define $B_1(F_v) = \left\{ \begin{pmatrix} a & b\\ 0 & 1 \end{pmatrix}: a \in F_v^{\times}, b \in F_v \right\}$.
	\begin{proposition}
		There are only two types of unitary irreducible representations of $B_1(F_v)$:\\
		\begin{itemize}
			\item[1.]	A character of $F_v^{\times} \simeq B_1(F_v) / N_v$;
			\item[2.]	The representation of $B_1(F_v)$ on $L^2(F_v^{\times})$ defined by the formula : $\begin{pmatrix} a & b\\ 0 & 1 \end{pmatrix} f (x) = \psi(bx) f(ax)$, where $\psi$ is a nontrivial character of $F_v$.
		\end{itemize}
		Moreover, for the second type, different $\psi$ give equivalent representations. In particular, there is only one non one-dimensional unitary irreducible representation of $B_1(F_v)$. 
	\label{UniqueRep}
	\end{proposition}
	A riguous proof of this proposition, in the case of an archimedean field, can be found in Page 34 (29), \cite{LL}; and in the case of a non archimedean filed, can be found in Chapter 8, \cite{BH2}.
	
	We finally deduce:
	\begin{proposition}
		Let $\pi$ be the local component on $v$ of a generic automorphic representation. For a $W \in W_{\pi}^{\infty}$, one actually has
		$$
			\int_{F_v^{\times} \times K_v} |W(a(y)k)|^2 d^{\times} y dk = \int_{F_v^{\times}} |W(a(y))|^2 d^{\times} y.
		$$
		As a consequence, the formula in Lemma \ref{NormIden} becomes
		$$ \lVert \varphi \rVert_{X(F)}^2 = 2 \Lambda(1,\pi,{\rm Ad}) \prod_{v\in V_F} \frac{ \zeta_v(2) \int_{F_v^{\times}} |W_{\varphi,v}(a(y))|^2 d^{\times} y }{L(1,\pi_v \times \bar{\pi}_v)}. $$
	\label{NormIdenKirillov}
	\end{proposition}
	\begin{remark}
		The norm identifications actually justify the notations $W_{\pi_v}^{\infty}$ and $K_{\pi_v}^{\infty}$ as smooth vectors in their completions $W_{\pi_v}$ and $K_{\pi_v}$.
	\end{remark}
	\begin{remark}
		The normalizations of local Whittaker functions $W_{\varphi,v}$ are different according that $\varphi$ is cuspidal or not. If it is cuspidal, the normalization is such that $W_{\varphi,v}(1) = 1$ for almost all $v$; while if it is (unitary) Eisenstein, the normalization imposes $f_v(1)=1$. This explains the missing of the factor $L(1,\pi_v \times \bar{\pi}_v)$ in the Eisenstein case.
	\end{remark}

	\subsection{Spectral Decomposition}

	The spectral decomposition, in the $L^2$ sense, is established in the first four sections of \cite{GJ}, which gives
	\begin{equation}
		R = \bigoplus_{\pi \text{ cuspidal}} \pi \oplus \int_{-i\infty}^{i\infty} \bigoplus_{\xi\in \widehat{F^{\times} \backslash \A^{(1)}}} \pi_{s,\xi} \frac{ds}{4\pi i} \oplus \bigoplus_{\chi \in \widehat{F^{\times} \backslash \A^{\times}},\chi^2 = \omega} \chi \circ \det,
	\label{RepsSpectralDecomp}
	\end{equation}
where $\pi_{s,\xi} = \pi(\xi |\cdot|^s, \omega\xi^{-1}|\cdot|^{-s})$. Note that $\pi_{s,\xi} \simeq \pi_{-s, \omega\xi^{-1}}$. According to Proposition I.1.4 of \cite{C}, the above spectral decomposition has an analogue for smooth vectors, namely
	\begin{equation}
		R^{\infty} = \bigoplus_{\pi \text{ cuspidal}} \pi^{\infty} \oplus \int_{-i\infty}^{i\infty} \bigoplus_{\xi\in \widehat{F^{\times} \backslash \A^{(1)}}} \pi_{s,\xi}^{\infty} \frac{ds}{4\pi i} \oplus \bigoplus_{\chi \in \widehat{F^{\times} \backslash \A^{\times}},\chi^2 = \omega} \chi \circ \det
	\label{SmoothSpectralDecomp}
	\end{equation}
with convergence for the topology of $R^{\infty}$. We are going to establish
	\begin{theorem}
		Suppose $\varphi \in R^{\infty}$, viewed as a function on $G(\A)$, then the following decomposition
		$$ \varphi(g) =  \sum_{\chi \in\widehat{F^{\times} \backslash \A^{\times}},\chi^2 = \omega} \frac{\langle \varphi, \chi \circ \det \rangle}{{\rm Vol}(X(F))} \chi \circ \det(g) + \sum_{\pi \ {\rm cuspidal}} \sum_{e \in \mathcal{B}(\pi)} \langle \varphi,e \rangle e(g) $$
		$$ + \sum_{\xi\in\widehat{F^{\times} \backslash \A^{(1)}}} \int_{-i\infty}^{i\infty} \sum_{\Phi \in \mathcal{B}(\pi_{0,\xi})} \langle \varphi, E(s,\Phi) \rangle E(s,\Phi)(g) \frac{ds}{4\pi i} $$
		converges absolutely and uniformly on any compact subset, where $\mathcal{B}(*)$ means taking an orthonormal basis of $*$ consisting of $K$-isotypical pure tensors. We may assume that if $\varphi$ is $K_v[n_v]$-invariant, then every function appearing on the right hand side is $K_v[n_v]$-invariant at any finite place $v$. $K_v$ need not be the standard maximal compact subgroup of $G_v$.
	\label{AnaSpDecompTrad}
	\end{theorem}
	\begin{remark}
		Therefore, the sum $\sum_{\xi\in\widehat{F^{\times} \backslash \A^{(1)}}}$ is actually finite and the number depends only on $F$ and $n_v$'s.
	\end{remark}
	If we consider the theory of Whittaker model as a theory of spectral decomposition with respect to the left action of $N(\A)$, then we further have
	\begin{theorem}
		Conditions are the same as in the above theorem. Any $\varphi \in R^{\infty}$, as a function on $G(\A)$, admits the following decomposition:
		$$ \varphi(g) = \varphi_N(g) + \sum_{\pi \ {\rm cuspidal}} \sum_{e \in \mathcal{B}(\pi)} \langle \varphi,e \rangle \sum_{\alpha \in F^{\times} }W_e(a(\alpha)g) + $$
		$$ \sum_{\xi\in\widehat{F^{\times} \backslash \A^{(1)}}} \int_{-i\infty}^{i\infty} \sum_{\Phi \in \mathcal{B}(\pi_{0,\xi})} \langle \varphi, E(s,\Phi) \rangle \sum_{\alpha \in F^{\times} } W_{\Phi,s}(a(\alpha)g) \frac{ds}{4\pi i}. $$
		The convergence is absolute and uniform on any Siegel domain.
	\label{AnaSpDecompWhittaker}
	\end{theorem}
	
	\begin{remark}
		In practice, the basis $\mathcal{B}(*)$ will be chosen so that the components of its elements at some archimedean place $v$ are $K_v$-isotypic where $K_v$ is the standard maximal compact subgroup of $G_v$.
	\label{ArchiBasis}
	\end{remark}
	We begin with some local Sobolev type analysis.
	
		\subsubsection{Local Bounds of $K$-isotypical Functions}
		
	\begin{lemma}
		Let $v$ be a finite place, and $\pi$ a unitary irreducible representation of $G_v$. Suppose $W \in W_{\pi}^{\infty}$, the smooth Whittaker model of $\pi$ w.r.t. $\psi_v$, is invariant by $K_v[m]$, then we have the following Sobolev inequality
		$$ |W(na(y)k)|^2 {\rm Vol}(1+\varpi_v^m\vO_v) \leq \lVert W \rVert^2 1_{v(y) \geq v(\psi)-m}, n\in N_v, y\in F_v^{\times}, k\in K_v $$
with the convention $1+\varpi_v^0\vO_v = \vO_v^{\times} $. At an unramified place $(m=0)$, we recall that
		$$ W(na(\varpi_v^l)k) = q_v^{-l/2} \frac{\alpha_1^{l+1}-\alpha_2^{l+1}}{\alpha_1 - \alpha_2} 1_{l\geq 0} W(1) $$
for some $\alpha_1, \alpha_2$ satisfying $|\alpha_1\alpha_2| = 1,  q_v^{-\theta} \leq |\alpha_1| \leq q_v^{\theta}$.
	\end{lemma}
	\proof Since $W$ is invariant by $K_v[m]$, for $x \in \varpi_v^m \vO_v$ we have
	$$ \psi_v(xy)W(y) = n(x)W(y) = W(y). $$
	But $\psi_v(xy)$ is not constantly $1$ for such $x$ if $v(y) < v(\psi) - m$, therefore $W(y) = 0$. We also have
	$$ W(yu) = a(u)W(y) = W(y), \forall u\in 1+\varpi_v^m \vO_v. $$
	We deduce
	$$ |W(y)|^2 {\rm Vol}(1+\varpi_v^m\vO_v) = \int_{t \in y(1+\varpi_v^m\vO_v)} |W(t)|^2 d^{\times}t \leq \lVert W \rVert^2, $$
and the Sobolev inequality follows by replacing $W$ by $k.W$, which is also $K_v[m]$-invariant, in the above argument. \endproof
	\begin{lemma}
		Let $v$ be a real place, and $\pi$ a unitary irreducible representation of $G_v$ with central character $\omega$. If $W\in W_{\pi}^{\infty}$, then
		$$ \forall n \in N(\R), y \in \R^{\times}, k \in {\rm SO}_2(\R), N \equiv 1 \pmod{2}, N > 0, $$
		$$ |W(na(y)k)| \ll_{N, \epsilon} |y|^{-N} \max(|y|^{\epsilon}, |y|^{-\epsilon}) S^{\pi}_{N+1}(W). $$
		Suppose further $W \in W_{\pi}^{{\rm fin}}$ transforms under the action of $K_v = {\rm SO}_2(\R)$ accroding to the character
		$$ \kappa_{\alpha} = \begin{pmatrix} \cos \alpha & \sin \alpha \\ -\sin \alpha & \cos \alpha \end{pmatrix} \mapsto e^{im\alpha}. $$
		Then we have the following Sobolev inequality, uniform in $m$,
		$$ |W(na(y)k)| \ll_{N, \omega, \theta, \epsilon} |y|^{-N} \max(|y|^{\epsilon}, |y|^{-\epsilon}) \lambda_W^{N'} \lVert W \rVert, $$
where $\lambda_W$ is the eigenvalue for $W$ of the elliptic operator $\Delta_v = -\mathcal{C}_{G_v} + 2 \mathcal{C}_{K_v}$, and $N'$ depends only on $N$ and $\omega$.
	\end{lemma}
	\proof Let $U = \begin{pmatrix} 1 &  \\   & 0 \end{pmatrix}$, $T = \begin{pmatrix} 0 & 1 \\  0 & 0 \end{pmatrix}$ be elements in the Lie algebra of ${\rm GL}_2(\R)$, then
	$$ T.W(a(y)) = -2\pi i y W(a(y)), U.W(a(y)) = y \frac{\partial}{\partial y} W(a(y)). $$
	We may only consider the case $y \in \R_+^{\times}$. Then $\forall x, y \in \R_+^{\times}$, we have
	$$ (-2\pi i y)^N W(a(y)) = T^N.W(a(x)) + \int_x^y UT^N. W(a(u)) d^{\times}u. $$
	Note that
	\begin{align*}
		\left\lvert \int_x^y UT^N. W(a(u)) d^{\times}u \right\rvert &\leq (\int_x^y |UT^N. W(a(u))|^2 d^{\times}u)^{1/2} (\int_x^y d^{\times}u)^{1/2} \\
		&\leq \lVert UT^N.W \rVert |\log (y/x)|^{1/2}.
	\end{align*}
	Thus
	$$ |(-2\pi i y)^N W(a(y))| \leq |T^N.W(a(x))| + \lVert UT^N.W \rVert |\log (y/x)|^{1/2}. $$
	Integrating against $\min(x,1/x) dx/x$ for $0<x<\infty$, using Cauchy-Schwarz and $\sqrt{1/2}(\sqrt{a} + \sqrt{b}) \leq \sqrt{a+b}$, we get
	$$ 2|(-2\pi i y)^N W(a(y))| \leq \lVert T^N.W \rVert + \lVert UT^N.W \rVert \int_0^{\infty} \min(x,1/x) |\log (y/x)|^{1/2} d^{\times}x. $$
	Using the bound $|\log t| \ll_{\epsilon} \max(t^{\epsilon},t^{-\epsilon})$, we get
	$$ |(-2\pi i y)^N W(a(y))| \ll_{\epsilon} \lVert T^N.W \rVert + \lVert UT^N.W \rVert \max(|y|^{\epsilon}, |y|^{-\epsilon}). $$
	Thus the first inequality follows for $k =1$. The general case follows by noting $S^{\pi}_{N+1}(k.W) \ll_N S^{\pi}_{N+1}(W)$, since the adjoint action of $K$ on $\mathfrak{g}$ has bounded coefficients.
	
	The second inequality follows from the equivalence of two systems of Sobolev norms. One is $S^{\pi}_d$'s, the other is defined with $\Delta_v$ and $I \in Z(\mathfrak{g})$. The proof is technical. We give it in the next section ( Theorem \ref{SobEquiv} ). \endproof
	
	Before proceeding to the complex place case, let's first recall that the irreducible representations of ${\rm SU}_2(\C)$ are parametrized by $m \in \mathbb{N}$, denoted by $(\rho_m,V_m)$. Here $V_m$ is the space of homogeneous polynomials in $\C[z_1,z_2]$ of degree $m+1$, equipped with the inner product
	$$ \langle P_1, P_2 \rangle = \int_{|z_1|^2 + |z_2|^2 \leq 1} P_1(z_1,z_2)\overline{P_2(z_1,z_2)} dz_1dz_2. $$
	The action of ${\rm SU}_2(\C)$ is given by
	$$ u.P(z_1,z_2) = P((z_1,z_2).u). $$
	Let $P_{m,k}(z_1,z_2)$ be a multiple of $z_1^{m-k}z_2^k$, normalized such that they form an orthonormal basis of $V_m$. Now let $\pi$ be a unitary irreducible representation of $G(\C)$. Let $W_{m,k} \in W_{\pi}^{{\rm fin}}$ span the $\rho_m$-isotypical subspace, with $W_{m,k}$ corresponding to $P_{m,k}$. Since $\rho_m$ is unitary, we have the following relation
	$$ \sum_{k=0}^m |W_{m,k}(gu)|^2 = \sum_{k=0}^m |W_{m,k}(g)|^2, \forall u\in {\rm SU}_2(\C). $$
	Therefore, we only need to bound $W_{m,k}(a(y))$ in order to bound $W_{m,k}(g)$. This works exactly the same as in the real place case. We omit the proof.
	\begin{lemma}
		Let $v$ be a complex place, and $\pi$ be a unitary irreducible representation of $G_v$ with central character $\omega$. If $W\in W_{\pi}^{\infty}$, then
		$$ \forall n \in N(\C), y \in \C^{\times}, k \in {\rm SU}_2(\C), N \in \mathbb{N}, $$
		$$ |W(na(y)k)| \ll_{N, \epsilon} |y|_v^{-N} \max(|y|_v^{\epsilon}, |y|_v^{-\epsilon}) S^{\pi}_{2N+2}(W). $$
		Suppose further $W \in W_{\pi}^{{\rm fin}}$ transforms under the action of $K_v = {\rm SU}_2(\C)$ accroding to $\rho_m$ and corresponds to some $P_{m,k}$. Then we have the following Sobolev inequality, uniformly in $m$,
		$$ |W(na(y)k)| \ll_{N,\omega, \theta } |y|_v^{-N} \max(|y|_v^{\epsilon}, |y|_v^{-\epsilon}) \lambda_W^{N'} \lVert W \rVert, $$
where $\lambda_W$ is the eigenvalue for $W$ of the elliptic operator $\Delta_v = -\mathcal{C}_{G_v} + 2 \mathcal{C}_{K_v}$, and $N'$ depends only on $N$ and $\omega$.
	\end{lemma}

		\subsubsection{Proof of Theorems \ref{AnaSpDecompTrad}, \ref{AnaSpDecompWhittaker}}
	
	We first deal with the cuspidal parts in the equations of Theorems \ref{AnaSpDecompTrad}, \ref{AnaSpDecompWhittaker}.
	
	Let $e \in \pi \subset R_0$ be a $K$-isotypic vector, with local Whittaker model $W_{e,v}$. Denote by $n_v$ the $K_v$-type of $W_{e,v}$, i.e.
	
	-- if $v < \infty$, then $W_{e,v}$ is $K_v[n_v]$-invariant. For almost all $v$, $n_v = 0$.
	
	-- if $v$ is a real place, then $W_{e,v}$ transforms under ${\rm SO}_2(\R)$ as $e^{in_v\alpha}$.
	
	-- if $v$ is a complex place, then $W_{e,v}$ transforms under ${\rm SU}_2(\C)$ as some $P_{n_v,k}$.
	
	Collecting all the estimations in the previous subsection, using Lemma \ref{NormIden} or Proposition \ref{NormIdenKirillov} with $\lVert e \rVert = 1$ and $C_{\infty}(\pi) \ll \lambda_{e,\infty} = \prod_{v| \infty} \lambda_{e,v}, C_f(\pi) \leq \prod_{v<\infty} q_v^{n_v}$, we obtain
	\begin{align*}
		W_e(na(y)k) &\ll_{F,N,\epsilon} |y|_{\infty}^{-N} \lambda_{e,\infty}^{N'} (\prod_{v<\infty} q_v^{n_v})^{\epsilon} \\
		&\  \cdot \prod_{v<\infty, n_v \neq 0} \left( L(1,\pi_v \times \bar{\pi}_v) \text{Vol}(1+\varpi_v^{n_v}\vO_v)^{-1} \right)^{1/2} \prod_{v<\infty} 1_{v(y) \geq v(\psi)-n_v} ,
	\end{align*}
where $|y|_{\infty} = \prod_{v|\infty} |y|_v$. The term $\prod_{v<\infty, n_v \neq 0} L(1,\pi_v \times \bar{\pi}_v) \text{Vol}(1+\varpi_v^{n_v}\vO_v)^{-1}$ can be bounded from above by a constant depending only on $n_v, v < \infty$, we thus get
	$$ W_e(na(y)k) \ll_{F,N,\epsilon,(n_v)_{v<\infty}} \lambda_{e,\infty}^{N'} |y|_{\infty}^{-N} \prod_{v<\infty} 1_{v(y) \geq v(\psi)-n_v}. $$
	Now since
	$$ e(na(y)k) = \sum_{\alpha\in F^{\times} } W_e(a(\alpha)na(y)k) = \sum_{\alpha\in F^{\times} } W_e(n'a(\alpha y)k), n' = a(\alpha)na(\alpha)^{-1}, $$
we have
	$$ \sum_{\alpha\in F^{\times} } |W_e(a(\alpha)na(y)k)| \ll_{F,N,\epsilon} C(n_v,v<\infty)\lambda_{e,\infty}^{N'} \sum_{\alpha\in F^{\times} } |\alpha y|_{\infty}^{-N} \prod_{v<\infty} 1_{v(\alpha y) \geq v(\psi)-n_v}. $$
Consider the splitting $\A^{\times} \simeq \A^1 \times \R_+$ and write $y=y_1t$ such that $y_1\in\A^1$ and $t \in \R_+ \hookrightarrow \A^{\times}$ with trivial component at finite places. We need only consider $y_1$ in a fundamental domain of $F^{\times} \backslash \A^1$. Since the quotient $F^{\times} \backslash \A^1$ is compact, we may assume that there exist $0<c<C$ such that for any place $v$, $c \leq |y_{1,v}|_v \leq C$ and for a.e. $v $, say $\forall v > v_0$, $|y_{1,v}|_v=1$. So the condition imposed in $\prod_{v<\infty}$ implies $|\alpha|_v \leq c^{-1}q_v^{n_v-v(\psi)}$ and $|\alpha|_v \leq 1, \forall v > v_0$ (one may choose $v_0$ big enough depending only on $n_v$'s) in order to get a non zero contribution. Thus, $\alpha$ runs over the non zero elements in a lattice of $F_{\infty}$ depending only on $n_v$'s. Therefore
	$$ \sum_{\alpha\in F^{\times} } |\alpha y|_{\infty}^{-N} \prod_{v<\infty} 1_{v(\alpha y) \geq v(\psi)-n_v} \ll_{n_v,v<\infty} |y|_{\infty}^{-N} \ll_{F,N} |y|^{-N}. $$
We conclude
	\begin{equation}
		\sum_{\alpha\in F^{\times} } |W_e(a(\alpha)na(y)k)| \ll_{F,N,n_v,v<\infty} \lambda_{e,\infty}^{N'} |y|^{-N}.
	\label{GlobalCuspWhittaker}
	\end{equation}
	
	Now let's turn to the Eisenstein parts of Theorems \ref{AnaSpDecompTrad}, \ref{AnaSpDecompWhittaker}.
	
	Using Lemma \ref{NormIdenEis} instead of \ref{NormIden} in the above argument, we get
	\begin{equation}
		\sum_{\alpha\in F^{\times} } |W_{\Phi,s}(a(\alpha)na(y)k)| \ll_{F,N,n_v,v<\infty} \lambda_{\Phi,s,\infty}^{N'} |y|^{-N}.
	\label{GlobalEisWhittaker}
	\end{equation}
We have an expression for the constant term
	$$ E(s,\Phi)_N(g) = \Phi_s(g) + M(s)\Phi_s(g). $$
	$\Phi_s |_K$ belongs to some irreducible component $\sigma$ of $\text{Res}_K^{G(\A)} \pi_{s,\xi} = \text{Ind}_{K\cap B(\A)}^K (\xi, \omega \xi^{-1})$. From basic representation theory, an orthonormal basis of functions on the compact group $K$ is given by matrix coefficients. So
	$$ \Phi_s(k) = \sqrt{\dim \sigma } <\sigma(k).v, v_0>_{\sigma}$$
with $v,v_0 \in \sigma$ of norm $1$, and
	$$ \sigma(b).v_0 = (\xi, \omega \xi^{-1})(b).v_0. $$
Thus follows the bound (recall that we are dealing with $\Re (s) = 0$)
	$$ |\Phi_s(na(y)k)| = |y|^{1/2} |\Phi_s(k)| \leq |y|^{1/2} \sqrt{\dim \sigma } \ll_{n_v,v<\infty } |y|^{ 1/2} \lambda_{K_{\infty}}(\Phi)^{1/2}, $$
where $\lambda_{K_{\infty}}(\Phi)$ is the eigenvalue of $\Phi$ for the Casimir of $K_{\infty}$. Note that $M(s)$ is unitary for $s \in i\R$ and doesn't change the $K$-type, thus
	$$ |M(s)\Phi_s(na(y)k)| \ll_{n_v,v<\infty } |y|^{ 1/2} \lambda_{K_{\infty}}(\Phi)^{1/2}. $$
Hence
	\begin{equation}
		| E(s,\Phi)_N(na(y)k) | \ll_{n_v,v<\infty } |y|^{ 1/2} \lambda_{K_{\infty}}(\Phi)^{1/2} \leq |y|^{ 1/2} \lambda_{\Phi_s,\infty}^{1/2}.
	\label{GlobalEisConstant}
	\end{equation}
Theorems \ref{AnaSpDecompTrad} \& \ref{AnaSpDecompWhittaker} will be established by using the following generalized Weyl's law, which is an immediate consequence of the Ph.D thesis \cite{Ma} of Marc R. Palm at G\"ottingen. 
	\begin{theorem}
		Given a sequence of non-negative integers $\bar{n} = (n_v)_{v < \infty}$ with $ n_v = 0 $ for a.e.$v$. Define
		$$ K_{{\rm f}}[\bar{n}] = \prod_{v<\infty} K_v[n_v] $$
and consider the space $R^{K_{{\rm f}}[\bar{n}]} = L^2(G(F) \backslash G(\A), \omega)^{K_{{\rm f}}[\bar{n}]}$. It is actually a representation of $G(F_{\infty}) \times K_{{\rm f}} $. The operator $\Delta_{\infty} = \prod_{v \mid \infty} \Delta_v$ is self-dual and commutes with the action of $K$. Then $\Delta_{\infty}^{-1-\epsilon}$ is of trace class in $R^{K_{{\rm f}}[\bar{n}]}$. More precisely,
		$$ \sum_{\pi'}\sum_{e} |\lambda_{e,\infty}|^{-1-\epsilon} + \sum_{\xi} \int_{-\infty}^{\infty} \sum_{\Phi} |\lambda_{\Phi_{i \tau}, \infty}|^{-1-\epsilon} \frac{d\tau}{4\pi} = O_{\epsilon}({\rm Vol}(Z(\A) G(F) \backslash G(\A) / K_{{\rm f}}[\bar{n}])). $$
		Here $\lambda_{e,\infty}$ runs over the discrete spectrum of $\Delta_{\infty}$, and $\lambda_{\Phi_{i \tau}, \infty}$ runs over the continuous spectrum of $\Delta_{\infty}$.
	\label{WeylLaw}
	\end{theorem}
	\begin{remark}
		We only need a weaker version here. Namely, we only need $\Delta_{\infty}^{-N}$ to be of trace class for some $N > 0 $.
	\end{remark}
	\begin{remark}
		If instead of $K_{{\rm f}}[\bar{n}]$ we consider $K_{\infty} \times K_{{\rm f}}[\bar{n}]$, the above theorem would coincide with the traditional geometrical Weyl's law. Note that this kind of Weyl's law was already used to establish theorems like \ref{AnaSpDecompTrad} for $K_{\infty}$-fixed case, e.g. \cite{CU}. Weyl's law is at the heart of the theory of analytical spectral decomposition.
	\end{remark}
	\begin{definition}
		(c.f. \cite{Ca}, Page 292) The Schwartz function space $R^s$ is the space of smooth functions $\varphi$ in ${\rm Ind}_{Z(\A)G(F)}^{G(\A)} \omega $, which are rapidly decreasing in any given Siegel domain, as well as $X.\varphi$ for any $X \in U(\mathfrak{g})$.	
	\end{definition}
The above argument also gives
	\begin{corollary}
		We have $R_0^{\infty} \subset R^s \subset R^{\infty}$.
	\label{Schwartz}
	\end{corollary}
	\begin{remark}
		If we take into account the central character, namely, if we write $R_{\omega}$ instead of $R$, we have $R^s_{\omega}R^s_{\omega'} \subset R^s_{\omega \omega'}$. In particular, if the central character is the trivial character $\omega_0$, $R^s_{\omega_0}$ is a ring for the pointwise multiplication.
	\end{remark}

	\subsection{Two Sobolev Norm Systems}

	Let $v$ be an archimedean place, and $\pi$ a unitary irreducible representation of $G_v$ with a fixed central character $\omega$. Let $\left\{ I_1,...,I_r \right\}$ be a basis of $Z(\mathfrak{g}_v)$. In our case, $r=1$ if $v$ is a real place, and $r=2$ if $v$ is a complex place. We define the Sobolev norm system
	$$H^{\pi}_d(v) = \max_{i_1+ \cdots +i_r + 2j = d } \lVert I_1^{i_1} \cdots I_r^{i_r} \Delta_v^j v \rVert . $$
	\begin{theorem}
		The Sobolev norm system $ H^{\pi}_d $ is equivalent to the Sobolev system $S^{\pi}_d$ for $\pi$ a local component of an automorphic representation. If the parameter $s$ of $\pi$ belongs to $ i\R \cup [-\theta, \theta]$ with $0 \leq \theta < 1/2$, then the implicit constants in the above equivalence can be taken independent of $\theta$.
	\label{SobEquiv}
	\end{theorem}
	The rest of this section is devoted to the proof of Theorem \ref{SobEquiv}.
	
		\subsubsection{$v$ a real place}

	The Hecke algebra $\mathcal{H}_v = U(\mathfrak{g})\oplus \underline{\epsilon} * U(\mathfrak{g})$, where $\underline{\epsilon}$ is the Dirac measure at $\begin{pmatrix} -1 & 0 \\ 0 & 1 \end{pmatrix}$. There is a classification of unitarizable irreducible $(\mathcal{H}_v, K_v)$-modules (c.f. for example 4.A \cite{G}). Each such module $\pi(\mu_1,\mu_2)$ is parametrized by $s_1,s_2 \in \C, m_1,m_2\in \left\{0,1\right\}$ with $\mu_i(t) = |t|^{s_i} {\rm sgn} (t)^{m_i}, i=1,2$. Put $s=s_1-s_2,t=s_1+s_2\in i\R, m=m_1-m_2$. There are three different series:
	\begin{itemize}
		\item[1.]	$s \in i\R$;
		\item[2.]	$0<s<1$ but only $s<2\theta$ is possible for the local component of an automorphic representation;
		\item[3.]	$0<s=p\in\Z,s-m$ is an odd integer.
	\end{itemize}
	In each series, there is an orthogonal, not necessarily normalized, basis consisting of $K_v$-isotypical vectors, $\left\{ e_k \right\}$. In the first two cases, $k$ runs through $k \equiv m \pmod{2}$, and in the last case, $|k|\geq p+1, k \equiv p+1 \pmod{2}$. There is a basis of $\mathfrak{g}_{\C}$,
	$$ \left\{ H=\begin{pmatrix} 0 & 1\\ -1 & 0 \end{pmatrix}, V_+=\begin{pmatrix} 1 & i\\ i & -1 \end{pmatrix}, V_-=\begin{pmatrix} 1 & -i\\ -i & -1 \end{pmatrix}, J=id \right\}$$
	with explicit action given as
	$$ H.e_k = ike_k; V_+.e_k = (s+1+k)e_{k+2}; V_-.e_k = (s+1-k)e_{k-2}; J.e_k = te_k; $$
	$$ \Delta_v.e_k = (\frac{1-s^2}{8}+\frac{k^2}{4})e_k. $$
Consider a general vector $v = \sum_k a_ke_k, a_k\in\C$. In the first series, Theorem 2.6.2 of \cite{B} implies $\lVert e_k \rVert = 1$. We easily deduce
	$$ \lVert H.v \rVert^2,\lVert V_+.v \rVert^2,\lVert V_-.v \rVert^2 \leq 16 \lVert \Delta_v^{1/2}.v \rVert^2. $$
In the second series, $\lVert e_k \rVert^2 = \left| \sqrt{\pi} \frac{\Gamma((s+1)/2) \Gamma(s/2)}{\Gamma((s+1+k)/2) \Gamma((s+1-k)/2)} \right|$ according to the proof of Theorem 2.6.4 of \cite{B}. As a consequence,
	$$ \frac{\lVert e_{k+2} \rVert^2}{\lVert e_k \rVert^2} = \left| \frac{s-1-k}{s+1+k} \right|, \frac{\lVert e_{k-2} \rVert^2}{\lVert e_k \rVert^2} = \left| \frac{s-1+k}{s+1-k} \right|. $$
	We get, for some implicit absolute constant,
	$$ \lVert H.v \rVert^2,\lVert V_+.v \rVert^2,\lVert V_-.v \rVert^2 \ll \lVert \Delta_v^{1/2}.v \rVert^2. $$
In the last series, it can be inferred from Theorem 2.6.5 of \cite{B} that $\pi(\mu_1,\mu_2)$ has the following model: Let $\mathbb{H}^+$ be the Poincar\'e half plane, and $\mathbb{H}^-$ its opposite. The space is, with the coordinates $z = x+iy$,
	$$ L^2(\mathbb{H}^{\pm}) = \left\{ f: \mathbb{H}^{\pm} \to \C \text{, holomorphic} : \int_{y\neq 0} |f(z)|^2 y^{p+1} \frac{dxdy}{|y|^2} < \infty \right\}. $$
Therefore one may take, for $|k| \geq p+1$,
	$$ e_k(z) = (z-i)^{-(k+p+1)/2}(z+i)^{(k-p-1)/2} 1_{\text{sgn}(k)\text{sgn}(y) < 0}. $$
Changing to the Poincar\'e disk model, one calculates easily, with $B(\cdot, \cdot)$ the Beta function,
	$$ \lVert e_k \rVert^2 = \pi 4^{-p} B((|k|-p-1)/2+1, p). $$
Consequently,
	$$ \frac{\lVert e_{k+2} \rVert^2}{\lVert e_k \rVert^2} \ll k^2,  \frac{\lVert e_{k-2} \rVert^2}{\lVert e_k \rVert^2} \ll k^2,$$
	$$ \lVert H.e_k \rVert^2,\lVert V_+.e_k \rVert^2,\lVert V_-.e_k \rVert^2 \ll \lVert \Delta_v^{1/2}.e_k \rVert^2. $$
We conclude that in all cases, by Cauchy-Shwarz and Weyl's law
	$$ \lVert H.v \rVert,\lVert V_+.v \rVert,\lVert V_-.v \rVert \ll \lVert \Delta_v.v \rVert + \lVert v \rVert, $$
thus $ S_d^{\pi} \ll_d H_d^{\pi} \ll S_{2d}^{\pi} $, and the two systems are equivalent.
	
		\subsubsection{$v$ a complex place}

	The unitary irreducible series $\pi(\mu_1, \mu_2)$ is parametrized by $s_1,s_2 \in \C, k_1,k_2 \in \Z$
	\begin{itemize}
		\item[-]	either with $t= s_1+s_2 \in i\R, s = s_1-s_2 \in i\R $ and $\mu_j(\rho e^{i \alpha}) = \rho^{2s_j} e^{i k_j \theta}, j=1,2$;
		\item[-]	or with $t= s_1+s_2 \in i\R, 0 < s = s_1-s_2 < 2\theta, k_1 = k_2 $.
	\end{itemize}
	Let $n_0 = k_1-k_2$. We may suppose $n_0 \geq 0 $ after exchanging $\mu_1$ and $\mu_2$ if necessary. The representation $\pi(\mu_1, \mu_2)$ has an orthogonal basis $\left\{ e_{n,k}^{(n_0)}: 0\leq k \leq n, n \geq |n_0|, n \equiv |n_0| \pmod{2} \right\}$ determined by
	$$ e_{n,k}^{(n_0)}(\begin{pmatrix} y_1 & x \\ 0 & y_2 \end{pmatrix} g ) = \mu_1(y_1)\mu_2(y_2) |y_1/y_2| e_{n,k}^{(n_0)}(g), \forall g \in G_v, $$
	$$ e_{n,k}^{(n_0)}(\begin{pmatrix} e^{i\alpha_1} & 0 \\ 0 & e^{-i\alpha_1} \end{pmatrix} u \begin{pmatrix} e^{i\alpha_2} & 0 \\ 0 & e^{-i\alpha_2} \end{pmatrix} ) = e^{in_0 \alpha_1} e^{i(n-2k)\alpha_2}, \forall u \in K_v = {\rm SU}_2(\C), $$
	$$ e_{n,k}^{(n_0)}(\begin{pmatrix} \cos \beta & \sin \beta \\ -\sin \beta & \cos \beta \end{pmatrix} ) = (\cos \beta)^{\frac{n+n_0}{2}-k} (\sin \beta)^{k-\frac{n-n_0}{2}} P_{\frac{n-n_0}{2}}^{(\frac{n_0-n}{2}+k, \frac{n_0+n}{2}-k)}(\cos 2\beta), $$
where $P_k^{(\alpha, \beta)}$ are the Jacobi polynomials. Alternatively,
	$$ e_{n,k}^{(n_0)} = \frac{\langle \rho_n(u)z_1^{n-k}z_2^k, z_1^{n-k_0}z_2^{k_0} \rangle_{\rho_n}}{\langle z_1^{n-k_0}z_2^{k_0},z_1^{n-k_0}z_2^{k_0} \rangle_{\rho_n}}, n-2k_0=n_0. $$
It will also be convenient to extend by $0$ to all integers $n,k$.
	
	The (complexified) Lie algebra $\mathfrak{su}_2$ has a basis
	$$ H_2 = \begin{pmatrix} i & 0 \\ 0 & -i \end{pmatrix}, X_{\pm} = \pm \begin{pmatrix} 0 & -1/2 \\ 1/2 & 0 \end{pmatrix} - i\begin{pmatrix} 0 & i/2 \\ i/2 & 0 \end{pmatrix}, $$
which act as
	$$ H_2.e_{n,k}^{(n_0)} = i(n-2k)e_{n,k}^{(n_0)}, X_+.e_{n,k}^{(n_0)} = (n-k)e_{n,k+1}^{(n_0)}, X_-.e_{n,k}^{(n_0)}=ke_{n,k-1}^{(n_0)}, $$
	$$ \Delta_v.e_{n,k}^{(n_0)} = ((1-s^2-n_0^2)/8+n(n+2)/4)e_{n,k}^{(n_0)}. $$
It is then obvious that $\Delta_v^{-1-\epsilon}$ is of trace class in $\pi(\mu_1, \mu_2)$. A standard argument then shows that it suffices to prove Theorem \ref{SobEquiv} for vectors of an orthonormal basis. The Cartan complement $\mathfrak{p}$ of $\mathfrak{su}_2$ has a basis (we ignore the center)
	$$ H_1 = \begin{pmatrix} 1 & 0 \\ 0 & -1 \end{pmatrix}, Y_+ = ad(X_+)(H_1), Y_-= ad(X_-)(H_1). $$
	Using the recurrence relations of Jacobi polynomials (c.f. \cite{A}),
	$$ (1-x^2)\frac{d}{dx}P_n^{(\alpha,\beta)}(x) = \frac{n[\alpha-\beta-(2n+\alpha+\beta)x]}{2n+\alpha+\beta}P_n^{(\alpha,\beta)}(x) + \frac{2(n+\alpha)(n+\beta)}{2n+\alpha+\beta}P_{n-1}^{(\alpha,\beta)}(x), $$
	\begin{align*}
		xP_n^{(\alpha,\beta)}(x) &= \frac{2(n+1)(n+1+\alpha+\beta)}{(2n+1+\alpha+\beta)(2n+2+\alpha+\beta)}P_{n+1}^{(\alpha,\beta)}(x) + \\
		& \frac{\beta^2-\alpha^2}{(2n+2+\alpha+\beta)(2n+\alpha+\beta)} P_n^{(\alpha,\beta)}(x) + \\
		& \frac{2(n+\alpha)(n+\beta)}{(2n+1+\alpha+\beta)(2n+\alpha+\beta)} P_{n-1}^{(\alpha,\beta)}(x) 
	\end{align*}
we can find for $n > 0$
	\begin{align*}
		H_1. e_{n,k}^{(n_0)} &= \frac{(s+n/2+1)(n-n_0+2)(n+n_0+2)}{(n+1)(n+2)} e_{n+2,k+1}^{(n_0)} \\
		&\ + \frac{2sn_0(n-2k)}{n(n+2)} e_{n,k}^{(n_0)} + \frac{(s-n/2)4k(n-k)}{n(n+1)} e_{n-2,k-1}^{(n_0)};
	\end{align*}
while for $n=0$
	$$ H_1.e_{0,0}^{(0)} = 2(s+1) e_{2,1}^{(0)}. $$
Since
	$$ Y_+.e_{n,k}^{(n_0)} = X_+H_1.e_{n,k}^{(n_0)} - H_1X_+.e_{n,k}^{(n_0)}, Y_-.e_{n,k}^{(n_0)} = X_-H_1.e_{n,k}^{(n_0)} - H_1X_-.e_{n,k}^{(n_0)}, $$
	we can only consider the actions of $H_1,H_2,X_+,X_-$ if we don't want to optimize.
	\begin{itemize}
		\item[Case 1:]	$s \in i\R$. Then we are in the unitary principal series case and the norm structure is the standard $L^2$-norm on ${\rm SU}_2(\C)$,
	$$ \lVert e_{n,k}^{(n_0)} \rVert^2 = \frac{(n-k)! k!}{(\frac{n-n_0}{2})!(\frac{n+n_0}{2})!(n+1)}. $$
	One easily verifies, if $n \neq 0$, $ \lVert X.e_{n,k}^{(n_0)} \rVert \ll \lVert \Delta_v^{1/2}. e_{n,k}^{(n_0)} \rVert, X=H_1,H_2,X_+,X_- $, hence
	$$ \lVert X.v \rVert \ll \lVert \Delta_v^2.v \rVert + \lVert \Delta_v.v \rVert + \lVert v \rVert, \forall v \in \pi^{\infty}, X=H_1,H_2,X_{\pm},Y_{\pm}. $$
	
		\item[Case 2:]	$0 < s < 2\theta < 1$. Then $n_0 = 0$, thus $n \equiv 0 \pmod{2}$. Let's write $e^{(s,0)}_{n,k} = e^{(0)}_{n,k}$ to emphasize the dependence on $s$. The norm satisfies
	$$ \lVert e^{(s,0)}_{n,k} \rVert^2 = (-1)^{n/2} \pi \frac{(s-1)\cdots (s-n/2)}{s(s+1)\cdots (s+n/2)} \frac{(n-k)! k!}{(\frac{n}{2})!(\frac{n}{2})!(n+1)}, $$
	which will be given by Lemma \ref{ComplexIsotypicalNorm}. With this, we easily see
	$$ \lVert X.v \rVert \ll \lVert \Delta_v^2.v \rVert + \lVert \Delta_v.v \rVert + \lVert v \rVert, \forall v \in \pi^{\infty}, X=H_1,H_2,X_{\pm},Y_{\pm}. $$
	\end{itemize}
	In the last case, the norm structure is defined via the intertwining operator (with analytic continuation for $s < 0$),
	$$ M(s)e^{(s,0)}_{n,k}(g) = \int_{ \C} e^{(s,0)}_{n,k}(n(x)g) dx = \lambda_{n,k}(s) e^{(-s,0)}_{n,k}(g). $$
	\begin{lemma}
		We have
		$$ \lambda_{n,k}(s) = (-1)^{n/2} \pi \frac{(s-1)\cdots (s-n/2)}{s(s+1)\cdots (s+n/2)}, $$
		Therefore,
		$$ \lVert e^{(s,0)}_{n,k} \rVert^2 = (-1)^{n/2} \pi \frac{(s-1)\cdots (s-n/2)}{s(s+1)\cdots (s+n/2)} \frac{(n-k)! k!}{(\frac{n}{2})!(\frac{n}{2})!(n+1)}. $$
	\label{ComplexIsotypicalNorm}
	\end{lemma}
	\proof We first consider $n=2k$. We know $e^{(s,0)}_{2k,k}(\begin{pmatrix} 1 & 0 \\ 0 & 1 \end{pmatrix}) = P_k^{(0,0)}(1) = 1$, so
	$$ \lambda_{2k,k}(s) = M(s)e^{(s,0)}_{n,k}(\begin{pmatrix} 1 & 0 \\ 0 & 1 \end{pmatrix}) 
	= \frac{\pi}{2} \int_{-1}^1 (\frac{1-t}{2})^{s-1} P_k^{(0,0)}(t)  dt. $$
	Now we can use the recurrence relation of Legendre polynomials
	to establish 
	$$ \lambda_{2k+2,k+1}(s) = \frac{2(2k+1)}{s} \lambda_{2k,k}(s+1) + \lambda_{2k-2,k-1}(s). $$
	The first two values are easy to obtain:
	$$ \lambda_{0,0}(s) = \frac{\pi}{s}, \lambda_{2,1}(s) = - \frac{\pi(s-1)}{s(s+1)}. $$
	By induction, we get
	$$ \lambda_{2k,k}(s) = (-1)^k \pi \frac{(s-1)\cdots (s-k)}{s(s+1)\cdots (s+k)}. $$
	Since $M(s)$ commutes with the action of $G_v$, it commutes with the action of $X_+,X_-$. It follows that for any $k$, $\lambda_{n,k}(s) = \lambda_{n,n/2}(s)$. This proves the above lemma and concludes the proof of Theorem \ref{SobEquiv}. \endproof

	\subsection{Construction of Automorphic Forms from Local Kirillov Models}
	\label{ConstAutoForm}

	The norm identifications tell us that, given a pure tensor $\varphi \in \otimes_v' \pi_v^{\infty}$, resulting from (\ref{WhittakerIntertwiner}), the $W_{\varphi,v}$ or the $K_{\varphi,v} $ must be a smooth vector in $W_{\pi_v}$ or $K_{\pi_v}$. Conversely, if we are given $K_v \in K_{\pi_v}^{\infty}$, which uniquely determine corresponding $W_v \in W_{\pi_v}^{\infty}$, and form $W(g) = \prod_v W_v(g_v)$, and $\varphi $ by (\ref{WhittakerInversion}), are we sure to get an element in $\pi^{\infty}$? The converse theorem, as is discussed in Section 5.2 of \cite{C}, gives an affirmative answer. Note that, to determine $W_v$ from $K_v$ at an archimedean place $v$, a concrete way is to apply the Casimir element $\mathcal{C} $ of ${\rm GL}_2(\R)$ in the real case, or the two embedded Casimir elements of ${\rm GL}_2(\R)$ in ${\rm GL}_2(\C)$ to get partial differential equations, since these elements should act as scalars depending only on $\pi_v$, then solve the corresponding Dirichlet problems.
	
	Alternatively, maybe also more naturally and directly, if one wants to avoid the converse theorem, one may decompose $W$ as an infinite sum of $K$-isotypical Whittaker functions, then change the order of summation to show that $\varphi$ is a convergent (thanks to the local and global estimations in the above sections) infinite sum of $K$-isotypical functions in $\pi$, with rapidly decreasing spectral parameter for $K$, thus is itself in $\pi^{\infty}$.

	\subsection{Decay of Matrix Coefficients: General Theory}
	\label{DefOfXi}

	At a place $v$, let $\pi_{\lambda}$ be the complementary series representation of $G_v$ with parameter $\lambda /2$ and with trivial central character. It has a unique $K_v$ invariant unit vector $w^0$. The elementary spherical function associated with $\pi_{\lambda}$ is defined to be (following Harish-Chandra's notation)
	$$ \varphi_{v, \lambda}(g) = \langle \pi_{\lambda}(g) w^0, w^0 \rangle. $$
	Its limit when $\lambda \to 0$, denoted by $\varphi_{v,0} = \Xi_v$, is the Harish-Chandra function. They are all positive and bi-$K_v$-invariant.
	\begin{theorem}
		Let $\pi$ be any unitary irreducible representation of $G_v$. Let $x_1,x_2$ be two $K_v$-finite vectors in $\pi$. Then
		\begin{itemize}
		\item[1] If $\pi$ is tempered, then
			$$ \langle \pi(g)x_1,x_2 \rangle \leq \dim (K_v x_1)^{1/2} \dim (K_v x_2)^{1/2} \Vert x_1 \Vert \cdot \Vert x_2 \Vert \Xi_v(g). $$
		\item[2] If $\pi$ is in the complementary series with parameter $\lambda /2$, then for any $\epsilon > 0$, there is a $A_v(\epsilon) > 0$
			$$ \langle \pi(g)x_1,x_2 \rangle \leq A_v(\epsilon) \dim (K_v x_1)^{1/2} \dim (K_v x_2)^{1/2} \Vert x_1 \Vert \cdot \Vert x_2 \Vert \Xi_v(g)^{1-\lambda-\epsilon}. $$
		\end{itemize}	
		Here $\dim (K_v x) = \dim {\rm span}(K_v \cdot x)$ is the dimension of the span of $x$ by $K_v$-action.
	\label{MatrixCoeffDecay}
	\end{theorem}
	\proof The tempered case is well known in \cite{CHH}. The non-tempered case, first proved in Theorem 2.11 \cite{S2} for real case, then recaptured in Lemma 9.1 \cite{V}, essentially is based on the following estimation
	\begin{equation}
		A_v(\epsilon)^{-1} \varphi_{v, 0}^{1-\lambda+\epsilon} \leq \varphi_{v, \lambda} \leq \varphi_{v,0}^{1-\lambda}.
	\label{SphericalEstimation}
	\end{equation}
	In fact, we have (c.f. 5.2 \cite{CU})
	$$ \varphi_{v, \lambda}(g) = f_v(\lambda, g) / f_v(1, g). $$
	\begin{itemize}
		\item[1.] If $v$ is complex, then
		$$ f_v(\lambda, \begin{pmatrix} T & 0 \\ 0 & T^{-1} \end{pmatrix} ) = \frac{T^{\lambda}-T^{-\lambda}}{\lambda} =  \log T \int_{-1}^1 T^{u\lambda} du, \forall T \geq 1. $$
		\item[2.] If $v$ is real, then
		$$ f_v(\lambda, \begin{pmatrix} e^{r/2} & 0 \\ 0 & e^{-r/2} \end{pmatrix} ) = \int_0^{2\pi} (\cosh r + \sinh r \cos u)^{(\lambda-1)/2} du, \forall r \geq 0 .$$
		\item[3.] If $v$ is finite, let $\varpi$ be a uniformizer, and $q$ the cardinality of the residue field, then
		$$ f_v(\lambda, m) = f_v(\lambda, \begin{pmatrix} \varpi^m & 0 \\ 0 & 1 \end{pmatrix} ) = q^{\lambda m/2} + q^{-\lambda m /2} + (1-q^{-1}) \sum_{k=1}^{m-1} q^{(k-m/2)\lambda}, \forall m \in \mathbb{N}. $$
	\end{itemize}
	The upper bound of (\ref{SphericalEstimation}) follows from the convexity of $\log f_v(\lambda, g)$ in $\lambda$. The lower bound follows by taking the major term in $f_v(\lambda,g)$. For example, in the case of a finite place, we use $f_v(\lambda, m) \geq q^{\lambda m/2}$ to get
	$$ \frac{\varphi_{v, \lambda}(m)}{(\varphi_{v, 0}(m))^{1-\lambda+\epsilon}} = \frac{f_v(\lambda,m)}{f_v(0,m)^{1-\lambda+\epsilon}f_v(1,m)^{\lambda-\epsilon}} \geq q^{\epsilon m/2} (\frac{1+q^{-1}}{m+1})^{1-\lambda+\epsilon} (1+q^{-1})^{-1}. $$
	Thus one may take $A_v(\epsilon)^{-1} = (\frac{\epsilon \log q}{2(1+\epsilon)})^{1+\epsilon} (1+q^{-1})^{\epsilon} q^{\frac{1+\epsilon}{\log q}-\frac{\epsilon}{2}}$ to conclude Theorem \ref{MatrixCoeffDecay}. \endproof

\section{Outline of the Proof}

The departure point of the proof is Jacquet-Langlands' generalization of Hecke's integral representation of $L$-functions, namely equation (\ref{JL}) that we copy here
$$ L(1/2,\pi\otimes\chi) = \left[ \prod_{v|\infty} \zeta(1/2,W_{\varphi,v},\chi_v,\psi_v)^{-1} \cdot \prod_{v<\infty} \frac{L(1/2,\pi_v \otimes \chi_v)}{\zeta(1/2,W_{\varphi,v},\chi_v,\psi_v)} \right] \cdot \zeta(1/2,\varphi,\chi). $$
Here $\varphi \in \pi^{\infty}$ is a pure tensor and smooth vector. We are going to establish the following proposition, which obviously implies Theorem \ref{Main}:

\begin{proposition}
	There is a pure tensor $\varphi \in \otimes_v' \pi_v^{\infty}$ such that for any $\epsilon > 0$,
	\begin{equation}
		\prod_{v|\infty} \zeta(1/2,W_{\varphi,v},\chi_v,\psi_v)^{-1} \cdot \prod_{v<\infty} \frac{L(1/2,\pi_v \otimes \chi_v)}{\zeta(1/2,W_{\varphi,v},\chi_v,\psi_v)} \ll_{\epsilon,F} Q^{1/2+\epsilon},
	\label{LocalGoal}
	\end{equation}
	where $Q=C(\chi)$ is the analytic conductor of $\chi$.
	
	There is an absolute constant $\delta > 0$ such that for any $\epsilon > 0$
	\begin{equation}
		\zeta(1/2,\varphi,\chi) \ll_{\epsilon,\pi} Q^{-\delta +\epsilon}.
	\label{GlobalGoal}
	\end{equation}
	We may choose $\delta = \frac {1 - 2 \theta}{8} $, or $\frac{25}{256}$ using the best known result of \cite{BB} i.e. $\theta = \frac{7}{64}$.
\label{Est}
\end{proposition}

The construction of $\varphi$ has its origin in an idea of P.Sarnak \cite{Sa} in the archimedean aspect. We consider the following family of test vectors of the form
$$ \varphi = n(t).\varphi_0, $$
where $\varphi_0 \in \pi^{\infty}$ is a fixed pure tensor and $t \in \A$. With this choice, the study of local zeta-functions shows that, under some technical conditions on $\varphi_0$, each local integral reaches its natural asymptotic lower bound for some $t_v = T_v$ with $|T_v|_v \asymp_{\epsilon} C(\chi_v)^{1 \pm \epsilon}$. We may see later that we can take $T_v=0$ for almost all $v$. We take $\varphi = n(T).\varphi_0$ with $T=(T_v)_v$ chosen above, then we get the estimation of the product of local terms in (\ref{JL}).

Recall the global zeta-function defined by
$$  \zeta(1/2, \varphi, \chi) = \int_{\A^{\times}} \left( \varphi - \varphi_N \right) (a(y)) \chi(y) d^{\times} y, a(y) = \begin{pmatrix} y & 0 \\ 0 & 1 \end{pmatrix}, $$
where the constant term $\varphi_N = 0$ since $\pi$ is cuspidal. We want to bound the global zeta-function by some negative power of $C(\chi)$. To deal with the fact that $F^{\times} \backslash \A^{\times}$ is non-compact, we then truncate the integral $\int_{F^{\times} \backslash \A^{\times}}$ as $\int_{F^{\times} \backslash \A^{\times}}^* := \int_{F^{\times} \backslash \A^{\times}} h(|y|)$, where $h: \R_+ \to [0,1]$ is a smooth function with compact support which will be described explicitly later. We remark that
$$ \left\lvert \int_{F^{\times} \backslash \A^{\times}}^* \varphi(a(y))\chi(y) d^{\times}y \right\rvert \leq \left( \int_{F^{\times} \backslash \A^{\times}}^* 1 d^{\times}y \right)^{1/2} \left( \int_{F^{\times} \backslash \A^{\times}}^* \left( n(T).|\varphi_0|^2 \right)(a(y)) d^{\times}y \right)^{1/2}. $$
The translation $n(t)$ on $|\varphi_0|^2$ is the same as translating the domain of integration $a(F^{\times} \backslash \A^{\times})$ into $a(F^{\times} \backslash \A^{\times})n(t)$. In the classical case ($F = \Q$ and $\varphi_0$ is spherical), the translated domain is the same as the semi straight-line $\{ yt + yi: y > 0 \}$. As $t \to \infty$, the slope of the line tends to $0$. The line becomes equidistributed on the modular surface ${\rm SL}_2(\Z) \backslash \mathbb{H}$. As a consequence the $n(t)$ (or $n(T)$) translation ``kills'' the portion of $|\varphi_0|^2$ orthogonal to the $1$-dimensional representations. Intuitively,
\begin{equation}
	\int_{a(F^{\times} \backslash \A^{\times})n(T)}^* |\varphi_0(t)|^2 dt \to \int_{Z(\A)G(F) \backslash G(\A)} |\varphi_0(g)|^2 dg = \langle \varphi_0, \varphi_0 \rangle .
\label{EffetEquidistribution}
\end{equation}
In order to diminish the right hand side, we amplify $\varphi_0$ by defining, for $E$ equal to some positive power of $Q$ to be chosen later, the following average of Dirac measures :
$$ \sigma = \frac{1}{M_E^2} \sum_{v,v'\in I_E} \delta_{|\varpi_v|_v |\varpi_{v'}^{-1}|_{v'}}  $$
with
$$ I_E = \left\{ v \mid q_v \in [E,2E], T_v = 0, \pi_v \text{ is unramified} \right\}, M_E = |I_E| \gg E/\log E, $$
and take, with $\varpi_v$ denoting a uniformiser at the place $v$,
$$ \varphi_0' = \frac{1}{M_E^2} \sum_{v,v'\in I_E} \chi(\varpi_v \varpi_{v'}^{-1}) a (\varpi_v \varpi_{v'}^{-1}).\varphi_0 = \sigma_{\chi}' * \varphi_0, $$
where $\sigma_{\chi}' = \frac{1}{M_E^2} \sum_{v,v'\in I_E} \chi(\varpi_v \varpi_{v'}^{-1}) \delta_{a (\varpi_v \varpi_{v'}^{-1})}$ is the adjoint measure of $\sigma$, i.e.
$$ \int_{F^{\times} \backslash \A^{\times}} h(|y|) n(T).\varphi_0'(a(y))\chi(y) d^{\times}y = \int_{F^{\times} \backslash \A^{\times}} \left(\sigma * h\right) (|y|) n(T).\varphi_0(a(y))\chi(y) d^{\times}y. $$
Here we have used the fact that the translations $a(\varpi_v \varpi_{v'}^{-1})$ commute with the translation $n(T)$. Instead of $\varphi_0$, we put $\varphi_0'$ into the above argument. This modification does not change the quality of truncation on integral. But in (\ref{EffetEquidistribution}), we get $\langle \varphi_0', \varphi_0' \rangle$ on the right hand side instead, which is some weighted average of
\begin{equation}
	\langle a(\frac{\varpi_{v_1}}{\varpi_{v_1'}}) \varphi_0, a(\frac{\varpi_{v_2}}{\varpi_{v_2'}}) \varphi_0 \rangle, v_1,v_1',v_2,v_2' \in I_E.
\label{Amplification1}
\end{equation}
Since the decay of matrix coefficients is of local nature, (\ref{Amplification1}) must be of size some negative power of $E$ when $v_1,v_1',v_2,v_2'$ are distinct. When $v_1,v_1',v_2,v_2'$ are not distinct, (\ref{Amplification1}) is bounded by $O(1)$, and the total contribution of this case is killed by the big denominator $M_E^4$. Of course this modification will increase the contribution of non one-dimensional parts of $|\varphi_0|^2$ by some positive power of $E$ as a factor.

Finally, we optimize the choice of $E$ and the truncation on integral to get (\ref{GlobalGoal}).

Let's discuss (\ref{GlobalGoal}) in more detail. In order to simplify notations and for further convenience, we introduce a functional on automorphic representations:
	$$ \varphi \mapsto l^{\chi |\cdot|^s}(\varphi) = \int_{F^{\times} \backslash \A^{\times}} \varphi(a(y)) \chi(y) |y|^s d^{\times}y = \zeta(s+1/2,\varphi,\chi), $$
	so that (\ref{GlobalGoal}) is equivalent to $$l^{\chi}(\varphi) \ll_{\epsilon,\pi} Q^{-\delta +\epsilon}. $$
There is a local analogue of this functional :
	$$
		W_{\varphi, v} \mapsto l^{\chi_v |\cdot|^s}(W_{\varphi,v}) = \int_{F_v^{\times}} W_{\varphi,v}(a(y)) \chi(y) |y|^s d^{\times}y.
	$$
 The truncation function $h \in C_c^{\infty}(\R_+)$ is made from a fixed function $h_0$ such that $h$ is supported in $[Q^{-\kappa-1},Q^{\kappa-1}]$. Here, $\kappa \in (0,1)$ is a parameter to be chosen later.
	\begin{lemma}
		We have
		$$
			l^{\chi}(\varphi) = \int_{F^{\times} \backslash \A^{\times}} \sigma * h(|y|) \varphi(a(y)) \chi(y) d^{\times}y + O_{h_0, \varphi_0, \epsilon}(Q^{-\kappa / 2+\epsilon}).
		$$
	\label{Tronc}
	\end{lemma}
Define another functional:
	$$
		\varphi \mapsto l^{\chi,h}(\varphi) = \int_{F^{\times} \backslash \A^{\times}} h(|y|) \varphi(a(y)) \chi(y) d^{\times}y.
	$$
We are reduced to examining :
	$$
		l^{\chi, \sigma * h}(\varphi) = l^{\chi,h}(\sigma_{\chi}' * \varphi) = \int_{F^{\times} \backslash \A^{\times}} h(|y|) \sigma_{\chi}' * \varphi(a(y)) \chi(y) d^{\times}y.
	$$
The inequality of Cauchy-Schwarz gives
	\begin{eqnarray}
		|l^{\chi,h}(\sigma_{\chi}' * \varphi)|^2 &\leq& \int_{F^{\times} \backslash \A^{\times}} h(|y|) d^{\times}y \cdot \int_{F^{\times} \backslash \A^{\times}} |\sigma_{\chi}' * \varphi(a(y))|^2 h(|y|) d^{\times}y,	\label{CS} \\
		&=& \int_{F^{\times} \backslash \A^{\times}} h(|y|) d^{\times}y \cdot l^h(n(T)|\sigma_{\chi}' * \varphi_0|^2), \nonumber
	\end{eqnarray}
where the second equality is due to the fact that $a(\varpi_v \varpi_{v'}^{-1})$ commute with $n(T)$. We then spectrally decompose $|\sigma_{\chi}' * \varphi_0|^2$ in $L^2(G(F) \backslash G(\A),1)$ as in Theorem \ref{AnaSpDecompWhittaker}, which is possible because $\varphi_0 \in R^s$. Setting $l^h = l^{1,h}$, we can interchange integrals as
	$$ l^h(n(T)|\sigma_{\chi}' * \varphi_0|^2) = \Sigma_1 + \Sigma_2 + \Sigma_3, $$
	where
	\begin{align*}
		\Sigma_1 &= l^h(n(T)|\sigma_{\chi}' * \varphi_0|^2_N ), \\
		\Sigma_2 &= \sum_{\pi' \ \text{cuspidal}} l^h(n(T) P_{\pi'}(|\sigma_{\chi}' * \varphi_0|^2)), \\
		\Sigma_3 &= \frac{1}{4\pi}\sum_{\xi \in \widehat{F^{\times} \backslash \A^{(1)}}} \int_{-\infty}^{\infty} l^h(n(T) (P_{\xi, i\tau}(|\sigma_{\chi}' * \varphi_0|^2) - P_{\xi, i\tau}(|\sigma_{\chi}' * \varphi_0|^2)_N) ) d\tau.
	\end{align*}
	This is verified by Theorem \ref{AnaSpDecompWhittaker}. In every summand of $\Sigma_2$ (resp. $\Sigma_3$), $P_{\pi'}$(resp. $P_{\xi, i\tau}$) denotes the projector onto the space of $\pi'$(resp. $\pi(\xi |\cdot|^{i\tau}, \xi^{-1}|\cdot|^{-i\tau})$). The function
		$$
			|\sigma_{\chi}' * \varphi_0|^2 = \frac{1}{M_E^4} \sum_{v_1,v_1',v_2,v_2'\in I_E}  \chi(\frac{\varpi_{v_1}}{\varpi_{v_1'}} ) \chi^{-1}(\frac{\varpi_{v_2}}{\varpi_{v_2'}} )
			a(\frac{\varpi_{v_1}}{\varpi_{v_1'}}) \varphi_0 a(\frac{\varpi_{v_2}}{\varpi_{v_2'}}) \overline{\varphi}_0.
		$$
	Let's write
	$$ S_{\text{cusp}}(v_1,v_1',v_2,v_2') =   \sum_{\pi'\  \text{cuspidal}} l^h( n(T) P_{\pi'}(a(\frac{\varpi_{v_1}}{\varpi_{v_1'}}) \varphi_0 a(\frac{\varpi_{v_2}}{\varpi_{v_2'}}) \overline{\varphi}_0 )), $$
	hence
	$$ \Sigma_2 =  \frac{1}{M_E^4} \sum_{v_1,v_1',v_2,v_2'\in I_E}  \chi(\frac{\varpi_{v_1}}{\varpi_{v_1'}} ) \chi^{-1}(\frac{\varpi_{v_2}}{\varpi_{v_2'}} ) S_{\text{cusp}}(v_1,v_1',v_2,v_2'). $$
	Define
	$$ S_{\text{cst}}(v_1,v_1',v_2,v_2') = l^h(n(T) (a(\frac{\varpi_{v_1}}{\varpi_{v_1'}}) \varphi_0 a(\frac{\varpi_{v_2}}{\varpi_{v_2'}}) \overline{\varphi}_0 )_N ) =  l^h( (a(\frac{\varpi_{v_1}}{\varpi_{v_1'}}) \varphi_0 a(\frac{\varpi_{v_2}}{\varpi_{v_2'}}) \overline{\varphi}_0 )_N ), $$
	$$ S_{\text{Eis}}(v_1,v_1',v_2,v_2') = \sum_{\xi \in \widehat{F^{\times} \backslash \A^{(1)}}} \int_{-\infty}^{\infty} l^h(n(T) P_{\xi, i\tau}(a(\frac{\varpi_{v_1}}{\varpi_{v_1'}}) \varphi_0 a(\frac{\varpi_{v_2}}{\varpi_{v_2'}}) \overline{\varphi}_0) - P_{\xi, i\tau}(a(\frac{\varpi_{v_1}}{\varpi_{v_1'}}) \varphi_0 a(\frac{\varpi_{v_2}}{\varpi_{v_2'}}) \overline{\varphi}_0)_N) d\tau. $$
	Therefore,
	\begin{equation}
		\Sigma_1 =  \frac{1}{M_E^4} \sum_{v_1,v_1',v_2,v_2'\in I_E}  \chi(\frac{\varpi_{v_1}}{\varpi_{v_1'}} ) \chi^{-1}(\frac{\varpi_{v_2}}{\varpi_{v_2'}} ) S_{\text{cst}}(v_1,v_1',v_2,v_2'),
	\label{AverageCst}
	\end{equation}
	$$ \Sigma_3 = \frac{1}{4\pi M_E^4} \sum_{v_1,v_1',v_2,v_2'\in I_E}  \chi(\frac{\varpi_{v_1}}{\varpi_{v_1'}} ) \chi^{-1}(\frac{\varpi_{v_2}}{\varpi_{v_2'}} )  S_{\text{Eis}}(v_1,v_1',v_2,v_2'). $$
	\begin{remark}
		Not every cuspidal representation $\pi'$ (resp. not every character $\xi$) has a non-trivial contribution in this decomposition. Only the ones which have ``smaller'' conductors than $\sigma_{\chi}' * \varphi_0$ at every place $v$ do. The exact choice of the base for spectral decomposition is a subtle matter. It will be described in Section \ref{BaseChoice}. Similarly, the number of $\xi$'s with non-zero contribution is also finite and depends on $F$ and $\varphi_0$.
	\end{remark}
	\begin{lemma}
		We have
		$$ \Sigma_1 \ll_{\epsilon, F, \pi } \kappa E^{\epsilon - 2} Q^{(2+\kappa) \epsilon }.  $$
	\label{ConstantL}
	\end{lemma}	
	\begin{proposition}
		A full list of the patterns of positions of $v_1,v_1',v_2,v_2'$ is
		\begin{itemize}
		\item	\textbf{Type 1:} $v_1,v_1',v_2,v_2'$ are distinct.
		\item	\textbf{Type 2:} $v_1 = v_2$ or $v_1' = v_2'$, and there are $3$ elements in $\left\{ v_1,v_1',v_2,v_2' \right\}$.
		\item	\textbf{Type 3:} $v_1 = v_2'$ or $v_1' = v_2 $ and there are $3$ elements in $\left\{ v_1,v_1',v_2,v_2' \right\}$.
		\item	\textbf{Type 4:} $v_1 = v_1'$ or $v_2 = v_2' $ and there are $3$ elements in $\left\{ v_1,v_1',v_2,v_2' \right\}$.
		\item	\textbf{Type 5:} $v_1 = v_2$ and $v_1' = v_2'$, and there are $2$ elements in $\left\{ v_1,v_1',v_2,v_2' \right\}$.
		\item	\textbf{Type 6:} $v_1 = v_2'$ and $v_1' = v_2 $ and there are $2$ elements in $\left\{ v_1,v_1',v_2,v_2' \right\}$.
		\item	\textbf{Type 7:} $v_1 = v_1'$ and $v_2' = v_2 $ and there are $2$ elements in $\left\{ v_1,v_1',v_2,v_2' \right\}$.
		\item	\textbf{Type 8:} $v_1 = v_1' = v_2$ or $v_1 = v_1' = v_2'$ or $v_2 = v_2' = v_1$ or $v_2 = v_2' = v_1'$ and there are $2$ elements in $\left\{ v_1,v_1',v_2,v_2' \right\}$.
		\item	\textbf{Type 9:} $v_1 = v_1' = v_2 = v_2'$.
		\end{itemize}
	\label{ATranslationTypes}
	\end{proposition}
\textbf{Type 1} is dominant in the sense that there are $\simeq M_E^4$ possibilities for this case but $O(M_E^3)$ for the other cases. Therefore it is also called to be {\it typical}. 

Recall that, $\theta$ is such that no complementary series representation with parameter $> \theta$ appears as a local component of a cuspidal representation. Let $\lambda_{e,\infty}$ (resp. $\lambda_{\Phi_{i \tau}, \infty}$) be the eigenvalue for $e$ (resp. $E(i\tau, \Phi)$) of $\Delta_{\infty}$, for $e$ (resp. $\Phi$) runing through an orthonormal base $\mathcal{B}(\pi')$ (resp. $\mathcal{B}(\pi (\xi, \xi^{-1})$), consisting of pure tensors of $\pi'$ (resp. $\pi(\xi, \xi^{-1})$). For the portion $\Sigma_2+\Sigma_3$, we need an adelic version of Weyl's law Theorem \ref{WeylLaw} and Lindel\"of's hypothesis on average. From it we deduce
	\begin{lemma}
		For a typical term, we have
		$$ S_{\rm cusp}(v_1,v_1',v_2,v_2') \ll_{\epsilon, F, \pi, \theta, \kappa, h_0} E^{2+\epsilon} Q^{-(1/2-\theta)+\epsilon }. $$
		In general, we have
		$$
			\Sigma_2 \ll_{\epsilon, F, \pi, \theta, \kappa, h_0} E^{2+\epsilon} Q^{-(1/2-\theta)+\epsilon }.
		$$
	\label{CuspL}
	\end{lemma}
	\begin{lemma}
		For a typical term, we have
		$$
			S_{{\rm Eis}}(v_1,v_1',v_2,v_2') \ll_{\epsilon, F, \pi, \kappa, h_0} E^{1+\epsilon} Q^{(\kappa-1)/2+\epsilon}.
		$$
		In general, we have
		$$ \Sigma_3 \ll_{\epsilon, F, \pi, \kappa, h_0} E^{1+\epsilon} Q^{(\kappa-1)/2+\epsilon}. $$
	\label{EisL}
	\end{lemma}
Lemmas \ref{ConstantL} to \ref{EisL} immediately imply
	\begin{lemma}
		We have
		$$
			l^h(n(T)|\sigma_{\chi}' * \varphi_0|^2) \ll_{\pi,\kappa, \epsilon} E^{\epsilon - 2} Q^{(2+\kappa )\epsilon } + E^{2+\epsilon} Q^{-(1/2-\theta)+\epsilon } + E^{1+\epsilon} Q^{(\kappa-1)/2+\epsilon}.
		$$
	\label{SpecFinal}
	\end{lemma}
	\begin{remark}
		A comparison between the eigenvalues appearing here and those appearing in the trace of $\Delta_{\infty}^{-l}$ should be taken into account, where $l > 1$ will be specified. We'll see this in detail later.
	\label{CompareTrace}
	\end{remark}
	\begin{remark}
		We should consider all types in Proposition \ref{ATranslationTypes} and recollect their effects to get the second assertions in Lemmas \ref{CuspL} and \ref{EisL}. But it turns out that the contribution of \textbf{Type 1} is always no less than that of other types.
	\end{remark}
Now it is clear that (\ref{GlobalGoal}) follows from Lemma \ref{Tronc}, (\ref{CS}) and Lemma \ref{SpecFinal}, by solving the equation 
$$
	\min_{\kappa, E} \max (E^{-1}, EQ^{-1/4+\theta /2}, Q^{-\kappa/2}, E^{1/2}Q^{(\kappa-1)/4} ) = Q^{-\frac{1-2\theta}{8}}.
$$ 
An optimal choice is $$E = Q^{\frac{1-2\theta}{8}}, \kappa = \frac{1}{4}+\frac{\theta}{6}. $$

\begin{remark}
	If we apply the $n(T)$ translation before the projections in $\Sigma_2$ and $\Sigma_3$, and use a more general result concerning the decay of matrix coefficients, then we find ourselves in the exact setting of \cite{MV}, where all the technical calculations are folded in the ``Ergodic Principle'' in Section 2.5.3.
\label{DifferenceFromMV}
\end{remark}

\section{Choice of $\varphi_0$ and Local Estimation}

	In this section we define the vector $\varphi $ of Proposition \ref{Est}. Recall that it is of the shape $\varphi = n(T)\varphi_0 $. Here $\varphi_0 \in \pi$ is a pure tensor corresponding to $W_0(g) = \prod_v W_{0,v}(g_v)$ in the Kirillov model of $\pi$. Recall also that we only need to specify $W_{0,v}$ for every place $v \in V_F$.

	\subsection{Archimedean places}
	We first make the notion ``Analytic Conductor'' precise. The general definition, for both ${\rm GL}_1$ and ${\rm GL}_2$ representations, is given in 3.1.8 \cite{MV}. In this paper, we're particularly interested in ${\rm GL}_1$ case. Using the notations from 3.1.8 \cite{MV} and from Chapter \Rmnum{14} $\mathcal{x}$ 4 \cite{L}, one easily sees that if $F_v = \R$ and $\chi_v(a) = \text{sgn} (a)^m |a|^{i\varphi}$, then $\mu_{\chi_v} = \frac{i\varphi+m}{2}, m \in \left\{ 0,1 \right\}$, and we may define
	$$C(\chi_v) = 2 + |\frac{i\varphi+m}{2}|.$$
	If $F_v = \C$ and $\chi_v(a) = (\frac{a}{|a|})^m |a|^{i2\varphi}$, then $\mu_{\chi_v} = i\varphi + |m|/2$, and we may define
	$$C(\chi_v) = (2 + |i\varphi + |m|/2|)^2.$$
	\begin{lemma}
		Let $\phi \in S(F_v^{\times})$ (i.e. $\phi$ as well as all its derivatives decay faster than any polynomial of $|t^{-1}|$ as $|t| \to +\infty$ and more rapidly than any polynomial of $|t|$ as $|t| \to 0$). Let $C = C(\chi_v)$ be the analytic conductor of $\chi_v$. Set, for $t \in F_v^{\times}$, the generalized Gauss sum
		$$
			G_{\phi}(\chi_v, t) = \int_{F_v } \phi(x) \psi_v(tx) \chi_v(x) dx.
		$$
		Then for any $N\in \mathbb{N}, 1/2 \leq \alpha < \beta <1$,
		$$
			|G_{\phi}(\chi_v,t)| \ll_{\phi, N, \alpha, \beta} \min \left( \left(\frac{1+|t|}{C}\right)^N, \left(\frac{C}{|t|}\right)^N, C^{1/2-\alpha}|t|^{\alpha - \beta} \right).
		$$
	\end{lemma}
	This is essentially the Lemma 3.1.14 of \cite{MV}. Let's recall the proof:
	\proof Note that $C$ is comparable with the maximal absolute value among eigenvalues of $\chi_v$ for a fixed $F_v^{\times}$-invariant basis of differential operators of degree $[F_v : \R]$. The first two bounds then follow from two different kinds of integration by parts. For the third bound, applying the local functional equation as in Tate's thesis, we obtain
	$$
		G_{\phi}(\chi_v,t) = \frac{\int_{F_v} \Phi(x+t) \chi_v^{-1}(x)|x|^{\alpha}d^{\times}x}{\gamma(\chi_v, \psi_v,1-\alpha)},
	$$
	where $\Phi = \widehat{\phi |\cdot|^{\alpha}} \in S(F_v)$ is the Fourier transform of $\phi(x) |x|^{\alpha}$. Recall if we fix a small $\epsilon > 0$, and let $\alpha \in [1/2, 1-\epsilon]$, by (3.5) of \cite{MV}, and the third property after Theorem 3 of \cite{L}, $\mathcal{x}3$
	$$
		|\gamma(\chi_v, \psi_v,1-\alpha)| \simeq_{\epsilon} C^{\alpha-1/2}.
	$$
	Then after some evident change of variables, one gets
	$$
		|G_{\phi}(\chi_v,t)| \simeq_{\epsilon} C^{1/2-\alpha} |t|^{\alpha} \left\lvert \int_{F_v} \Phi(tx) |x-1|^{\alpha-1} \chi^{-1}(x-1) dx \right\rvert.
	$$
	But for any $\beta > 0$, $ \Phi(x) \ll_{\alpha, \beta, \phi} |x|^{-\beta} $, thus
	$$
		|G_{\phi}(\chi_v,t)| \ll_{\epsilon, \alpha, \beta, \phi} C^{1/2-\alpha} |t|^{\alpha-\beta} \int_{F_v} |x|^{-\beta} |x-1|^{\alpha-1} dx.
	$$
	The integral converges if $1/2 \leq \alpha < \beta < 1$. Under this condition, we get
	$$ |G_{\phi}(\chi_v,t)| \ll_{\alpha, \beta, \phi} C^{1/2-\alpha} |t|^{\alpha-\beta}. $$
	\endproof
	
	\begin{corollary}
		For any $\epsilon > 0$, there is a constant $C_0$ depending only on $\phi$ and $\epsilon$, such that for $C \geq C_0$, there exists $t$ with $|t| \in [C^{1-\epsilon}, C^{1+\epsilon}]$, and $|G_{\phi}(\chi, t)| \gg_{\phi, \epsilon} C^{-1/2-\epsilon}$.
	\label{ArchGauss}
	\end{corollary}
	
	\proof Apply the Plancherel formula for $L^2(F_v)$
	\begin{align*}
		\int_{F_v} |\phi(x)|^2 dx &= \int_{F_v} |G_{\phi}(\chi_v,t)|^2 dt \ll_{\phi, N } \int_{|t| \leq C^{1-\epsilon}} (\frac{1+|t|}{C})^{2N} dt + \int_{|t| \geq C^{1+\epsilon}} (\frac{C}{|t|})^{2N} \\
		&\  + (C^{1+\epsilon} - C^{1-\epsilon}) \max_{|t| \in [C^{1-\epsilon}, C^{1+\epsilon}]} |G_{\phi}(\chi_v,t)|^2.
	\end{align*}
	The result follows by taking $N = 1+\lceil \frac{1}{2 \epsilon} \rceil $ ($N > 1/2+\frac{1}{2 \epsilon}$ suffices) for example. \endproof
	
	We choose $W_{0,v} \in S(F_v^{\times})$ and $T_v = t$ as in the above corollary, such that
	\begin{equation}
		\zeta(1/2,n(T_v)W_{0,v},\chi_v,\psi_v) \gg_{\epsilon, W_{0,v}} C(\chi_v)^{-1/2-\epsilon}.
	\label{LocalEst1}
	\end{equation}
	
	\begin{corollary}
		For any $0 < \epsilon < 1/2$, and any $\sigma \in \R$ varying in a compact set, we have
		$$ |G_{\phi}(\chi_v |\cdot|_v^{\sigma},t)| \ll_{\epsilon, \phi} \min(C^{-1/2 + \epsilon}, |t|^{-1/2 + \epsilon}).	$$
	\label{ArchiUpperEst}
	\end{corollary}	
	\proof In the case $\sigma = 0$, we have $ |G_{\phi}(\chi_v,t)| \ll_{\alpha, \beta, N, \phi} \min (C^{1/2-\alpha} |t|^{\alpha-\beta}, |t|^N C^{-N}) \leq C^{-\frac{N(\beta-1/2)}{N+\beta-\alpha}} $. Taking $\alpha = 1/2$, $\beta$ approaching $1$ and $N$ big enough gives the result. The general case follows by considering $G_{\phi}(\chi_v |\cdot|_v^{\sigma},t) = G_{\phi |\cdot|_v^{\sigma}}(\chi_v,t)$. \endproof
	
	\begin{remark}
		Note that $(\chi_v, t) \mapsto G_{\phi}(\chi_v,t)$ is a continuous function on $\widehat{F_v^{\times}} \times F_v$, which for each fixed $\chi_v$ is analytic in the variable $t$ and is not identically $0$. Hence it doesn't vanish identically for $t \in [1,2]$. Since $C(\chi_v) \leq C_0$ defines a compact region for $\chi_v$, a routine argument of compactness gives that, for any $\chi_v $ such that $C(\chi_v) \leq C_0$ there is a $t$ such that $|t| \in [1,2]$ and $G_{\phi}(\chi_v,t) \gg_{C_0,\phi} 1 $, hence $G_{\phi}(\chi_v,t) \gg_{\phi, \epsilon} C(\chi_v)^{-1/2-\epsilon} $. Note that $|t| \in [1,2]$ is included in $C(\chi_v)^{1-\epsilon} / C_0 \leq |t| \leq 2C(\chi_v)^{1+\epsilon} $. We obtain in all cases the existence of some $C^{1-\epsilon} \ll_{\phi, \epsilon} |t| \ll C^{1+\epsilon}$ such that $G_{\phi}(\chi_v,t) \gg_{\phi, \epsilon} C(\chi_v)^{-1/2-\epsilon} $, which suffices for our application.
	\end{remark}

	\subsection{Non-Archimedean places}
	We study the analog of the generalized Gauss sum as in the previous subsection at a finite place. Let's first recall some basic properties of Gauss sums.
	\begin{definition}
		Let $\chi$ be a character of $\vO^{\times}$ with $v(\chi)=r>0$, and $\psi$ be an additive character of $F_v$ with $v(\psi) = l$. The Gauss sum associated with $\chi$ and $\psi$ is defined by
		$$ G(\chi, \psi) = \int_{\vO_v^{\times}} \chi(t) \psi(t) d^{\times}t. $$
	\end{definition}
	Let's write $U_v^n = 1+\varpi_v^n \vO$ for $n \in \N$ with the convention $U_v^0 = \vO_v^{\times}$.
	\begin{proposition}
		The Gauss sum $G(\chi, \psi)$ is non zero only when $l=r$, in which case its absolute value is given by
		$$ |G(\chi, \psi)| = q_v^{-r/2} \zeta_v(1), $$
		where the factor $\zeta_v(1)$ is due to the measure normalization $\int_{\vO^{\times}} d^{\times}t = 1$.
	\label{ClassicalGauss}
	\end{proposition}
	\proof If $l<r$, then for $\forall x \in \vO_v$, we have
	\begin{eqnarray*}
		G(\chi, \psi) &=& \int_{\vO_v^{\times}} \psi(t(1+\varpi_v^{r-1}x)) \chi(t(1+\varpi_v^{r-1}x)) d^{\times}t \\
		&=& \chi(1+\varpi_v^{r-1}x) \int_{\vO_v^{\times}} \psi(t) \chi(t) d^{\times}t = \chi(1+\varpi_v^{r-1}x) G(\chi, \psi),
	\end{eqnarray*}
	where we have used $\psi(\varpi_v^{r-1}tx) = 1$ since $\varpi_v^{r-1}tx \in \varpi_v^l \vO_v$. But $\chi(1+\varpi_v^{r-1}x)$ is not identically $1$, hence $G(\chi, \psi) = 0$.
	
	If $l>r$, then we have
	\begin{align*}
		G(\chi,\psi) &= \sum_{a \in \vO_v^{\times} / U_v^r} \chi(a) \int_{U_v^r} \psi(at) d^{\times} t \\
		&= q_v^{d_v/2} \sum_{a \in \vO_v^{\times} / U_v^r} \chi(a) q_v^{-r} \int_{\vO_v} \psi(a(1+\varpi_v^r x)) dx \\
		&= q_v^{d_v/2} \sum_{a \in \vO_v^{\times} / U_v^r} \chi(a) \psi(a) q_v^{-r} \int_{\vO_v} \psi(\varpi_v^r x) dx.
	\end{align*}
	But $x \mapsto \psi(\varpi_v^r x)$ is a non-trivial additive character of $\vO_v$, hence
	$$ \int_{\vO_v} \psi(\varpi_v^r x) = 0, $$
	and we deduce that $G(\chi,\psi) = 0$.
	
	We finally assume $l \geq r$ and calculate
	\begin{align*}
		|G(\chi,\psi)|^2 &= \int_{\vO_v^{\times}} \int_{\vO_v^{\times}} \psi(t_1 - t_2) \chi(t_1 t_2^{-1}) d^{\times} t_1 d^{\times} t_2 \\
		&= \int_{\vO_v^{\times}} \int_{\vO_v^{\times}} \psi((t-1)t_2) \chi(t) d^{\times} t d^{\times} t_2 \\
		&= \sum_{a \in \vO_v^{\times} / U_v^l} \chi(a) \int_{\vO_v^{\times}} \int_{U_v^l} \psi((at-1)t_2) d^{\times} t d^{\times} t_2.
	\end{align*}
	Note that $(t-1)t_2 \in \varpi_v^l \vO_v \subset \varpi_v^r \vO_v$, hence $\psi((t-1)t_2) = 1$ in the above equation. We therefore continue
	\begin{align*}
		|G(\chi,\psi)|^2 &= \sum_{a \in \vO_v^{\times} / U_v^l} \chi(a) \int_{\vO_v^{\times}} \int_{U_v^l} \psi((a-1)tt_2) d^{\times} t d^{\times} t_2 \\
		&= {\rm Vol}(U_v^l) \sum_{n=0}^l \sum_{a \in (U_v^n - U_v^{n+1}) / U_v^l} \chi(a) \int_{\vO_v^{\times}} \psi(\varpi_v^n t_2) d^{\times}t_2 \\
		&= q_v^{-l+d_v/2} \zeta_v(1)^2 \sum_{n=0}^l \sum_{a \in (U_v^n - U_v^{n+1}) / U_v^l} \chi(a) \left( \int_{O_v} \psi(\varpi_v^n t_2) dt_2 - q_v^{-1} \int_{\vO_v} \psi(\varpi_v^{n+1} t_2) dt_2 \right) \\
		&= q_v^{-l} \zeta_v(1)^2 \sum_{n=0}^l \sum_{a \in (U_v^n - U_v^{n+1}) / U_v^l} \chi(a) (1_{n \geq l} - q_v^{-1} 1_{n+1 \geq l} ) \\
		&= q_v^{-l} \zeta_v(1)^2 \left( 1 - q_v^{-1} \sum_{a \in U_v^{l-1} / U_v^l} \chi(a) \right) \\
		&= q_v^{-l} \zeta_v(1)^2 \times \begin{cases} 1 & \text{if } l=r \\ 0 & \text{if } l > r, \end{cases}
	\end{align*}
	which confirms the result of the last case and concludes the proposition. \endproof
	
	Recall that the conductor of $\psi_v$ is $\varpi_v^{-d_v} \vO_v$. Take the convention $n(\varpi_v^0) = 1$. 
	\begin{lemma}
		Let $W$ transform as $\omega_v$ under the action of $a(\vO_v^{\times})$. Suppose the conductor of $\omega_v \chi_v$ is $1+\varpi_v^r\vO_v$. Then if $r >0, l \neq 0$,
		$$ |\zeta(s+1/2, n(\varpi_v^{-l})W, \chi_v, \psi_v)| = \zeta_v(1) q^{-r/2} q^{-\Re(s) (l-r-d_v)} |W(a(\varpi_v^{l-r-d_v}))|. $$
		If $r>0, l=0$, $$ \zeta(s+1/2, W, \chi_v, \psi_v)= 0. $$
		If $r=0,l \neq 0$,
		$$ \zeta(s+1/2, n(\varpi_v^{-l})W, \chi_v, \psi_v) = \sum_{k=l-d_v}^{\infty} W(a(\varpi_v^k))\chi_v(\varpi_v)^k q_v^{-sk}$$
		$$ - \frac{1}{q_v-1} W(a(\varpi_v^{l-d_v-1}))\chi_v(\varpi_v)^{l-d_v-1}q_v^{-s(l-d_v-1)}. $$
		If $r=0,l=0$,
		$$ \zeta(s+1/2, W, \chi_v, \psi_v) = \sum_k W(a(\varpi_v^k))\chi_v(\varpi_v)^k q_v^{-sk}. $$
	\label{NonArchiEstG}
	\end{lemma}
	\proof We notice that, for $l \neq 0$,
	\begin{eqnarray}
		& & \zeta(s+1/2, n(\varpi_v^{-l})W, \chi_v, \psi_v) = \int_{F_v^{\times}} W(a(y)) \psi_v(\varpi_v^{-l} y) \chi_v(y) |y|_v^s d^{\times} y \nonumber \\
		&=& \sum_{m \in \Z} W(a(\varpi_v^m)) \chi_v(\varpi_v)^m q_v^{-ms} G(\omega_v\chi_v, \psi_v(\varpi_v^{m-l} \cdot));
	\label{GenGaussWithTrans}
	 \end{eqnarray}
	while for $l=0$,
	\begin{equation}
		\zeta(s+1/2, W, \chi_v, \psi_v) = \sum_{m \in \Z} W(a(\varpi_v^m)) \chi_v(\varpi_v)^m q_v^{-ms} \int_{\vO_v^{\times}} \omega_v \chi_v(y) d^{\times} y.
	\label{GenGaussWithoutTrans}
	\end{equation}
	
	Let's consider the first case i.e. $r \geq 1$ and $l \neq 0$. We apply (\ref{GenGaussWithTrans}) with Proposition \ref{ClassicalGauss} and obtain
	$$ \zeta(s+1/2, n(\varpi_v^{-l})W, \chi_v, \psi_v) = W(a(\varpi_v^{l-r-d_v})) \chi_v(\varpi_v)^{l-r-d_v} q_v^{-(l-r-d_v)s} G(\omega_v \chi_v, \psi_v(\varpi_v^{-r-d_v} \cdot) ). $$
	The first assertion follows from Proposition \ref{ClassicalGauss}.
	
	Next, we consider the second case $r>0$ and $l=0$. We notice $ \int_{\vO_v^{\times}} \omega_v \chi_v(y) d^{\times} y = 0 $ since $\omega_v \chi_v$ is a non trivial multiplicative character on $\vO_v^{\times}$. The second assertion follows by applying (\ref{GenGaussWithoutTrans}).
	
	In the third case $r=0,l \neq 0$, we analyze
	\begin{eqnarray}
	& &	\int_{\vO_v^{\times}} \psi_v(\varpi_v^{m-l}y) \omega_v \chi_v(y) d^{\times} y = \int_{\vO_v^{\times}} \psi_v(\varpi_v^{m-l}y) d^{\times} y \nonumber \\
	&=&  q_v^{d_v/2} \zeta_v(1) \left( \int_{\vO_v} \psi_v(\varpi_v^{m-l}y) dy -  q_v^{-1} \int_{\vO_v} \psi_v(\varpi_v^{m-l+1}y) dy \right) \nonumber \\
	&=& \zeta_v(1) (1_{m-l \geq -d_v} - q_v^{-1} 1_{m-l+1 \geq -d_v} ).  \label{GaussSum}	
	\end{eqnarray}
	The third assertion then follows by applying (\ref{GenGaussWithTrans}) again.
	
	The last case is an easy consequence of (\ref{GenGaussWithoutTrans}). \endproof
	
	The following corollary is essentially Lemma 11.7 of \cite{V}.
	\begin{corollary}
		Let $r_v$ be the conductor of $\omega_v \chi_v$. We take $W_{0,v}(y)$ to be the new vector, or ``vecteur essentiel'' in the sense of (4.4) of \cite{JPS}, of $\pi_v$. Then if $r_v + d_v >0$,
		$$ |\zeta(s,n(\varpi_v^{-(r_v+d_v)})W_{0,v},\chi_v,\psi_v)| = q_v^{-r_v/2} \zeta_v(1). $$
		If $r_v + d_v =0$, $$\zeta(s,W_{0,v},\chi_v,\psi_v) = L(s, \pi_v\otimes \chi_v). $$
	\label{NonArchiEst}
	\end{corollary}
	\proof We only need to remark that Th\'eor\`eme 5.(\rmnum{2}) of \cite{JPS} implies for $\omega_v=1$, $r_v=0$ and $d_v=0$,
	$$ \zeta(s, W_{0,v}, \chi_v, \psi_v) = L(s, \pi_v \times \chi_v). $$
	The second assertion follows by noting our measure normalization and the general definition of the ``new vector'' for $\text{GL}_2$, in Theorem 4.24 of \cite{G} for example. Specializing the above equation to $\chi_v=1$ and comparing it with (\ref{GenGaussWithoutTrans}), we deduce
	$$ W_{0,v}(1) = 1. $$
	We then use the first assertion of Lemma \ref{NonArchiEstG} to get the first assertion. \endproof
	
	As a consequence
	\begin{align}
		\prod_{v<\infty} \left\lvert \frac{L(1/2, \pi_v\otimes \chi_v)}{\zeta(1/2,n(T_v)W_{0,v},\chi_v,\psi_v)} \right\rvert &\leq \prod_{v<\infty, \omega_v \chi_v \ {\rm ramified}} \frac{C(\omega_v \chi_v)^{1/2}}{(1-q_v^{-1/2+\theta})^2} \nonumber \\
		&\ll_{\epsilon, \pi, F} \prod_{v < \infty} C(\chi_v)^{1/2+\epsilon}. \label{LocalEst2}
	\end{align}
	Note that (\ref{LocalGoal}) is established by (\ref{LocalEst1}) and (\ref{LocalEst2}) once $T = (T_v)_v$ and $\varphi = n(T) \varphi_0$ are chosen, where $\varphi_0$ corresponds to $(W_{0,v})_v$, $T_v = \varpi_v^{-(r_v+d_v)}$.
	
	\begin{proposition}
		The function $\varphi_0$ corresponding to $\prod_v W_{0,v}$ in the Kirillov model of $\pi$ verifies
		$ \varphi_0 \in R^s. $
	\end{proposition}
	\proof This is an obvious consequence of the discussion in Section \ref{ConstAutoForm}. In fact it is easy to verify $\varphi_0 \in R_0^{\infty}$, then we apply Corollary \ref{Schwartz}. \endproof

	\subsection{A Calculation for Unitary Principal Series}
	We are interested in consequences of Lemma \ref{NonArchiEstG} in the case of a unitary principal series representation. Recall $v(\psi)=d_v$. For simplicity of notations, we omit the subscript $v$. Assume that the representation takes the form $\pi = \pi(\xi, \xi^{-1})$ for some unramified unitary character $\xi$ of $F^{\times}$. For an integer $m \geq 0$, we are interested in vectors of $\pi$ invariant by $K^0[m]$. Let $W_{\pi}$ denote the Whittaker model of $\pi$.
	\begin{proposition}
		If $W \in W_{\pi}$ is invariant by $K^0[m]$, then we have
		$$ |W(a(y))| \ll (v(y)+1)(m+1) q^{(m-d_v)/2} \lVert W \rVert |y|^{1/2} 1_{v(y)\geq -d_v} $$
		with the implicit constant being absolute.
	\end{proposition}
	\proof If we write $ \vO_n = \varpi^n\vO - \varpi^{n+1}\vO, n \geq 1$, then
	$$ D_0 =  B(\vO)wN(\vO) = \left\{ \begin{pmatrix} a & b \\ c & d \end{pmatrix} \in K: c \in \vO^{\times} \right\}, $$
	$$ D_n = B(\vO)N_-(\vO_n) =  \left\{ \begin{pmatrix} a & b \\ c & d \end{pmatrix} \in K: c \in \vO_n \right\}, n \geq 1, $$
	are double cosets of $G(O)$ w.r.t. $B(\vO)$ with mass (assuming ${\rm Vol}(G(\vO)) = 1$)
	\begin{equation}
		{\rm Vol}(D_0) = \frac{q}{q+1}, {\rm Vol}(D_n) = \frac{q^{-(n-1)}(1-q^{-1})}{q+1}.
	\label{MassDoubleCosets}
	\end{equation}
	Therefore for any $f \in \pi^{B(\vO)}$, there is a sequence of complex numbers $f_n, n \geq 0$ characterizing $f$ by
	$$ f\mid_{D_n} = f_n. $$
	Therefore, if
	$$ W(a(y)) = W_f(a(y)) = \xi^{-1}(y)|y|^{1/2} \int_F f(wn(x)) \psi(-xy) dx $$
	denotes the Whittaker function of $f$, then we obtain, with $t = \xi(\varpi)$,
	\begin{eqnarray*}
		W(a(y)) &=& \xi^{-1}(y)|y|^{1/2} \left( \int_{\vO} f\left(w \begin{pmatrix} 1 & x \\ & 1 \end{pmatrix} \right) \psi(-xy) dx \right. \\
		& & +\left. \sum_{n=1}^{\infty} \int_{\vO_n^{-1}} f\left( \begin{pmatrix} 1/x & -1 \\ & x \end{pmatrix} \begin{pmatrix} 1 & \\ 1/x & 1 \end{pmatrix} \right) \psi(-xy)  dx \right) \\
		&=& \xi^{-1}(y)|y|^{1/2} \left( f_0 \int_{\vO} \psi(-xy) dx + \sum_{n=1}^{\infty} f_n \int_{\vO_n^{-1}} \xi^{-2}(x) \psi(-xy) \frac{dx}{|x|} \right) \\
		&=& \xi^{-1}(y)|y|^{1/2} \left( q^{-d_v/2} f_0 1_{v(y) \geq -d_v} + \sum_{n=1}^{\infty} f_n t^{2n} \int_{\vO^{\times}} \psi(-\varpi^{-n}xy) \frac{dx}{|x|} \right) \\
		&=& q^{-d_v/2} \xi^{-1}(y)|y|^{1/2} \left( f_0 1_{v(y) \geq -d_v} + \sum_{n=1}^{\infty} f_n t^{2n} (1_{v(y)-n \geq -d_v} - q^{-1} 1_{v(y)-n+1 \geq -d_v}) \right),
	\end{eqnarray*}
	where we have used (\ref{GaussSum}). Hence
	\begin{eqnarray}
		W(a(y)) &=& q^{-d_v/2} \xi^{-1}(y)|y|^{1/2} 1_{v(y)\geq -d_v} \cdot \nonumber \\
		& & \left( f_0 - q^{-1} f_{v(y)+d_v+1} t^{2(v(y)+d_v+1)} + (1-q^{-1}) \sum_{n=1}^{v(y)+d_v} f_n t^{2n} \right). \label{v(y)small}
	\end{eqnarray}
	By the discussion in Section 3.1.6 of \cite{MV}, we have
	\begin{eqnarray}
		\lVert W \rVert^2 &=& \int_{G(\vO)} |f(k)|^2 dk = \frac{q^{d_v/2}}{1+q^{-1}} \int_F |f(wn(x))|^2 dx \nonumber \\
		&=& \frac{q^{d_v/2}}{1+q^{-1}} \left( |f_0|^2 \int_{\vO} dx + \sum_{n=1}^{\infty} |f_n|^2 \int_{\vO_n^{-1}} \frac{dx}{|x|^2} \right) \nonumber \\
		&=& (1+q^{-1})^{-1} \left( |f_0|^2 + \sum_{n=1}^{\infty} |f_n|^2 q^{-n}(1-q^{-1}) \right). \label{fGeneralNorm}
	\end{eqnarray}

	Now we assume in addition that $f$ is invariant by $K^0[m]$. Since
	$$ D_0, D_n, 1 \leq n \leq m-1, \cup_{n=m}^{\infty} D_n $$
	are the double cosets of $G(\vO)$ w.r.t. $B(\vO)$ and $K_0[m]$,	
	we have equivalently
	$$ f_n = f_m, \forall n \geq m. $$
	Consequently, (\ref{fGeneralNorm}) becomes
	\begin{equation}
		\lVert W \rVert^2 = (1+q^{-1})^{-1} \left( |f_0|^2 + \sum_{n=1}^{m-1} |f_n|^2 q^{-n}(1-q^{-1}) + |f_m|^2 q^{-m} \right).
	\label{fK0mInvNorm}
	\end{equation}
	Thus if $v(y)+d_v \geq m$, we rewrite
	\begin{eqnarray}
		W(a(y)) &=& q^{-d_v/2} \xi^{-1}(y)|y|^{1/2} 1_{v(y)\geq -d_v} \left(f_0 + (1-q^{-1}) \sum_{n=1}^{m-1} f_n t^{2n} + \right. \nonumber \\
		& & \left. \left(\frac{t^{2m}-t^{2(v(y)+d_v+1)}}{1-t^2}-q^{-1} \frac{t^{2m}-t^{2(v(y)+d_v+2)}}{1-t^2} \right)f_m \right). \label{v(y)big}
	\end{eqnarray}
	We apply Cauchy-Schwarz and compare (\ref{fK0mInvNorm}) with (\ref{v(y)small}) if $v(y)+d_v < m$, or (\ref{v(y)big}) if $v(y)+d_v \geq m$. The proposition follows. \endproof
	
	We apply the third case of Lemma \ref{NonArchiEstG} to the above $W=W_f$ and obtain for $\Re(s) = 1/2, \epsilon > 0$,
	\begin{align}
		|l^{|\cdot|^s}(n(\varpi^{-l})W)| &\leq \left( \left(\sum_{k=l-d_v}^{\infty} |W(a(\varpi^k))|^2 \right)^{1/2} \left(\sum_{k=l-d_v}^{\infty} q^{-k} \right)^{1/2} \right. \nonumber \\
		&\ \left. + |W(a(\varpi^{l-d_v-1}))| \frac{q^{-(l-d_v-1)/2}}{q-1} \right) \nonumber \\
		&\ll_{\epsilon} (m+1) q^{(m-d_v)/2} q^{-(l-d_v)(1-\epsilon)} \lVert W \rVert. \label{LocalZetaDecay}
	\end{align}
	
	The analogue of (\ref{LocalZetaDecay}) at an infinite place is just a consequence of integration by parts. Take the case of a real place for example. If $W \in W_{\pi}^{\infty}$ then we know that $W(a(y))$ is of rapid decay as $|y| \to \infty$, controlled by $|y|^{1/2-\epsilon}, \forall \epsilon > 0$ as $|y|\to 0$, as well as $X.W$ for any $X$ in the enveloping algebra of $G$. Consequently, for $\Re(s) = 1/2+\epsilon$, we have
	$$ l^{|\cdot|^s}(n(t)W) = -\frac{1}{t} \int_{F} \left( n(t).U.W(a(y)) |y|^{s-2} + (s-1)n(t).W(a(y))|y|^{s-2} \right) dy. $$
	The right side converges thanks to the upper bounds of $W(a(y)), U.W(a(y))$, where $U=\begin{pmatrix} 1 & \\ & 0 \end{pmatrix}$. We then use the local functional equation to see
	$$ l^{|\cdot|^{s-1}}(n(t)W) = \gamma(s-1/2, \pi, \psi)^{-1} l^{|\cdot|^{1-s}}(w.n(t).W). $$
	The gamma factor $\gamma(s-1/2, \pi, \psi) = \gamma(s-1/2, \xi, \psi)\gamma(s-1/2, \xi^{-1}, \psi)$ is of size $\asymp_{\epsilon} C(\xi)^{1-2\epsilon} |s|^{1-2\epsilon}$, while the integral is bounded, separating the contributions from $|y| \leq 1$ and from $|y| > 1$, as $\ll_{\epsilon} \lVert w.n(t).W \rVert + \lVert T.w.n(t).W \rVert $, with $T = \begin{pmatrix} 0 & 1 \\ 0 & 0  \end{pmatrix} $ in the Lie algebra of $G$. We do similar estimation for $n(t).U.W$. Using Theorem \ref{SobEquiv}, We thus find
	\begin{equation}
		|l^{|\cdot|^s}(n(t)W) | \ll_{\epsilon} |t|^{-1} |s|^{-2+\epsilon}C(\xi)^{-1+\epsilon} \lVert \Delta.W \rVert, \forall \epsilon > 0, \Re (s)=1/2+\epsilon.
	\label{LocalZetaDecayInfinite}
	\end{equation}
	Finally, note that if $v$ is a complex place, the proof of Theorem \ref{SobEquiv} given here implies that we should replace $\lVert \Delta.W \rVert$ by $\lVert \Delta^8.W \rVert$ in (\ref{LocalZetaDecayInfinite}).

\section{Some Special Decay of Matrix Coefficients}

	\subsection{Branching Law at a Finite Place}
	\label{BranchingLawSection}

	Although Theorem \ref{MatrixCoeffDecay} is general and convenient to use, it fails to be optimal in many cases of application. At least in its application to our situation, the dimension factor ``$\dim (K_v x_1)^{1/2} \dim (K_v x_2)^{1/2}$'' could be supressed in many places. This is in particular true for ``classical vectors'' and for diagonal matrices in $G_v$, at a finite place $v$.
	
	\textbf{Since we stick to a finite place, we shall omit the subscript $v$ everywhere in this section.}

	Write $\vp = \varpi \vO$. Let $\varepsilon_0$ be a character of $\vO^{\times}$ of conductor $\vp^{N_0}$. Then $\varepsilon_0$ determines as well a character of $B(\vO) = \bigcap_{N \geq 1} K^0[N]$ and a character of each $K^0[N]$ with $N \geq N_0$, in each case taking $\begin{pmatrix} a & b\\ c & d \end{pmatrix}$ to $\varepsilon_0(a)$. Define right regular $K$-representations
	$$ {\rm Ind}(\varepsilon_0) = {\rm Ind}_{B(\vO)}^K \varepsilon_0; {\rm Ind}_N(\varepsilon_0) = {\rm Ind}_{K^0[N]}^K \varepsilon_0, \forall N \geq N_0. $$
Clearly ${\rm Ind}_N(\varepsilon_0)$ is naturally embedded into ${\rm Ind}_{N+1}(\varepsilon_0)$ for any $N \geq N_0$. Let $u_{N_0}(\varepsilon_0) = {\rm Ind}_{N_0}(\varepsilon_0)$, and $u_N(\varepsilon_0)$ be the ortho-complement of ${\rm Ind}_{N-1}(\varepsilon_0)$ in ${\rm Ind}_N(\varepsilon_0)$. The following characterization can be found in \cite{Ca3}.
	\begin{proposition}
		$u_N(\varepsilon_0)$ is the unique irreducible representation of $K$ which
		\begin{itemize}
			\item[(1)] is trivial on $K[N]$ but not on $K[N-1]$;
			\item[(2)] has a vector $v \neq 0$ such that $u_N(g).v = \varepsilon_0(a)v$ for each $g=\begin{pmatrix} a & b \\ c & d \end{pmatrix} \in K^0[N]$.
		\end{itemize}
		$v$ satisfying (2) is unique up to scalar. We call it $v_N(\varepsilon_0)$.
	\label{Charu_N}
	\end{proposition}
	\begin{theorem}
		Let $\pi$ be a unitary irreducible representation of $G={\rm GL}_2(F)$ with conductor $v(\pi) = c$ and central character $\varepsilon$ whose restriction to $\vO^{\times}$ is $\varepsilon_0$ of conductor $N_0$. Then we have a decomposition of $K$-representations
		$$ {\rm Res}_K^G \pi = \pi^{K[c-1]} \bigoplus \bigoplus_{N \geq c} u_N(\varepsilon_0), $$
with the convention $\pi^{K[-1]} = \{ 0 \}$.
		Furthermore, $u_N(\varepsilon_0)$ can be identified with the ortho-complement of $\pi^{K[N-1]}$ in $\pi^{K[N]}$.
	\label{BranchingLaw}
	\end{theorem}
	\proof If $\pi$ is supercuspidal, this is just Theorem 1 of \cite{Ca3}. Assume first that $\pi = \pi(\chi_1, \chi_2)$ is in the principal or complementary series. We may assume $v(\chi_1) \geq v(\chi_2)$ after exchanging $\chi_1, \chi_2$ if necessary. For any fixed $N \geq \max (v(\chi_1), v(\chi_2))$, write $\overline{K} = K / K[N], \overline{B} = B(\vO) / B(\vO) \cap K[N]$. We naturally have identifications of $K$-representations
	$$ {\rm Res}_K^G \pi = {\rm Ind}_{B(\vO)}^K (\chi_1, \chi_2); \pi^{K[N]} \simeq {\rm Ind}_{\overline{B}}^{\overline{K}} (\chi_1, \chi_2). $$
By Frobenius reciprocity, we also have
	$$ {\rm Hom}_K ( \pi^{K[N]}, \pi^{K[N]} ) \simeq  {\rm Hom}_{\overline{K}} ( \pi^{K[N]}, \pi^{K[N]} ) \simeq {\rm Hom}_{\overline{B}} ((\chi_1, \chi_2), {\rm Ind}_{\overline{B}}^{\overline{K}} (\chi_1, \chi_2)), $$
which is the space of functions in ${\rm Ind}_{\overline{B}}^{\overline{K}} (\chi_1, \chi_2)$ transforming as $(\chi_1, \chi_2)$ under the right translation by $\overline{B}$. We denote its dimension by $d_N$, and note the double coset decomposition
	$$ \overline{K} = \overline{B} \coprod \overline{B} w \overline{B} \coprod \coprod_{k=1}^{N-1} \overline{B} n_-(\varpi^k) \overline{B}, w = \begin{pmatrix} & -1 \\ 1 & \end{pmatrix}, n_-(\varpi^k) = \begin{pmatrix} 1 & \\ \varpi^k & 1 \end{pmatrix}. $$
The contribution of $\overline{B}$ to $d_N$ is $1$; the contribution of $\overline{B} w \overline{B}$ to $d_N$ is $1_{\chi_1 = \chi_2}$. The contribution of $\overline{B} n_-(\varpi^k) \overline{B}$ to $d_N$ is $1$ iff
	$$ \chi_1(a) \chi_2(d) = \chi_1(a-b\varpi^k) \chi_2(d+b\varpi^k) $$
for all $a,d\in \vO^{\times}, b \in \vO$ satisfying $a-d-b\varpi^k \in \vp^N$, which is equivalent to $v(\chi_1 \chi_2^{-1}) \geq k$. We deduce that
	\begin{equation}
		d_N = N+1-v(\chi_1 \chi_2^{-1}), \forall N \geq \max(v(\chi_1), v(\chi_2)) = v(\chi_1).
	\label{FrobDimF}
	\end{equation}
If we write $u_N'$ to be the $K$-representation which is the ortho-complement of $\pi^{K[N-1]}$ in $\pi^{K[N]}$, the above formula shows that $u_N'$ is irreducible if $N > v(\chi_1)$. Since $\pi^N$ has a subspace of dimension $N-c+1$ of functions transforming as $(\varepsilon_0,1)$ under $B(\vO)$ by the theory of conductor (c.f.  Theorem 1 of \cite{Ca2} or Theorem 4.24 of \cite{G}), $u_N'$ has one such nonzero function as long as $N \geq c=v(\chi_1) + v(\chi_2) \geq N_0$. By Proposition \ref{Charu_N}, $u_N' \simeq u_N(\varepsilon_0)$ if $N \geq \max(c, v(\chi_1)+1)$. If $c \geq 1+v(\chi_1)$, then we are done. Otherwise, we must have $v(\chi_2) = 0$. Consequently $c = v(\chi_1) = v(\chi_1\chi_2^{-1}) = N_0$, and $d_c = 1$, which shows that $\pi^{K[c]}$ is irreducible. But $\pi^{K[c]}$ contains a nonzero vector $v_0$ transforming as $(\varepsilon_0,1)$ under $B(\vO)$ by definition of conductor, hence $\pi^{K[c]} \simeq u_c(\varepsilon_0)$ by Proposition \ref{Charu_N} if we can prove $\pi^{K[c]} \neq \pi^{K[c-1]}$. If $c=0$ we are done by convention. Otherwise,  $v_0 \notin \pi^{K[c-1]}$ since $v(\varepsilon_0) = c$. Hence $\pi^{K[c-1]} \varsubsetneq \pi^{K[c]}$, and we are done for $\pi$ in principal series.

	Assume at last $\pi = \pi(\chi_1, \chi_2)$ is a special representation, with $\chi_1 \chi_2^{-1} = |\cdot|^{-1}$. Then as $K$-representations, ${\rm Res}_K^G \pi$ is $\tilde{\pi} = {\rm Ind}_{B(\vO)}^K (\chi_1,\chi_2)$ quotient by the one dimensional subspace spanned by $0 \neq f_0 \in \tilde{\pi}$ defined by
	$$ f_0(k) = \chi_1(\det k), \forall k \in K. $$
Note that, up to scalar, $f_0$ spans $\chi_1 \circ \det \in \widehat{K}$, lies in $\tilde{\pi}^{K[N]}$ for any $N \geq \max(v(\chi_1), v(\chi_2)) = v(\chi_1)$, and $\tilde{\pi}$ is semi-simple as a $K$-representation. Hence we can apply (\ref{FrobDimF}) to $\tilde{\pi}$ to get
	$$ d_N = {\rm Hom}_K(\pi^{K[N]}, \pi^{K[N]}) = N, \forall N \geq v(\chi_1). $$
But $c=\max(1, 2v(\chi_1)) \geq v(\chi_1)+1$ in this case, hence $u_N'$ similarly defined as in the previous case makes sense for $N\geq c$ and is irreducible by the above formula. It also contains a nonzero vector transforming as $(\varepsilon_0,1)$ under $B(\vO)$ by Theorem 4.24 of \cite{G}, hence $u_N' \simeq u_N(\varepsilon_0)$ and we are done. \endproof
	\begin{remark}
		In the case $v(\chi_1\chi_2^{-1}) = \max(v(\chi_1), v(\chi_2)) = c'$, the above argument actually gives the complete branching law. For $\pi$ in principal series, we get
		$$ {\rm Res}_K^G \pi = \bigoplus_{c' \leq N < c} u_N' \bigoplus \bigoplus_{N \geq c} u_N(\varepsilon_0), $$
with each component $K$-irreducible. For $\pi$ special, $c'=0, c=1$ and $\varepsilon_0 = 1$, we get
		$$ {\rm Res}_K^G \pi = \bigoplus_{N \geq 1} u_N(1). $$
		In particular, in the general case, we get complete branching law for $\pi \otimes \chi_2^{-1}$, i.e. branching law up to twisting by a character. The multiplicity one holds for both principal and special series.
	\end{remark}
	\begin{definition}
		Let $\pi$ be as in Theorem \ref{BranchingLaw}. A vector $v \in \pi$ is called classical if
		$$ \pi(g).v = \varepsilon_0(a).v, \forall g = \begin{pmatrix} a & b \\  & d \end{pmatrix} \in B(\vO). $$
	\end{definition}
It is easy to see that the space of classical vectors is spanned by $v_N(\varepsilon_0), N\geq c$ under the isomorphism in Theorem \ref{BranchingLaw}. We write the corresponding vectors in $\pi$ by $v_N(\pi), N \geq c$.

	\subsection{Matrix Coefficients for Classical Vectors}

	\begin{proposition}
		Let $\pi$ be as in Theorem \ref{BranchingLaw}. If $\pi$ is tempered, we have
		$$ \langle a(y).v_N(\pi) , v_N(\pi) \rangle \leq \lVert v_N(\pi) \rVert^2 \Xi(a(y)), \forall y \in F^{\times}, a(y) = \begin{pmatrix} y & \\ & 1 \end{pmatrix}. $$
		If $\pi$ is not tempered, then for any $\epsilon > 0$ we have
		$$ \langle a(y).v_N(\pi) , v_N(\pi) \rangle \ll_{\epsilon} \lVert v_N(\pi) \rVert^2 \Xi(a(y))^{1-2\theta - \epsilon}, \forall y \in F^{\times}, a(y) = \begin{pmatrix} y & \\ & 1 \end{pmatrix}. $$
		Here $\Xi = \Xi_v$ is the Harish-Chandra's function defined in section \ref{DefOfXi}.
	\label{MCDwoD}
	\end{proposition}
	We are going to prove Proposition \ref{MCDwoD} by giving an explicit description of $v_N(\pi)$ in some suitable model of $\pi$. Recall (c.f. Proposition \ref{UniqueRep}) the (completed) Kirillov model $K_{\pi}^{\psi}$ of $\pi$ is the space of functions in $L^2(F^{\times}, d^{\times}x)$ with the action of $B=B(F)$ given by
	$$ \begin{pmatrix} a & b \\ 0 & 1 \end{pmatrix}.f(x) = \psi(bx)f(ax), \begin{pmatrix} u & \\ & u \end{pmatrix}.f(x) = \varepsilon(u) f(x), \forall a,u \in F^{\times}, b \in F. $$
Since $F^{\times} \simeq \vO^{\times} \times \Z$, we get a model by doing partial Fourier transform on $\vO^{\times}$.
	\begin{definition}
		The dual Kirillov model $\widehat{K}_{\pi}^{\psi}$ of $\pi$ is the space of functions $F:\widehat{\vO^{\times}} \to \C[[t]]$ s.t. for any $\nu \in \widehat{\vO^{\times}}$, $F(\nu,t) = \sum_{n \in \Z} F_n(\nu)t^n$ with $F_n(\nu) \in \C$, and
		$$ \lVert F \rVert^2 = \sum_{\nu \in \widehat{\vO^{\times}}} \sum_{n \in \Z} |F_n(\nu)|^2 < \infty. $$
	\end{definition}
In fact, in order to pass from $K_{\pi}^{\psi}$ to $\widehat{K}_{\pi}^{\psi}$, we take a $f \in L^2(F^{\times})$ and define
	$$ F_n(\nu) = \int_{\vO^{\times}} f(\varpi^n u) \nu(u) d^{\times} u. $$
	\begin{remark}
		The dual Kirillov model is extensively used in \cite{JL} and \cite{Ca3}. More precisely, the model they used is the subspace of $\widehat{K}_{\pi}^{\psi}$ of smooth vectors, which we shall, by abus of language, still call the dual Kirillov model.
	\end{remark}
In the dual Kirillov model, we have
	$$ \begin{pmatrix} \delta \varpi^l & 0 \\ 0 & 1 \end{pmatrix}.F(\nu,t) = t^{-l} \nu(\delta)^{-1} F(\nu,t), $$
	$$ w.F(\nu,t) = C(\nu,t)F(\nu^{-1}\varepsilon_0^{-1},t^{-1}z_0^{-1}), z_0=\varepsilon(\varpi), w=\begin{pmatrix} 0 & -1 \\ 1 & 0 \end{pmatrix}, $$
where $C(\nu,t) \in \C[[t]][1/t]$ characterize $\pi$. We give $C(\nu,t)$ w.r.t. different series to which $\pi$ belongs, which can be found in \cite{JL} and is essentially the local functional equations. The following observation, which is just Lemma 2 of \cite{Ca3}, is important for our discussion.
	\begin{lemma}
		The group $K[N]$ is generated by $B[N] = Z_vN_vA_v \cap K[N]$ and the Weyl element $w$. Hence, $v \in \pi$ is fixed by $K[N]$ if and only if $v$ and $w.v$ are fixed by $B[N]$.
	\end{lemma}

		\subsubsection{$\pi$ is supercuspidal}

	The case of a supercuspidal representation is treated detailly in \cite{Ca3}. We recall the main results without proof. Only the last two corollaries are not in \cite{Ca3}.
	\begin{lemma}
		There is $n_{\nu} \in \Z$ with $n_{\nu} = -v(\pi \otimes \nu) \leq -2$ s.t. for some $C_0(\nu) \in \C^{\times}$,
		$$ C(\nu, t) = C_0(\nu) t^{n_{\nu}}. $$
	\end{lemma}
	\begin{corollary}
		Let $N \geq N_0 = v(\varepsilon_0)$. The space $\pi^{K[N]}$ corresponds to functions $F(\nu, t)$ in the dual Kirillov model satisfying
		\begin{itemize}
			\item[(1)] $F(\nu,t) = 0$ unless $v(\nu) \leq N$;
			\item[(2)] $F_n(\nu) = 0$ unless $-N \leq n \leq n_{\nu} + N$.
		\end{itemize}
	\label{Behaviorat0-1}
	\end{corollary}
	\begin{corollary}
		The space $u_N(\varepsilon_0), N\geq c=v(\pi) (> N_0)$ as in the decomposition in Theorem \ref{BranchingLaw} corresponds to functions $F(\nu, t)$ in the dual Kirillov model satisfying
		\begin{itemize}
			\item[(1)] $F(\nu,t) = 0$ unless $v(\nu) \leq N$;
			\item[(2)] if $v(\nu) \leq N-1$, then $F_n(\nu) = 0$ unless $n=-N$ or $n_{\nu}+N$;
			\item[(3)] if $v(\nu) = N$, then $F_n(\nu) = 0$ unless $-N \leq n \leq n_{\nu}+N$.
		\end{itemize}
	\end{corollary}
	\begin{corollary}
		The unique classcial vector $v_N(\pi)$ of $u_N(\varepsilon_0)$ ($N \geq c$) corresponds to the function $F(\nu, t)$ in the dual Kirillov model satisfying
		\begin{itemize}
			\item[(1)] $F(\nu,t) = 0$ unless $\nu = \varepsilon_0^{-1}$;
			\item[(2)] $F(\varepsilon_0^{-1}, t) = C t^{N-c}$ for some $C\in \C$.
		\end{itemize}
	\end{corollary}
	\begin{corollary}
		If $\pi$ is supercuspidal, then for $N \geq c=c(\pi)$ we have
		$$ | \langle a(y).v_N(\pi), v_N(\pi) \rangle | = 1_{v(y)=0} \lVert v_N(\pi) \rVert^2. $$
	\label{MCDwoD-1}
	\end{corollary}

		\subsubsection{$\pi$ is a principal or complementary series}

	Assume $\pi = \pi(\mu_1, \mu_2)$ with $\mu_1, \mu_2$ quasi-characters of $F^{\times}$. We fix a $\psi$ s.t. $v(\psi)=0$. For any $\mu \in \widehat{\vO^{\times}}$ and $y \in F^{\times}$, define the Gauss sum as in \cite{JL},
	$$ \eta(\mu, y) = \int_{\vO^{\times}} \mu(x)\psi(xy) d^{\times}y. $$
We also define the root number $r(\mu)$ if $v(\mu) = n > 0$ as
	$$ r(\mu) = \mu(-1)\mu(\varpi)^n \eta(\mu, \varpi^{-n})^{-1} q^{-n/2}. $$
	\begin{lemma}
		The local functional equations imply:
		\begin{itemize}
			\item[(1)] If $v(\mu_2 \nu^{-1} \varepsilon^{-1}) = v(\nu^{-1} \mu_1^{-1}) = n_1 >0$ and $v(\mu_1 \nu^{-1} \varepsilon^{-1}) = v(\nu^{-1} \mu_2^{-1}) = n_2 >0$, then we have
			$$ C(\nu, t) = r(\nu^{-1} \mu_1^{-1}) r(\nu^{-1} \mu_2^{-1}) \nu(\varpi)^{n_1+n_2} t^{-n_1-n_2}. $$
			\item[(2)] If $v(\mu_2 \nu^{-1} \varepsilon^{-1}) = v(\nu^{-1} \mu_1^{-1}) = n_1 >0$ and $v(\mu_1 \nu^{-1} \varepsilon^{-1}) = v(\nu^{-1} \mu_2^{-1}) = 0$, then we have
			$$ C(\nu,t) = r(\nu^{-1}\mu_1^{-1}) \nu(\varpi)^{n_1} t^{-n_1} \frac{1-\mu_2(\varpi)^{-1} q^{-1/2} t^{-1} }{1-\mu_2(\varpi) q^{-1/2} t}. $$
			\item[(3)] If $v(\mu_2 \nu^{-1} \varepsilon^{-1}) = v(\nu^{-1} \mu_1^{-1}) = 0$ and $v(\mu_1 \nu^{-1} \varepsilon^{-1}) = v(\nu^{-1} \mu_2^{-1}) = n_2 > 0$, then we have
			$$ C(\nu,t) = r(\nu^{-1}\mu_2^{-1}) \nu(\varpi)^{n_2} t^{-n_2} \frac{1-\mu_1(\varpi)^{-1} q^{-1/2} t^{-1} }{1-\mu_1(\varpi) q^{-1/2} t}. $$
			\item[(4)] If $v(\mu_2 \nu^{-1} \varepsilon^{-1}) = v(\nu^{-1} \mu_1^{-1}) = 0$ and $v(\mu_1 \nu^{-1} \varepsilon^{-1}) = v(\nu^{-1} \mu_2^{-1}) = 0$, then we have
			$$ C(\nu,t) = \frac{1-\mu_1(\varpi)^{-1} q^{-1/2} t^{-1} }{1-\mu_1(\varpi) q^{-1/2} t} \frac{1-\mu_2(\varpi)^{-1} q^{-1/2} t^{-1} }{1-\mu_2(\varpi) q^{-1/2} t}. $$
		\end{itemize}
	\end{lemma}
	\begin{corollary}
		Let $N \geq \max (v(\mu_1), v(\mu_2))$. The space $\pi^{K[N]}$ corresponds to functions $F(\nu,t)$ in the dual Kirillov model satisfying
		\begin{itemize}
			\item[(1)] $F(\nu,t) = 0$ unless $v(\nu) \leq N$.
			\item[(2)] If $v(\nu^{-1}\mu_1^{-1}) = n_1 > 0$, and $v(\nu^{-1}\mu_2^{-1}) = n_2 > 0$, then $F_n(\nu) = 0$ unless $-N \leq n \leq N-n_1-n_2$.
			\item[(3)] If $v(\nu^{-1}\mu_1^{-1}) = n_1 > 0$, and $v(\nu^{-1}\mu_2^{-1}) = 0$, then $F_n(\nu) = 0$ unless $n \geq -N$, and we have
			$$ F_{k+N-n_1}(\nu) = \mu_1(\varpi)^k q^{-k/2} F_{N-n_1}(\nu), \forall k \geq 0. $$
			\item[(4)] If $v(\nu^{-1}\mu_1^{-1}) = 0$, and $v(\nu^{-1}\mu_2^{-1}) = n_2 > 0$, then $F_n(\nu) = 0$ unless $n \geq -N$, and we have
			$$ F_{k+N-n_2}(\nu) = \mu_2(\varpi)^k q^{-k/2} F_{N-n_2}(\nu), \forall k \geq 0. $$
			\item[(5)] If $v(\nu^{-1}\mu_1^{-1}) = 0$, and $v(\nu^{-1}\mu_2^{-1}) = 0$, then $F_n(\nu) = 0$ unless $n \geq -N$, and we have
			$$ F_{n+2}(\nu) - (\mu_1(\varpi) + \mu_2(\varpi)) q^{-1/2} F_{n+1}(\nu) + \mu_1(\varpi)\mu_2(\varpi) q^{-1} F_n(\nu) = 0, \forall n \geq N-1. $$
		\end{itemize}
	\label{Behaviorat0-2}
	\end{corollary}
In fact, the corollary follows from the lemma by applying the following (obvious) proposition.
	\begin{proposition}
		Let $F(t) \in t^{-N} \C[[t]]$. Suppose $P(t), Q(t) \in \C[t]$ with $P(0) \neq 0$ s.t. $\frac{Q(t^{-1})F(t^{-1})}{P(t)} \in t^{-N'}\C[[t]]$, where $Q(t^{-1})F(t^{-1})$ is viewed as in $\C[[t^{-1}]]$, $\frac{1}{P(t)}$ in $\C[[t]]$. Then $F(t) = t^{-N} \frac{S(t)}{Q(t)}$ for some $S(t) \in \C[t]$ with $\deg S \leq N+N'$.
	\end{proposition}
	\begin{corollary}
		The unique classcial vector $v_N(\pi)$ of $u_N(\varepsilon_0)$ ($N \geq c = v(\mu_1) + v(\mu_2)$) corresponds to the function $F(\nu, t)$ in the dual Kirillov model satisfying $F_n(\nu) = 0$ unless $\nu = \varepsilon_0^{-1}$ and $n \geq 0$, and
		\begin{itemize}
			\item[(1)] If $v(\mu_1) = n_1 > 0$ and $v(\mu_2) = n_2 > 0$, then up to a constant factor, $F(\varepsilon_0^{-1}, t) = t^{N-n_1-n_2}, \forall N \geq c = n_1+n_2$.
			\item[(2)] If $v(\mu_1) = n_1 > 0$ and $v(\mu_2) = 0$, then for $N=c=n_1$, up to a constant factor
			$$ F(\varepsilon_0^{-1}, t) = F(\mu_1^{-1},t) = \frac{1}{1-\mu_1(\varpi)q^{-1/2}t}; $$
while for $N > c$, up to a constant factor
			$$ F(\varepsilon_0^{-1}, t) = F(\mu_1^{-1},t) = -\frac{1-\overline{\mu_1(\varpi)}q^{-1/2}}{1-|\mu_1(\varpi)|^2q^{-1}}t^{N-n_1-1} + \frac{1}{1-\mu_1(\varpi)q^{-1/2}t} t^{N-n_1}. $$
			\item[(3)] If $v(\mu_1) = 0$ and $v(\mu_2) = n_2 > 0$, then for $N=c=n_2$, up to a constant factor
			$$ F(\varepsilon_0^{-1}, t) = F(\mu_2^{-1},t) = \frac{1}{1-\mu_2(\varpi)q^{-1/2}t}; $$
while for $N > c$, up to a constant factor
			$$ F(\varepsilon_0^{-1}, t) = F(\mu_2^{-1},t) = -\frac{1-\overline{\mu_2(\varpi)}q^{-1/2}}{1-|\mu_2(\varpi)|^2q^{-1}}t^{N-n_2-1} + \frac{1}{1-\mu_2(\varpi)q^{-1/2}t}t^{N-n_2}. $$
			\item[(4)] If $v(\mu_1) = v(\mu_2) = 0$ and $\mu_1 \neq \mu_2$, then for $N=0=c$, up to a constant factor
			$$ F(\varepsilon_0^{-1},t) = F(1,t) = \frac{1}{(1-\mu_1(\varpi)q^{-1/2}t)(1-\mu_2(\varpi)q^{-1/2}t)}; $$
for $N=1$, up to a constant factor, with
			$$ A = 1-|\mu_1(\varpi)|^2q^{-1}, B = 1-|\mu_2(\varpi)|^2q^{-1}, C = 1-\mu_1\overline{\mu_2}(\varpi)q^{-1}, $$
			$$ F(\varepsilon_0^{-1},t) = F(1,t) =\frac{AC}{1-\mu_1(\varpi)q^{-1/2}t} - \frac{B\overline{C}}{1-\mu_2(\varpi)q^{-1/2}t}; $$
while for $N > 1$, up to a constant factor, with $D = \overline{\mu_1}\overline{\mu_2}(\varpi)(\mu_1(\varpi) - \mu_2(\varpi))q^{-3/2}$,
			$$ F(\varepsilon_0^{-1},t) = F(1,t) = t^{N-2}\left( D + \frac{ACt}{1-\mu_1(\varpi)q^{-1/2}t} - \frac{B\overline{C}t}{1-\mu_2(\varpi)q^{-1/2}t} \right). $$
			\item[(5)] If $\mu_1 = \mu_2 = \mu$ with $v(\mu)=0$, then for $N=0=c$, up to a constant factor,
			$$ F(\varepsilon_0^{-1},t) = F(1,t) = \frac{1}{(1-\mu(\varpi)q^{-1/2}t)^2}; $$
for $N=1$, up to a constant factor,
			$$ F(\varepsilon_0^{-1},t) = F(1,t) = \frac{1+|\mu(\varpi)|^2q^{-1}}{1-\mu(\varpi)q^{-1/2}t} - \frac{1-|\mu(\varpi)|^2q^{-1}}{(1-\mu(\varpi)q^{-1/2}t)^2}; $$
while for $N > 1$, up to a constant factor,
			$$ F(\varepsilon_0^{-1},t) = F(1,t) = t^{N-2}\left( -\frac{\overline{\mu(\varpi)} |\mu(\varpi)|^2 q^{-3/2} }{1-|\mu(\varpi)|^2 q^{-1}} + \frac{1+|\mu(\varpi)|^2q^{-1}}{1-\mu(\varpi)q^{-1/2}t}t - \frac{1-|\mu(\varpi)|^2q^{-1}}{(1-\mu(\varpi)q^{-1/2}t)^2}t \right). $$
		\end{itemize}
	\end{corollary}
	\begin{corollary}
		If $\pi$ is a principal unitary series or a complementary series representation, then for $N \geq c=v(\pi)$ we have
		$$ | \langle a(y).v_N(\pi), v_N(\pi) \rangle | = 1_{v(y)=0} \lVert v_N(\pi) \rVert^2 $$
except in the following cases:
		\begin{itemize}
			\item[(1)]	$\pi = \pi(\mu_1, \mu_2)$ with $v(\mu_1)>0, v(\mu_2) = 0$ or $v(\mu_1)=0, v(\mu_2) > 0, N=c$, then
			$$ | \langle a(y).v_c(\pi), v_c(\pi) \rangle | = q^{-|v(y)|/2} \lVert v_N(\pi) \rVert^2. $$
			\item[(2)] $\pi = \pi(\mu_1, \mu_2)$ with $v(\mu_1)=v(\mu_2) = 0, \mu_1 \neq \mu_2$, then with $t_1=\mu_1(\varpi), t_2=\mu_2(\varpi)$, for $N=0=c$,
			$$ \frac{| \langle a(y).v_0(\pi), v_0(\pi) \rangle |}{\lVert v_0(\pi) \rVert^2} = \frac{q^{-|v(y)|/2}}{1+q^{-1}} \left| \frac{t_1^{|v(y)|+1} - t_2^{|v(y)|+1}}{t_1-t_2} -q^{-1}t_1t_2\frac{t_1^{|v(y)|-1} - t_2^{|v(y)|-1}}{t_1-t_2} \right|; $$
while for $N=1$,
			$$ \frac{| \langle a(y).v_1(\pi), v_1(\pi) \rangle |}{\lVert v_1(\pi) \rVert^2} = \frac{q^{-|v(y)|/2}}{1+q^{-1}} \left| q^{-1}\frac{t_1^{|v(y)|+1} - t_2^{|v(y)|+1}}{t_1-t_2} -t_1t_2\frac{t_1^{|v(y)|-1} - t_2^{|v(y)|-1}}{t_1-t_2} \right|. $$
			\item[(3)] $\pi = \pi(\mu_1, \mu_2)$ with $\mu_1=\mu_2=\mu, v(\mu)=0$, then for $N=0=c$,
			$$ \frac{| \langle a(y).v_0(\pi), v_0(\pi) \rangle |}{\lVert v_0(\pi) \rVert^2} = q^{-|v(y)|/2} \left( 1+|v(y)| \frac{1-q^{-1}}{1+q^{-1}} \right); $$
while for $N=1$,
			$$ \frac{| \langle a(y).v_1(\pi), v_1(\pi) \rangle |}{\lVert v_1(\pi) \rVert^2} = q^{-|v(y)|/2} \left( 1-|v(y)| \frac{1-q^{-1}}{1+q^{-1}} \right). $$
		\end{itemize}
	\label{MCDwoD-2}
	\end{corollary}

		\subsubsection{$\pi$ is a special representation}

	Write $\pi = \pi(\mu |\cdot|^{1/2}, \mu |\cdot|^{-1/2})$ with $\mu \in \widehat{F^{\times}}$.
	\begin{lemma}
		The local functional equations imply:
		\begin{itemize}
			\item[(1)] If $v(\nu^{-1}\mu^{-1}) = n > 0$, then we have
			$$ C(\nu,t) = r(\nu^{-1}\mu^{-1})^2 \nu(\varpi)^{2n} t^{-2n}. $$
			\item[(2)] If $v(\nu^{-1}\mu^{-1}) = 0$, then we have
			$$ C(\nu,t) = -\mu(\varpi)^{-1}t^{-1} \frac{1-\mu(\varpi)^{-1}q^{-1}t^{-1}}{1-\mu(\varpi)q^{-1}t}. $$
		\end{itemize}
	\end{lemma}
	\begin{corollary}
		Let $N \geq v(\mu)$. The space $\pi^{K[N]}$ corresponds to functions $F(\nu,t)$ in the dual Kirillov model satifying
		\begin{itemize}
			\item[(1)] $F(\nu,t) = 0$ unless $v(\nu) \leq N$.
			\item[(2)] If $v(\nu^{-1}\mu^{-1}) = l > 0$, then $F_n(\nu) = 0$ unless $-N \leq n \leq N-2l$.
			\item[(3)] If $v(\nu^{-1}\mu^{-1}) = 0$, then $F_n(\nu) = 0$ unless $ n \geq -N$, and we have
			$$ F_{k+N-1}(\nu) = \mu(\varpi)^k q^{-k} F_{N-1}(\nu), \forall k \geq 0. $$
		\end{itemize}
	\label{Behaviorat0-3}
	\end{corollary}
The proof is the same as for principal and complementary series.
	\begin{corollary}
		The unique classcial vector $v_N(\pi)$ of $u_N(\varepsilon_0)$ ($N \geq c = \max (2v(\mu), 1)$) corresponds to the function $F(\nu,t)$ in the dual Kirillov model satisfying $F_n(\nu) = 0$ unless $\nu=\varepsilon_0^{-1}$ and $n \geq 0$, and
		\begin{itemize}
			\item[(1)] If $v(\mu) = l > 0$, then up to a constant factor, $F(\varepsilon_0^{-1},t) = t^{N-2l}$.
			\item[(2)] If $v(\mu) = 0$, then for $N=1=c$, up to a constant factor,
			$$ F(\varepsilon_0^{-1},t) = F(1,t) = \frac{1}{1-\mu(\varpi)q^{-1}t}; $$
while for $N > 1$, up to a constant factor,
			$$ F(\varepsilon_0^{-1},t) = F(1,t) = - \frac{\mu(\varpi)^{-1}q^{-1}}{1-q^{-2}}t^{N-2} + \frac{t^{N-1}}{1-\mu(\varpi)q^{-1}t}. $$
		\end{itemize}
	\end{corollary}
	\begin{corollary}
		If $\pi$ is a special representation, then for $N \geq c=v(\pi)$ we have
		$$ |\langle a(y).v_N(\pi), v_N(\pi) \rangle | = 1_{v(y)=0} \lVert v_N(\pi) \rVert^2 $$
except in the case when $c = 1, N=1$, then
		$$ |\langle a(y).v_1(\pi), v_1(\pi) \rangle | = q^{-l} \lVert v_N(\pi) \rVert^2. $$
	\label{MCDwoD-3}
	\end{corollary}
Proposition \ref{MCDwoD} follows easily from Corollary \ref{MCDwoD-1}, \ref{MCDwoD-2} and \ref{MCDwoD-3}. The non-tempered case follows the same argument as in the proof of Theorem \ref{MatrixCoeffDecay}.
	\begin{remark}
		Note that Corollary \ref{Behaviorat0-1}, \ref{Behaviorat0-2} and \ref{Behaviorat0-3} also give a proof of (\ref{Whittakerat0}) in the case of a finite place. See Remark \ref{RWhittakerat0}.
	\end{remark}
	\begin{remark}
		Our discussion shows that the Gram-Schmidt procedure described in (38) of \cite{BH3} is simple, at least locally, i.e. there is $M \leq 2$ s.t. $v_{N+1}(\pi) = a(\varpi).v_N(\pi), \forall N \geq M$. It can be seen from the above explicit description of $v_N(\pi)$. But we wonder if a direct proof exists.
	\end{remark}

\section{Global Estimation}

	\subsection{Truncation}
	The goal of this section is to establish Lemma \ref{Tronc}.

Fix a function $h_0 \in C^{\infty}(\R^+)$ supported in $(0,2]$ such that $h_0 \mid_{(0,1]} = 1$ and $0 < h_0 < 1$. Denote by $\mathcal{M}(\cdot)$ the Mellin transform. For any $A>0$, let $h_{0,A}(t) = h_0(t/A)$. The following relation is immediate:
	$$
		|\mathcal{M}(\sigma * h_{0,Q^{-\kappa-1}})(s)| \leq 2^{|\Re (s)|}Q^{-(\kappa+1)\Re (s)} |\mathcal{M}(h_0)(s)|.
	$$
For any $t>0$, choose $y_t \in \A^{\times}$ such that $|y_t| = t$, and define
	$$
	 	f(t) = \int_{F^{\times} \backslash \A^{(1)}} \varphi(a(yy_t)) \chi(yy_t) d^{\times}y,
	$$
then
	$$
		l^{\chi,\sigma * h_{0,Q^{-\kappa-1}}}(\varphi) = \int_0^{+\infty} \sigma * h_{0,Q^{-\kappa-1}}(t) f(t) d^{\times} t.
	$$
Note that $\mathcal{M}(f)(s) = l^{\chi |\cdot|^s}(\varphi)$, Mellin inversion gives
	\begin{align*}
		|l^{\chi,\sigma * h_{0,Q^{-\kappa-1}}}(\varphi)| &=  \left\lvert \int_{\Re (s) = -1/2-\epsilon} \mathcal{M}(\sigma * h_{0,Q^{-\kappa-1}})(-s)  l^{\chi |\cdot|^s}(\varphi) \frac{ds}{2\pi i} \right\rvert \\
		&\ll Q^{-(\kappa+1)(1/2+\epsilon)}  \int_{\Re (s) = -1/2-\epsilon} |\mathcal{M}(h_0)(-s)  l^{\chi |\cdot|^s}(\varphi)| ds.
	\end{align*}
According to (\ref{JL}), one can write
	\begin{align*}
		l^{\chi |\cdot|^s}(\varphi) &= L(s+1/2, \pi \otimes \chi) \prod_{v | \infty} l^{\chi_v |\cdot |_v^s} (n(T_v)W_{0,v}) \prod_{v < \infty} \frac{l^{\chi_v |\cdot |_v^s} (n(T_v)W_{0,v})}{L(s+1/2, \pi_v \otimes \chi_v)} \\
		&= L^{(S)}(s+1/2, \pi \otimes \chi) \prod_{v \in S} l^{\chi_v |\cdot |_v^s} (n(T_v)W_{0,v}),
	\end{align*}
	where $S$ is the finite subset of places $v$ for which $T_v \neq 0$ or $\pi_v$ is ramified. From Corollary \ref{ArchiUpperEst} and Corollary \ref{NonArchiEst}, one sees that for each $v \in S$, $|l^{\chi_v |\cdot |_v^s} (n(T_v)W_{0,v})| \ll_{\epsilon, \varphi_0} C(\chi_v)^{-1/2 + \epsilon}$ and the product of the implicit constants tends to $0$ as $S$ increases. So
	$$
		\prod_{v \in S} l^{\chi_v |\cdot |_v^s} (n(T_v)W_{0,v}) \ll_{\epsilon, \varphi_0} Q^{-1/2 + \epsilon}.
	$$
	By the convexity bound together with bounds towards the Ramanujan-Petersson conjecture, we have
	$$
		L^{(S)}(s+1/2, \pi \otimes \chi) \ll_{\epsilon} (1+|s|)^2 C(\pi \otimes \chi)^{1/2+\epsilon}, \Re (s) = -1/2-\epsilon.
	$$
	Note that $C(\pi \otimes \chi) \ll C(\pi)C(\chi)^2$, we finally get
	$$
		l^{\chi,\sigma * h_{0,Q^{-\kappa-1}}}(\varphi) \ll_{\epsilon, \varphi_0, h_0} Q^{-\kappa /2 + \epsilon}.
	$$
	Similar argument, using Mellin inversion for $\Re s = 1/2+\epsilon$, gives
	$$
		l^{\chi,\sigma *(1- h_{0,Q^{\kappa-1}}) }(\varphi) \ll_{\epsilon, \varphi_0, h_0} Q^{-\kappa /2 + \epsilon}.
	$$
	Lemma \ref{Tronc} is proved by taking $h = h_{0,Q^{\kappa-1}} - h_{0,Q^{-\kappa-1}}$.
	
	We will need to exploit the Mellin transform of $h$ further. Since for any $h \in C_c(\R_+)$,
	$$ \mathcal{M}(h)(s) = (-1)^n \frac{\mathcal{M}(h^{(n)})(s+n)}{s(s+1)\dotsb (s+n-1) }, $$
	we have, for $h = h_{0,A}$,
	$$ \mathcal{M}(h^{(n)})(s) = A^{s-n} \mathcal{M}(h_0^{(n)})(s). $$
	For $h = h_{0,Q^{\kappa-1}} - h_{0,Q^{-\kappa-1}}$, we thus get for $n \geq 1$,
	$$ \mathcal{M}(h)(s) = (-1)^n \frac{(Q^{(\kappa -1)s}-Q^{-(\kappa +1)s})\mathcal{M}(h_0^{(n)})(s+n)}{s(s+1)\dotsb (s+n-1)}. $$
	Note that $h_0^{(n)}$ is supported in $[1,2]$, hence
	\begin{align}
		|\mathcal{M}(h)(s)| &\leq \frac{2\kappa |s|  \log Q \max(Q^{(\kappa -1) \Re (s)}, Q^{-(\kappa +1) \Re (s)}) \lVert h_0^{(n)} \rVert_{\infty} \int_1^2 t^{\Re(s)+n} d^{\times}t }{|s(s+1) \dotsb (s+n-1)|} \nonumber \\
		&\ll_{\Re(s)+n} \frac{2\kappa \log Q \lVert h_0^{(n)} \rVert_{\infty} Q^{(\kappa -1) \Re (s)}}{|(s+1) \dotsb (s+n-1)|}, \Re(s) \geq 0, \label{Mellin>0} \\
		&\ll_{\Re(s)+n} \frac{2\kappa \log Q \lVert h_0^{(n)} \rVert_{\infty} Q^{(-\kappa -1) \Re (s)}}{|(s+1) \dotsb (s+n-1)|}, \Re(s)<0.	\label{Mellin<0}
	\end{align}

	\subsection{Estimation of the Constant Contribution}
	
	Writing the Fourier expansion
	$$ \varphi_0(g) = \sum_{\alpha \in F^{\times}} W_0(a(\alpha)g), $$
	we obtain
	$$ (a(\frac{\varpi_{v_1}}{\varpi_{v_1'}}) \varphi_0 a(\frac{\varpi_{v_2}}{\varpi_{v_2'}})\overline{\varphi}_0)_N (g) = \sum_{\alpha \in F^{\times}} W_0(a(\alpha)g a(\frac{\varpi_{v_1}}{\varpi_{v_1'}}) ) \overline{W_0(a(\alpha)g a(\frac{\varpi_{v_2}}{\varpi_{v_2'}}))}. $$
	As a consequence, we get a Rankin-Selberg like equality for $\Re(s)$ large enough,
	\begin{equation}
		l^{|\cdot|^s}((a(\frac{\varpi_{v_1}}{\varpi_{v_1'}}) \varphi_0 a(\frac{\varpi_{v_2}}{\varpi_{v_2'}})\overline{\varphi}_0)_N) = \int_{\A^{\times}} W_0(a(y)a(\frac{\varpi_{v_1}}{\varpi_{v_1'}})) \overline{W_0(a(y)a(\frac{\varpi_{v_2}}{\varpi_{v_2'}}))} |y|^s d^{\times}y.
	\label{ConstantRS}
	\end{equation}
	This integral splits into product of local factors
	$$ \int_{\A^{\times}} W_0(a(y)a(\frac{\varpi_{v_1}}{\varpi_{v_1'}})) \overline{W_0(a(y)a(\frac{\varpi_{v_2}}{\varpi_{v_2'}}))} |y|^s d^{\times}y = \frac{L(s+1,\pi \times \bar{\pi})}{\zeta_F(2s+2)} \prod_{v|\infty} \int_{F_v^{\times}} |W_{0,v}(a(y))|^2 |y|_v^s d^{\times}y \cdot \prod_{v<\infty} \Sigma_v $$
	with
	$$ \Sigma_v = \frac{\zeta_v(2s+2)\int_{F_v^{\times}} W_{0,v}(a(y)a(u_v)) \overline{W_{0,v}(a(y)a(u_v'))} |y|_v^s d^{\times}y }{L(s+1,\pi_v \times \bar{\pi}_v)}. $$
	Here, $u_v,u_v'$ are suitably chosen according to $\left\{ v_1,v_1',v_2,v_2' \right\}$. For almost all $v$, $\Sigma_v$ equals $q_v^{-d_v/2}$. This identity admits a meromorphic continuation to $\C$ and is holomorphic for $\Re(s)>0$. By the convergence of $L(s, \pi \times \bar{\pi})$, we have
	$$ \frac{L(s+1,\pi \times \bar{\pi})}{\zeta_F(2s+2)} \ll_{\epsilon, \pi } 1, \text{ for } \Re (s) = \epsilon > 0. $$
	If $v$ is a ramified place of $\pi$, we can always say that the corresponding local factor is bounded by some constant depending only on $\Re(s),\pi$. So we may only consider unramified places of $\pi$. At such a place, $W_{0,v}$ is spherical and is the new vector (c.f. (\ref{Unramified})). If $\alpha_{1,v},\alpha_{2,v}$ are the Satake parameters ($ |\alpha_{1,v}\alpha_{2,v}| = 1$), then
	$$ W_{0,v} (a(\varpi_v^m)) = q_v^{-m/2}\frac{\alpha_{1,v}^{m+1}-\alpha_{2,v}^{m+1}}{\alpha_{1,v}-\alpha_{2,v}},m\geq 0, $$
	$$ W_{0,v} (a(\varpi_v^m)) = 0, m < 0. $$
	Hence the corresponding $\Sigma_v$ is explicitly computable. We should distinguish 8 cases. Denote ${\rm tr}_v = \alpha_{1,v} + \alpha_{2,v}, {\rm n}_v = \alpha_{1,v}\alpha_{2,v}$. If we write $\max(|\alpha_{1,v}|, |\alpha_{2,v}|) = q_v^{\theta_v}$, then $|{\rm tr}_v| \leq q_v^{\theta_v}+q_v^{-\theta_v} \ll q_v^{\theta_v}, |{\rm n}_v| = 1$.
	
	\textbf{Case 1:}	$u_v = \varpi_v, u_v' = 1$. We have
		\begin{eqnarray*}
			\Sigma_v &=& \frac{q_v^{-\frac{d_v}{2}}\zeta_v(2s+2)}{L(s+1,\pi_v \times \bar{\pi}_v)} \sum_{m=0}^{\infty} q_v^{-\frac{m+1}{2}} \frac{\alpha_{1,v}^{m+2}-\alpha_{2,v}^{m+2}}{\alpha_{1,v}-\alpha_{2,v}}q_v^{-\frac{m}{2}} \frac{\overline{\alpha_{1,v}}^{m+1}-\overline{\alpha_{2,v}}^{m+1}}{\overline{\alpha_{1,v}}-\overline{\alpha_{2,v}}} q_v^{-ms} \\
			&=& \frac{q_v^{-\frac{d_v+1}{2}}\zeta_v(2s+2)}{L(s+1,\pi_v \times \bar{\pi}_v) \lvert \alpha_{1,v} - \alpha_{2,v} \rvert^2} \left( \frac{\alpha_{1,v}^2 \overline{\alpha_{1,v}}}{1-\alpha_{1,v}\overline{\alpha_{1,v}}q_v^{-(s+1)}} - \frac{\alpha_{1,v}^2 \overline{\alpha_{2,v}}}{1-\alpha_{1,v}\overline{\alpha_{2,v}}q_v^{-(s+1)}} \right. \\
			& & \left. + \frac{\alpha_{2,v}^2 \overline{\alpha_{2,v}}}{1-\alpha_{2,v}\overline{\alpha_{2,v}}q_v^{-(s+1)}} - \frac{\alpha_{2,v}^2 \overline{\alpha_{1,v}}}{1-\alpha_{2,v}\overline{\alpha_{1,v}}q_v^{-(s+1)}} \right) \\
			&=& \frac{q_v^{-\frac{d_v+1}{2}}\zeta_v(2s+2)}{L(s+1,\pi_v \times \bar{\pi}_v) \lvert \alpha_{1,v} - \alpha_{2,v} \rvert^2} \left( \frac{\alpha_{1,v}^2(\overline{\alpha_{1,v}} - \overline{\alpha_{2,v}})}{\left( 1-\alpha_{1,v}\overline{\alpha_{1,v}}q_v^{-(s+1)} \right)\left( 1-\alpha_{1,v}\overline{\alpha_{2,v}}q_v^{-(s+1)} \right)} \right. \\
			& & \left. - \frac{\alpha_{2,v}^2(\overline{\alpha_{1,v}} - \overline{\alpha_{2,v}})}{\left( 1-\alpha_{2,v}\overline{\alpha_{2,v}}q_v^{-(s+1)} \right)\left( 1-\alpha_{2,v}\overline{\alpha_{1,v}}q_v^{-(s+1)} \right)} \right) \\
			&=& \frac{q_v^{-\frac{d_v+1}{2}}\zeta_v(2s+2)}{\alpha_{1,v} - \alpha_{2,v}} \left( \alpha_{1,v}^2 \left( 1-\alpha_{2,v} \overline{{\rm tr}_v} q_v^{-(s+1)} + \alpha_{2,v}^2 \overline{{\rm n}_v} q_v^{-2(s+1)} \right) \right. \\
			& & \left. - \alpha_{2,v}^2 \left( 1-\alpha_{1,v} \overline{{\rm tr}_v} q_v^{-(s+1)} + \alpha_{1,v}^2 \overline{{\rm n}_v} q_v^{-2(s+1)} \right) \right) \\
			&=& \frac{q_v^{-\frac{d_v+1}{2}}({\rm tr}_v-{\rm n}_v\overline{{\rm tr}_v}q_v^{-s-1})}{1-q_v^{-2s-2}}
		\end{eqnarray*}
		We get, for $\epsilon > 0$ small and $\Re (s) = \epsilon$,
		$$ \left\lvert \frac{q_v^{-\frac{d_v+1}{2}}({\rm tr}_v-{\rm n}_v \overline{{\rm tr}_v}q_v^{-s-1})}{1-q_v^{-2s-2}} \right\rvert \leq \frac{1+q_v^{-1-\epsilon}}{1-q_v^{-2-2\epsilon}}q_v^{-\frac{d_v+1}{2}} | {\rm tr}_v | \ll_{\epsilon} q_v^{-\frac{d_v+1}{2}} | {\rm tr}_v |. $$
		Hence we get the estimation for $\Re (s) = \epsilon$,
		\begin{equation}
			\lvert \Sigma_v \rvert \ll_{\epsilon} q_v^{-\frac{d_v+1}{2}} | {\rm tr}_v |.
		\label{LocalConstCase1}
		\end{equation}
		
	\textbf{Case 2:}	$u_v = 1, u_v' = \varpi_v$. We similarly have
		$$ \Sigma_v = \frac{q_v^{-\frac{d_v+1}{2}}(\overline{{\rm tr}_v}-\overline{{\rm n}_v}{\rm tr}_vq_v^{-s-1})}{1-q_v^{-2s-2}}, $$
		hence also the similar estimation for $\Re (s) = \epsilon$,
		\begin{equation}
			\lvert \Sigma_v \rvert \ll_{\epsilon} q_v^{-\frac{d_v+1}{2}} | {\rm tr}_v |.
		\label{LocalConstCase2}
		\end{equation}
		
	\textbf{Case 3:}	$u_v = \varpi_v^{-1}, u_v' = 1$. We similarly have
		$$ \Sigma_v = \frac{q_v^{-\frac{d_v+1}{2}-s}(\overline{{\rm tr}_v}-\overline{{\rm n}_v}{\rm tr}_vq_v^{-s-1})}{1-q_v^{-2s-2}}, $$
		hence also the similar estimation for $\Re (s) = \epsilon$,
		\begin{equation}
			\lvert \Sigma_v \rvert \ll_{\epsilon} q_v^{-\frac{d_v+1}{2}-\epsilon} | {\rm tr}_v |.
		\label{LocalConstCase3}
		\end{equation}
		
	\textbf{Case 4:}	$u_v = 1, u_v' = \varpi_v^{-1}$. We similarly have
		$$ \Sigma_v = \frac{q_v^{-\frac{d_v+1}{2}-s}({\rm tr}_v-{\rm n}_v\overline{{\rm tr}_v}q_v^{-s-1})}{1-q_v^{-2s-2}}, $$
		hence also the similar estimation for $\Re (s) = \epsilon$,
		\begin{equation}
			\lvert \Sigma_v \rvert \ll_{\epsilon} q_v^{-\frac{d_v+1}{2}-\epsilon} | {\rm tr}_v |.
		\label{LocalConstCase4}
		\end{equation}
		
	\textbf{Case 5:}	$u_v = \varpi_v, u_v' = \varpi_v^{-1}$. We similarly have
		$$ \Sigma_v = \frac{q_v^{-\frac{d_v}{2}-1-s} \left( {\rm tr}_v^2-{\rm n}_v-{\rm n}_v |{\rm tr}_v|^2 q_v^{-s-1} + {\rm n}_v |{\rm n}_v|^2 q_v^{-2(s+1)} \right)}{1-q_v^{-2s-2}}, $$
		hence also the similar estimation for $\Re (s) = \epsilon$,
		\begin{equation}
			\lvert \Sigma_v \rvert \ll_{\epsilon} q_v^{-\frac{d_v}{2}-1-\epsilon} (|{\rm tr}_v|^2+1).
		\label{LocalConstCase5}
		\end{equation}
		
	\textbf{Case 6:}	$u_v = \varpi_v^{-1}, u_v' = \varpi_v$. We similarly have
		$$ \Sigma_v = \frac{q_v^{-\frac{d_v}{2}-1-s} \left( \overline{{\rm tr}_v}^2-\overline{{\rm n}_v}-\overline{{\rm n}_v} |{\rm tr}_v|^2 q_v^{-s-1} + \overline{{\rm n}_v} |{\rm n}_v|^2 q_v^{-2(s+1)} \right)}{1-q_v^{-2s-2}}, $$
		hence also the similar estimation for $\Re (s) = \epsilon$,
		\begin{equation}
			\lvert \Sigma_v \rvert \ll_{\epsilon} q_v^{-\frac{d_v}{2}-1-\epsilon} (|{\rm tr}_v|^2+1).
		\label{LocalConstCase6}
		\end{equation}
		
	\textbf{Case 7:}	$u_v = \varpi_v, u_v' = \varpi_v$. We easily get
		$$ \Sigma_v = q_v^{-\frac{d_v}{2}+s}, $$
		hence also the similar estimation for $\Re (s) = \epsilon$,
		\begin{equation}
			\left\lvert \Sigma_v \right\rvert \leq q_v^{-\frac{d_v}{2}+\epsilon}.
		\label{LocalConstCase7}
		\end{equation}
	
	\textbf{Case 8:}	$u_v = \varpi_v^{-1}, u_v' = \varpi_v^{-1}$. We easily get
		$$ \Sigma_v = q_v^{-\frac{d_v}{2}-s}, $$
		hence also the similar estimation for $\Re (s) = \epsilon$,
		\begin{equation}
			\left\lvert \Sigma_v \right\rvert \leq q_v^{-\frac{d_v}{2}-\epsilon}.
		\label{LocalConstCase8}
		\end{equation}

	Note that at an archimedean place $v$, we have $W_{0,v}(a(y)) \in S(F_v^{\times})$, hence
	\begin{equation}
		\left\lvert \int_{F_v^{\times}} \left\lvert W_{0,v}(a(y)) \right\rvert^2 |y|^s d^{\times}y \right\rvert \ll_{\epsilon} 1, \Re(s)=\epsilon .
	\label{LocalConstArchi}
	\end{equation}
	\begin{lemma}
		We have Ramanujan conjecture on average, i.e.
		$$ \sum_{v \in I_E} | {\rm tr}_v |^2 \ll_{F, \epsilon, \pi} M_E E^{\epsilon}, \sum_{v \in I_E} | {\rm tr}_v | \ll_{F, \epsilon, \pi} M_E E^{\epsilon}. $$
	\label{RamanujanAverage}
	\end{lemma}
	The second inequality follows form the first. By the theory of Rankin-Selberg, $L(s, \pi \times \bar{\pi})$ is meromorphic and only has possible simple poles at $s=0,1$. This implies
	$$ \sum_{\substack{\alpha \text{ integral ideal of } F \\ N_F(\alpha) \leq N }} |\lambda_{\pi}(\alpha)|^2 \ll_{F, \epsilon, \pi} N^{1+\epsilon}, \forall \epsilon > 0. $$
	Here, $\lambda_{\pi}(\alpha)$ is the Hecke eigenvalues which coincides with ${\rm tr}_v$ when $\alpha$ is the prime ideal corresponding to $v$.
	
	We insert them into (\ref{ConstantRS}) and note that
	$$ S_{\text{cst}}(v_1,v_1',v_2,v_2') = \int_{\Re (s) = \epsilon } \mathcal{M}(h)(-s) l^{|\cdot|^s}((a(\frac{\varpi_{v_1}}{\varpi_{v_1'}}) \varphi_0 a(\frac{\varpi_{v_2}}{\varpi_{v_2'}})\overline{\varphi}_0)_N) \frac{ds}{2\pi i}, $$
	which with (\ref{Mellin<0}) gives, distinguishing w.r.t. types discribed in Proposition \ref{ATranslationTypes}, that for $\Re(s)=\epsilon$,
	\begin{itemize}
	\item[(1)]	In the \textbf{Type 1} of Proposition \ref{ATranslationTypes}, we use (\ref{LocalConstArchi}), (\ref{LocalConstCase1}),(\ref{LocalConstCase2}),(\ref{LocalConstCase3}),(\ref{LocalConstCase4}) to get
	$$ S_{\text{cst}}(v_1,v_1',v_2,v_2') \ll_{F, \epsilon, \pi, h_0} \kappa \log Q Q^{(1+\kappa) \epsilon } E^{-2+\epsilon} \prod_{v=v_1,v_1',v_2,v_2'} | {\rm tr}_v |. $$
	By Lemma \ref{RamanujanAverage}, the contribution of this case in (\ref{AverageCst}) is
	$$ \ll_{F, \epsilon, \pi} \kappa Q^{(2+\kappa) \epsilon } E^{-2+\epsilon} \frac{\left( \sum_{v \in I_E} |{\rm tr}_v| \right)^4}{M_E^4} \ll \kappa Q^{(2+\kappa) \epsilon } E^{-2+\epsilon}. $$
	
	\item[(2)]	In the \textbf{Type 2} of Proposition \ref{ATranslationTypes}, we use (\ref{LocalConstArchi}), (\ref{LocalConstCase7}), (\ref{LocalConstCase3}), (\ref{LocalConstCase4}) or (\ref{LocalConstArchi}), (\ref{LocalConstCase8}), (\ref{LocalConstCase1}), (\ref{LocalConstCase2}) to get
	$$ S_{\text{cst}}(v_1,v_1',v_2,v_2') \ll_{F, \epsilon, \pi, h_0} \kappa \log Q Q^{(1+\kappa) \epsilon } E^{-1+\epsilon} \prod_{v=v_1,v_1',v_2' {\rm or} v_1,v_2,v_2'} | {\rm tr}_v |. $$
	By Lemma \ref{RamanujanAverage}, the contribution of this case in (\ref{AverageCst}) is
	$$ \ll_{F, \epsilon, \pi} \kappa Q^{(2+\kappa) \epsilon } E^{-1+\epsilon} \frac{\left( \sum_{v \in I_E} |{\rm tr}_v| \right)^3}{M_E^4} \ll \kappa Q^{(2+\kappa) \epsilon } E^{-2+\epsilon}. $$
	
	\item[(3)]	In the \textbf{Type 3} of Proposition \ref{ATranslationTypes}, we use (\ref{LocalConstArchi}), (\ref{LocalConstCase5}), (\ref{LocalConstCase2}), (\ref{LocalConstCase3}) or (\ref{LocalConstArchi}), (\ref{LocalConstCase6}), (\ref{LocalConstCase1}), (\ref{LocalConstCase4}) to get
	\begin{eqnarray*}
		S_{\text{cst}}(v_1,v_1',v_2,v_2') &\ll_{F, \epsilon, \pi, h_0}& \kappa \log Q Q^{(1+\kappa) \epsilon } E^{-1+\epsilon} \\
		& & \cdot \left( \prod_{v=v_1,v_1',v_2,v_2'} | {\rm tr}_v | + \prod_{v=v_1',v_2 {\rm or} v_1,v_2'} | {\rm tr}_v | \right).
	\end{eqnarray*}
	By Lemma \ref{RamanujanAverage}, the contribution of this case in (\ref{AverageCst}) is
	$$ \ll_{F, \epsilon, \pi} \kappa Q^{(2+\kappa) \epsilon } E^{-1+\epsilon} \frac{\left( \sum_{v \in I_E} |{\rm tr}_v| \right)^2\left( \sum_{v \in I_E} |{\rm tr}_v|^2+1 \right)}{M_E^4} \ll \kappa Q^{(2+\kappa) \epsilon } E^{-2+\epsilon}. $$
	
	\item[(4)]	In the \textbf{Type 4} of Proposition \ref{ATranslationTypes}, we use (\ref{LocalConstArchi}), (\ref{LocalConstCase2}), (\ref{LocalConstCase4}), or (\ref{LocalConstArchi}), (\ref{LocalConstCase1}), (\ref{LocalConstCase3}) to get
	$$ S_{\text{cst}}(v_1,v_1',v_2,v_2') \ll_{F, \epsilon, \pi, h_0} \kappa \log Q Q^{(1+\kappa) \epsilon } E^{-1+\epsilon} \prod_{v=v_2,v_2' {\rm or} v_1,v_1'} | {\rm tr}_v |. $$
	By Lemma \ref{RamanujanAverage}, the contribution of this case in (\ref{AverageCst}) is
	$$ \ll_{F, \epsilon, \pi} \kappa Q^{(2+\kappa) \epsilon } E^{-1+\epsilon} \frac{M_E \left( \sum_{v \in I_E} |{\rm tr}_v| \right)^2}{M_E^4} \ll \kappa Q^{(2+\kappa) \epsilon } E^{-2+\epsilon}. $$
	
	\item[(5)]	In the \textbf{Type 5} of Proposition \ref{ATranslationTypes}, we use (\ref{LocalConstArchi}), (\ref{LocalConstCase7}), (\ref{LocalConstCase8}) to get
	$$ S_{\text{cst}}(v_1,v_1',v_2,v_2') \ll_{F, \epsilon, \pi, h_0} \kappa \log Q Q^{(1+\kappa) \epsilon } E^{\epsilon}. $$
	The contribution of this case in (\ref{AverageCst}) is
	$$ \ll_{F, \epsilon, \pi} \kappa Q^{(2+\kappa) \epsilon } E^{\epsilon} \frac{M_E^2}{M_E^4} \ll \kappa Q^{(2+\kappa) \epsilon } E^{-2+\epsilon}. $$
	
	\item[(6)]	In the \textbf{Type 6} of Proposition \ref{ATranslationTypes}, we use (\ref{LocalConstArchi}), (\ref{LocalConstCase5}), (\ref{LocalConstCase6}) to get
	$$ S_{\text{cst}}(v_1,v_1',v_2,v_2') \ll_{F, \epsilon, \pi, h_0} \kappa \log Q Q^{(1+\kappa) \epsilon } E^{-2+\epsilon} \prod_{v=v_1,v_2} (|{\rm tr}_v|^2+1). $$
	By Lemma \ref{RamanujanAverage}, the contribution of this case in (\ref{AverageCst}) is
	$$ \ll_{F, \epsilon, \pi} \kappa Q^{(2+\kappa) \epsilon } E^{-2+\epsilon} \frac{\left( \sum_{v \in I_E} |{\rm tr}_v|^2+1 \right)^2}{M_E^4} \ll \kappa Q^{(2+\kappa) \epsilon } E^{-4+\epsilon}. $$
	
	\item[(7)]	In the \textbf{Type 7} of Proposition \ref{ATranslationTypes}, we use (\ref{LocalConstArchi}) to get
	$$ S_{\text{cst}}(v_1,v_1',v_2,v_2') \ll_{F, \epsilon, \pi, h_0} \kappa \log Q Q^{(1+\kappa) \epsilon }. $$
	The contribution of this case in (\ref{AverageCst}) is
	$$ \ll_{F, \epsilon, \pi} \kappa Q^{(2+\kappa) \epsilon } \frac{M_E^2}{M_E^4} \ll \kappa Q^{(2+\kappa) \epsilon } E^{-2+\epsilon}. $$
	
	\item[(8)]	In the \textbf{Type 8} of Proposition \ref{ATranslationTypes}, we use (\ref{LocalConstArchi}), (\ref{LocalConstCase2}), (\ref{LocalConstCase4}) or (\ref{LocalConstArchi}), (\ref{LocalConstCase1}), (\ref{LocalConstCase3}) to get
	$$ S_{\text{cst}}(v_1,v_1',v_2,v_2') \ll_{F, \epsilon, \pi, h_0} \kappa \log Q Q^{(1+\kappa) \epsilon } E^{-1+\epsilon} \prod_{v=v_2,v_2' {\rm or} v_1,v_1'} |{\rm tr}_v|. $$
	By Lemma \ref{RamanujanAverage}, the contribution of this case in (\ref{AverageCst}) is
	$$ \ll_{F, \epsilon, \pi} \kappa Q^{(2+\kappa) \epsilon } E^{-1+\epsilon} \frac{\left( \sum_{v \in I_E} |{\rm tr}_v| \right)^2}{M_E^4} \ll \kappa Q^{(2+\kappa) \epsilon } E^{-3+\epsilon}. $$
	
	\item[(9)]	In the \textbf{Type 9} of Proposition \ref{ATranslationTypes}, we use (\ref{LocalConstArchi}) to get
	$$ S_{\text{cst}}(v_1,v_1',v_2,v_2') \ll_{F, \epsilon, \pi, h_0} \kappa \log Q Q^{(1+\kappa) \epsilon }. $$
	The contribution of this case in (\ref{AverageCst}) is
	$$ \ll_{F, \epsilon, \pi} \kappa Q^{(2+\kappa) \epsilon } \frac{M_E}{M_E^4} \ll \kappa Q^{(2+\kappa) \epsilon } E^{-3+\epsilon}. $$
	
	\end{itemize}
	The proof of Lemma \ref{ConstantL} is completed.

	\subsection{Estimation of the Cuspidal Constribution}
	\label{BaseChoice}
The goal of this section is to establish Lemma \ref{CuspL}. Recall that we are reduced to estimating 
$$ S_{\text{cusp}}(v_1,v_1',v_2,v_2') =   \sum_{\pi' \text{cuspidal}} l^h( n(T) P_{\pi'}(a(\frac{\varpi_{v_1}}{\varpi_{v_1'}}) \varphi_0 a(\frac{\varpi_{v_2}}{\varpi_{v_2'}}) \overline{\varphi}_0 )). $$
The projector $P_{\pi'}$ is realized by the choice of a basis of $\pi'$, denoted by $\mathcal{B}(\pi'; v_1,v_1',v_2,v_2')$. It is determined by the choices of local basis of $\pi_v'$, denoted by $\mathcal{B}_v(\pi'; v_1,v_1',v_2,v_2')$. When there is no confusion, we may write them shortly as $\mathcal{B}$ resp. $\mathcal{B}_v$. They are related with each other by
$$ \mathcal{B} = \Pi_v' \mathcal{B}_v, e \leftrightarrow (W_{e,v})_v. $$
Here, $W_{e,v}$ is the component at $v$ of $e$ in the Kirillov model. We may also write it as $e_v$ if there is no confusion. According to Remark \ref{ArchiBasis}, we only need to choose $\mathcal{B}_v$ for $v < \infty$.
\begin{definition}
	Denote, for any subgroup $H \subset G(F_v)$ and $g \in G(F_v)$, $H^g = gHg^{-1}$. Then the Harish-Chandra's function $\Xi_v^{g_0}$ associated to the Borel subgroup $B(F_v)^{g_0}$ is given by, with notations in Section \ref{DefOfXi}
	$$ \Xi_v^{g_0}(g) = \Xi_v(g_0^{-1}gg_0). $$
\end{definition}
\begin{definition}
	Suppose $v(\pi')=m$. For any integer $n$, recall that the space of $K_v^0[n]$-invariant vectors of $\pi_v'$ is of dimension $\max(n-m+1,0)$. A \textbf{standard basis} of level $n$ consists of, for each integer $l$ such that $m \leq l \leq n$, a vector invariant by $K_v^0[l]$ and orthogonal to all the vectors invariant by $K_v^0[l-1]$, and vectors orthogonal to the space of $K_v^0[n]$-invariant vectors. A \textbf{nice basis} of level $n$ w.r.t. $g\in G_v$ consists of the $g$ translates of the vectors of a standard basis of level $n$. Define \textbf{the maximal compact subgroup} $K_v^*$ of $G_v$ \textbf{associated with the above nice basis} to be $$K_v^* = K_v^g. $$
	If $\mathcal{B}_v$ is a standard or nice basis of level $n$, we write $\mathcal{B}_v^*$ to be the elements in $\mathcal{B}_v$ invariant by $K_v^0[n]$ or its corresponding translate. We also call the basis as in Remark \ref{ArchiBasis} standard. At an infinite place, we define $\mathcal{B}_v^* = \mathcal{B}_v$. We write
	$$ \mathcal{B}^* = \Pi_v' \mathcal{B}_v^*. $$
\end{definition}
\begin{remark}
	Note that, if $\mathcal{B}_v$ is a standard basis of level $n$, then $\mathcal{B}_v^*$ is just the set of $v_l(\pi')$, defined in Section \ref{BranchingLawSection}, with $v(\pi')=m \leq l \leq n$.
\end{remark}
We choose $\mathcal{B}_v$ and $K_v^*$ explicitly as follows:

\textbf{Case 1:}

\begin{tabular}{| m{3.6cm} || m{3cm} | c |}
	\hline
	Position of $v$ & $a(\frac{\varpi_{v_1}}{\varpi_{v_1'}}) \varphi_0 a(\frac{\varpi_{v_2}}{\varpi_{v_2'}}) \overline{\varphi}_0$ & $\mathcal{B}_v$ \\ \hline
	$v \notin \{v_1,v_1',v_2,v_2'\}$ or $v=v_1=v_1'=v_2=v_2'$ & $K_v^0[v(\varphi_0)]$-invariant & standard of level $v(\varphi_0)$ \\ \hline
	$v=v_1=v_2 \notin \{ v_1',v_2' \}$ & $K_v^0[v(\varphi_0)]^{a(\varpi_v)}$-invariant & nice of level $v(\varphi_0)$ w.r.t. $a(\varpi_v)$ \\ \hline
	$v=v_1'=v_2' \notin \{ v_1, v_2 \}$ & $K_v^0[v(\varphi_0)]^{a(\varpi_v^{-1})}$-invariant & nice of level $v(\varphi_0)$ w.r.t. $a(\varpi_v^{-1})$ \\ \hline
\end{tabular}

\textbf{Case 2:}

\begin{tabular}{| m{3.8cm} || m{3cm} | c |}
	\hline
	Position of $v$ & $a(\frac{\varpi_{v_1}}{\varpi_{v_1'}}) \varphi_0 a(\frac{\varpi_{v_2}}{\varpi_{v_2'}}) \overline{\varphi}_0$ & $\mathcal{B}_v$ \\ \hline
	$v=v_1 \notin \{ v_1',v_2,v_2' \}$ or $v=v_1=v_2=v_2' \neq v_1'$ or $v=v_2 \notin \{ v_1,v_1',v_2' \}$ or $v=v_2=v_1=v_1' \neq v_2'$ & $K_v^0[v(\varphi_0)+1]^{a(\varpi_v)}$-invariant & nice of level $v(\varphi_0)+1$ w.r.t. $a(\varpi_v)$ \\ \hline
	$v=v_1' \notin \{ v_1,v_2,v_2' \}$ or $v=v_1'=v_2=v_2' \neq v_1$ or $v=v_2' \notin \{ v_1,v_1',v_2 \}$ or $v=v_2'=v_1=v_1' \neq v_2$ & $K_v^0[v(\varphi_0)+1]$-invariant & standard of level $v(\varphi_0)+1$ \\ \hline
\end{tabular}

\textbf{Case 3:}

\begin{tabular}{| m{3.8cm} || m{3cm} | c |}
	\hline
	Position of $v$ & $a(\frac{\varpi_{v_1}}{\varpi_{v_1'}}) \varphi_0 a(\frac{\varpi_{v_2}}{\varpi_{v_2'}}) \overline{\varphi}_0$ & $\mathcal{B}_v$ \\ \hline
	$v=v_1=v_2' \notin \{ v_1', v_2 \}$ or $v=v_2=v_1' \notin \{ v_1, v_2' \}$ & $K_v^0[v(\varphi_0)+2]^{a(\varpi_v)}$-invariant & nice of level $v(\varphi_0)+2$ w.r.t. $a(\varpi_v)$ \\ \hline
\end{tabular}

Then we rewrite
\begin{equation}
	S_{{\rm cusp}}(v_1,v_1',v_2,v_2') =  \sum_{\pi'} \sum_{e\in \mathcal{B}^* }  C(\varphi_0, e; v_1,v_1',v_2,v_2') l^h(n(T)e),
\label{CuspSum}
\end{equation}
with
$$ C(\varphi_0, e; v_1,v_1',v_2,v_2') = \langle a(\frac{\varpi_{v_1}}{\varpi_{v_1'}}) \varphi_0 a(\frac{\varpi_{v_2}}{\varpi_{v_2'}}) \overline{\varphi}_0, e \rangle. $$
We have
\begin{equation}
	l^h(n(T)e) = \int_{\Re (s) = 0} \mathcal{M}(h)(-s) l^{|\cdot|^s}(n(T)e) \frac{ds}{2\pi i},
\label{CuspMellinInverse}
\end{equation}
and since the vector $e$ is a pure tensor, we have
$$ l^{|\cdot |^s}(n(T)e) = L(s+1/2, \pi') \prod_{v | \infty} l^{|\cdot|^s}(n(T_v)W_{e,v}) \prod_{v<\infty} \frac{l^{|\cdot|^s}(n(T_v)W_{e,v})}{L(s+1/2, \pi'_v)}. $$
\begin{lemma}
	We have for any $\epsilon > 0$,
	$$ |l^{|\cdot |^s}(n(T)e)| \ll_{\epsilon, \varphi_0} |L(s+1/2,\pi')| |T|^{-1/2+\theta+\epsilon} \lambda_{e,\infty}^{1/2+\epsilon}, s \in i\R. $$
\label{MBound}
\end{lemma}
To prove Lemma \ref{MBound}, we shall estimate the local terms case by case. This is technical and will be given in the subsequent subsections. In fact, Lemma \ref{MBound} will be a consequence of Corollary \ref{Case1LocalEst}, Lemma \ref{Case2LocalEst}, \ref{Case3LocalEst}, \ref{Case4LocalEst}, as well as Lemma \ref{NormIden} and Remark \ref{LocalBoundRSLFactorCusp} (with $\lVert e \rVert_{X(F)} = 1$). Thus we get
\begin{equation}
	|S_{{\rm cusp}}(v_1,v_1',v_2,v_2')| \ll_{\epsilon, \varphi_0} |T|^{-1/2+\theta+\epsilon} \int_{\Re (s) = 0} \mathcal{M}(h)(-s) S_{{\rm cusp}}^*(s; v_1,v_1',v_2,v_2') \frac{ds}{2\pi i},
\label{SBound}
\end{equation}
with
$$ S_{{\rm cusp}}^*(s; v_1,v_1',v_2,v_2') = \sum_{\pi'} \sum_{e\in \mathcal{B}^*} \lambda_{e,\infty}^{1/2+\infty} |L(s+1/2,\pi')| \left| C(\varphi_0, e; v_1,v_1',v_2,v_2') \right|. $$
\begin{theorem}
	Give $(n_v)_{v < \infty}$ with $n_v \in \N, n_v = 0$ for all but finitely many $v$. For any cuspidal representation $\pi'$, let $\pi_{\infty} = \otimes_{v | \infty} \pi_v'$ be the infinite part of $\pi'$, on which the Casimir element of $Z_{\infty} \backslash G_{\infty} = \prod_{v | \infty} Z_v \backslash G(F_v)$ acts by multiplication by $\lambda_{\pi',\infty}^*$. Then there is some constant $A,B > 0$ such that for $s \in i\R$
	$$ \sum_{\pi': v(\pi') \leq n_v, v < \infty} |L(s+1/2,\pi')|^4 (\lambda_{\pi',\infty}^*)^{-A} \ll_{\epsilon} (1+|s|)^B \left( \prod_{v < \infty} q_v^{n_v} \right)^{1+\epsilon}. $$
\label{LindelofAverage4th}
\end{theorem}
The above is a consequence of the main theorem of \cite{MV} without amplification. We sketch the proof as follows. Write $\mathcal{C} = \prod_{v < \infty} q_v^{n_v}$. We construct some (normalized) Eisenstein series $\varphi_1 \in \pi_1=\pi(1,1), \varphi_2 \in \pi_2=\pi(|\cdot|^s, |\cdot^{-s}|)$ with $v(\varphi_1) = 0, v(\varphi_2) = n_v, \forall v < \infty$, and prove the existence of some (normalized) $\varphi_3 \in \pi'^{\infty}$ such that for $\epsilon > 0$ and some constant $C > 0$ the triple product
$$ I(\varphi_1,\varphi_2,\varphi_3) = \int_{ZG(F) \backslash G(\A)} \varphi_1(g)\varphi_2(g)\varphi_3(g) \gg_{\epsilon} |L(1/2, \pi_1 \times \pi_2 \times \pi')| C_{\infty}(\pi_1 \times \pi_2 \times \pi')^{-C} \mathcal{C}^{-1-\epsilon}. $$
On the other hand, the triple product is just a coefficient of the projection of $\varphi_1 \varphi_2$ onto the space of $\pi'$, hence
$$ \sum_{\pi'} |I(\varphi_1,\varphi_2,\varphi_3)|^2 \leq \langle \varphi_1 \varphi_2, \varphi_1 \varphi_2 \rangle_{{\rm reg}} = \langle \varphi_1 \overline{\varphi}_1, \varphi_2 \overline{\varphi}_2 \rangle_{{\rm reg}}, $$
where $\langle \cdot, \cdot \rangle_{{\rm reg}}$ is some regularized inner product. $\varphi_1 \overline{\varphi}_1$ being spherical at all finite place, the right hand side is bounded by $O(1)$. We conclude the theorem by noticing that $L(1/2, \pi_1 \times \pi_2 \times \pi') = L(s+1/2,\pi')^2 L(-s+1/2,\pi')^2$ and $C_{\infty}(\pi_1 \times \pi_2 \times \pi')^{-C} \ll (\lambda_{\pi',\infty}^*)^{-A} (1+|s|)^{-B}$.
\begin{corollary}
	With notation as in Theorem \ref{LindelofAverage4th}, there is some constant $A,B > 0$ such that
$$ \sum_{\pi': v(\pi') \leq n_v, v < \infty} |L(s+1/2,\pi')|^2 (\lambda_{\pi',\infty}^*)^{-A} \ll_{\epsilon} (1+|s|)^B \left( \prod_{v < \infty} q_v^{n_v} \right)^{1+\epsilon}. $$
\label{LindelofAverage2nd}
\end{corollary}
\proof This is just a usual application of Cauchy-Schwarz inequality combined with Theorem \ref{WeylLaw}. \endproof
We apply Cauchy-Schwarz to get, for some constant $A' > 0$ large enough,
\begin{align}
	S_{{\rm cusp}}^*(s; v_1,v_1',v_2,v_2') &\leq \left( \sum_{\pi'} \sum_{e \in \mathcal{B}^*} \lambda_{e, \infty}^{1/2-A'+\epsilon} |L(s+1/2,\pi')|^2 \right)^{1/2} \lVert \Delta_{\infty}^{A'} \left( a(\frac{\varpi_{v_1}}{\varpi_{v_1'}}) \varphi_0 a(\frac{\varpi_{v_2}}{\varpi_{v_2'}}) \overline{\varphi}_0 \right) \rVert \nonumber \\
	&\ll_{\varphi_0} \left( \sum_{\pi': v(\pi') \leq n_v, v < \infty} |L(s+1/2,\pi')|^2 (\lambda_{\pi',\infty}^*)^{-A} \right)^{1/2} (\prod_{n_v \neq 0} n_v)^{1/2} \nonumber \\
	&\ll_{\epsilon} (1+|s|)^{B/2} \left( \prod_{v < \infty} q_v^{n_v} \right)^{1/2+\epsilon}, \label{S^*Bound}
\end{align}
where $n_v$ is just the level of $\mathcal{B}_v$ chosen for the spectral decomposition. Distinguishing the $9$ types described in Proposition \ref{ATranslationTypes}, we easily see
\begin{equation}
	\prod_{v < \infty} q_v^{n_v} \ll_{\varphi_0} E^4,
\label{AmplificationBound}
\end{equation}
and \textbf{Type 1} contributes $E^4$. Inserting (\ref{AmplificationBound}) and (\ref{S^*Bound}) into (\ref{SBound}), using (\ref{Mellin>0}) we get Lemma \ref{CuspL}.

We turn to the proof of Lemma \ref{MBound}.

		\subsubsection{\textbf{At $v$ such that $T_v \neq 0$}}
				
	In this case, $\mathcal{B}_v$ is given by the first case of \textbf{Case 1}, hence is standard. Note that for $s \in i\R$,
	$$ |l^{|\cdot |^s}(n(T_v)W_{e,v})|^2 = \int_{F_v^{\times}} \langle n(-T_v)a(y)n(T_v)W_{e,v},W_{e,v} \rangle |y|^s d^{\times}y. $$
	By Theorem \ref{MatrixCoeffDecay}, we get,
	\begin{equation}
		|l^{|\cdot |^s}(n(T_v)W_{e,v})|^2 \leq A_v(\epsilon) \dim (K_v e_v) \lVert W_{e,v} \rVert^2 \int_{F_v^{\times}} \Xi_v(n(-T_v)a(y)n(T_v))^{1-2\theta-\epsilon} d^{\times}y.
	\label{LbyXi}
	\end{equation}
	\begin{lemma}
		For any $\epsilon > 0$, we have
		$$ |l^{|\cdot |^s}(n(T_v)W_{e,v})| \ll_{\epsilon, \theta} |T_v|_v^{-1/2+\theta+\epsilon} \dim(K_v e_v)^{1/2} \lVert W_{e,v} \rVert, s \in i\R. $$
	\label{Case1Lemma}
	\end{lemma}
	\begin{corollary}
		There exists a constant $C(\theta, \epsilon)$ depending only on $\theta$ and $\epsilon$ such that:
		
		If $ v|\infty $, then we have
		$$ |l^{|\cdot |^s}(n(T_v)W_{e,v})| \leq C(\theta, \epsilon) \lambda_{e,v}^{1/2+\epsilon} |T_v|_v^{-1/2+\theta+\epsilon} \lVert W_{e,v} \rVert, s \in i\R. $$
		
		If $ v < \infty $, then 
		$$ \left| \frac{l^{|\cdot |^s}(n(T_v)W_{e,v})}{L(s+1/2,\pi_v')} \right| \leq C(\theta, \epsilon) |T_v|_v^{-1/2+\theta+\epsilon} q_v^{v(\varphi_0)/2+\epsilon} \frac{\zeta_v(2)^{1/2}\lVert W_{e,v} \rVert}{\sqrt{L(1, \pi_v' \times \bar{\pi}_v')}}, s \in i\R. $$
	\label{Case1LocalEst}
	\end{corollary}
	\begin{remark}
		By reducing $\epsilon$ to $\epsilon /2$, $C(\theta, \epsilon)=1$ is admissible for all $v < \infty$ outside a finite set of places depending only on $\theta$ and $\epsilon$. It will always be such case whenever $C(\theta, \epsilon)$ appears after. This ensures us that the product of $C(\theta, \epsilon)$ over all places is still bounded by some constant depending only on $\theta$ and $\epsilon$.
	\end{remark}
	Note that if $v | \infty$, $\dim(K_v e_v), C(\pi_v') \ll \lambda_{e,v}$; and if $v < \infty$, $e_v$ is $K_v^0[v(\varphi_0)]$ invariant by the choice of $\mathcal{B}_v^*$, $v(\pi') \leq v(\varphi_0)$. We deduce the corollary from the lemma by taking into account Remark \ref{LocalBoundRSLFactorCusp} and
	$$ [K_v : K_v^0[v(\varphi_0)] ] \ll q_v^{v(\varphi_0)}. $$
	
	We now prove Lemma \ref{Case1Lemma} place by place.
	
	\textbf{At a Real Place : $F_v = \R$}
	
	Recall that the (bi-$K_v$-invariant, $\begin{pmatrix} -1 & 0\\ 0 & 1 \end{pmatrix}$-invariant) Harish-Chandra's function as in \cite{CU}, 5.2.2 is given by some Legendre function as:
 	$$ \Xi_v(\begin{pmatrix} e^{r/2} & 0\\ 0 & e^{-r/2} \end{pmatrix}) = \mathfrak{P}_{-1/2}(\cosh r),r>0. $$
	For some absolute constants $\alpha, \beta > 0$, we have
	$$ \mathfrak{P}_{-1/2}(\cosh r) \leq e^{-r/2}(\alpha + \beta r).$$
	We make a change of variable $t = \frac{y+y^{-1}}{2}$ and get
	\begin{align*}
		&\ \int_{\R^{\times}} \Xi_v(n(-T_v)a(y)n(T_v))^{1-2\theta} d^{\times}y \\
		&\leq 2(1+T_v^2)^{-\frac{1-2\theta}{2}} (1+\log(1+T_v^2))^{1-2\theta} \int_1^{\infty} (t-1)^{-1/2+\theta}(\alpha' + \beta \log t)^{1-2\theta} \\
		&\ + t^{-1/2+\theta}(\alpha' + \beta \log (t+1))^{1-2\theta} \frac{dt}{\sqrt{t^2-1}} \\
		&\ll_{\theta} (1+T_v^2)^{-\frac{1-2\theta}{2}} (1+\log(1+T_v^2))^{1-2\theta}.
	\end{align*}
	We get the lemma at $v$ using (\ref{LbyXi}).
	
	\textbf{At a Complex Place : $F_v = \C$}
	
	The Harish-Chandra's function as in \cite{CU} 5.2.1 is given by:
	$$ \Xi_v(\begin{pmatrix} t & 0\\ 0 & t^{-1} \end{pmatrix}) = \frac{2 \log t}{t-t^{-1}}, t>0. $$
	When we evaluate it at $n(-T_v)a(y)n(T_v)$, the corresponding $t$ satisfies
	$$ t^2 + t^{-2} = |y| + |y|^{-1} + \frac{|T_v|^2 |y-1|^2}{|y|}. $$
	This expression being invariant by the change of variable $y \mapsto y^{-1}$, we get, with the change of variable $r = \frac{|y|+|y|^{-1}}{2}$,
	\begin{align*}
		&\ \int_{\C^{\times}} \Xi_v(n(-T_v)a(y)n(T_v))^{1-2\theta} d^{\times}y = 2\int_{|y|>1} (\frac{2 \log t}{t-t^{-1}})^{1-2\theta} d^{\times}y \\
		&\leq 2 (2(1+|T_v|^2))^{-\frac{1-2\theta}{2}} (\log 2(1+|T_v|^2))^{1-2\theta} \cdot 2\pi \int_1^{\infty} (\frac{1+\frac{\log (r+1)}{\log 2}}{\sqrt{r-1}})^{1-2\theta} \frac{dr}{\sqrt{r^2-1}} \\
		&\ll_{\theta} (1+|T_v|^2)^{-\frac{1-2\theta}{2}} (1+\log(1+|T_v|^2))^{1-2\theta}.
	\end{align*}
	We get the lemma at $v$ using (\ref{LbyXi}).
	
	\textbf{At a Non Archimedean Place}
	
	The values of the Harish-Chandra function associated with the standard Borel subgroup can be inferred from the Macdonald formula, i.e. Theorem 4.6.6 of \cite{B}, by letting $\alpha_1 \rightarrow 1, \alpha_2 = 1$,
	$$
		\Xi_v(n) = \Xi_v( \begin{pmatrix} \varpi_v^n & 0\\ 0 & 1 \end{pmatrix} ) = q_v^{-n/2} +  n q_v^{-n/2} \frac{1-q_v^{-1}}{1+q_v^{-1}}, n \geq 0.
	$$
	We apply (42) of \cite{CU} to the torus $\mathbb{T} = n(-T_v) A_v n(T_v)$. More precisely, using the notations as in \cite{CU}, we calculate in our situation $\mathbb{T}_c = n(-T_v) a(\vO_v^{\times}) n(T_v)$, $\mathbb{T}_1 = \mathbb{T} \cap K_v = n(-T_v) a(1+\varpi_v^d \vO_v) n(T_v)$, where $d=\max(0,-v(T_v))$ with convention $1+\varpi_v^0 \vO_v = \vO_v^{\times}$. Hence $\mathbb{T}_c / \mathbb{T}_1 \simeq \vO_v^{\times} / (1+\varpi_v^d \vO_v)$ and $\left\lvert \mathbb{T}_c / \mathbb{T}_1 \right\rvert = (q_v-1)q_v^d$ with convention $(q_v-1)q_v^0=1$. Following the proof of Lemma 5.6 of \cite{CU}, we see $\delta = d(p, \Gamma) = d$ with $x_0 = \text{proj}(p, \Gamma) = n(-T_v) a(\varpi_v^{-d})$, where $\Gamma$ is the unique geodesic in the Tits building of $G(F_v)$ fixed by $\mathbb{T}_c$ and $p$ is the similitude class of the standard lattice $\Lambda_0 = \vO_v \oplus \vO_v$. We identify $\tau \in \vO_v^{\times} / (1+\varpi_v^d \vO_v)$ with its corresponding element in $\mathbb{T}_c / \mathbb{T}_1$. The only terms remaining to determine is $\epsilon(\tau)$. By definition, we know
	$$ 2\epsilon(\tau) = d(n(-T_v) a(\tau) n(T_v).p, p) = d(\begin{pmatrix} \tau & T_v(\tau -1) \\ 0 & 1 \end{pmatrix}.p, p). $$
	But it is easy to see that $\varpi_v^{d-v(\tau-1)} \begin{pmatrix} \tau & T_v(\tau -1) \\ 0 & 1 \end{pmatrix}. \Lambda_0$ is the ``smallest'' lattice in the similitude class of lattices of $n(-T_v) a(\tau) n(T_v).p$ contained in $\Lambda_0$, hence
	$$ d(\begin{pmatrix} \tau & T_v(\tau -1) \\ 0 & 1 \end{pmatrix}.p, p) = v\left( \det \left( \varpi_v^{d-v(\tau-1)} \begin{pmatrix} \tau & T_v(\tau -1) \\ 0 & 1 \end{pmatrix} \right) \right) = 2(d-v(\tau -1)). $$
	The number of $\tau$'s such that $v(\tau -1)=n$ is $(q_v-1)q_v^{d-n-1}$ if $1 \leq n \leq d-1$; $(q_v -2)q_v^{d-1}$ if $n=0$; $q$ if $n=d$. We can therefore calculate and bound the local integral as,
	\begin{align*}
		&\ q_v^{d_v/2} \int_{F_v^{\times}} \Xi_v(n(-T_v)a(y)n(T_v))^{1-2\theta} d^{\times}y \\
		&= 2 \sum_{n> 2 d} \Xi_v(n)^{1-2\theta} + \sum_{n=1}^{d-1} \frac{ q_v^{d-n} - q_v^{d-n-1} }{q_v^d - q_v^{d-1}} \Xi_v(2(d-n))^{1-2\theta} \\
		&\ + \frac{1}{q_v^d - q_v^{d-1}} \Xi_v(0)^{1-2\theta} + \frac{q_v^d - 2q_v^{d-1}}{q_v^d - q_v^{d-1}} \Xi_v(2d)^{1-2\theta} \\
		&\ll C(\theta) \max(1,|T_v|)^{-(1-2\theta)} (1+ \max(1, \log |T_v|))^{2-2\theta}.
	\end{align*}
	We get the lemma at $v$ by using (\ref{LbyXi}) and conclude the lemma.
	
	We record the following estimation: for some constant $C'(\theta)$ depending only on $\theta$,
	\begin{equation}
		q_v^{d_v/2} \int_{F_v^{\times}} \Xi_v(a(y))^{1-2\theta} d^{\times}y \leq 2\sum_{n>0} (n+1)q_v^{-n(1/2-\theta)} + 1 \leq C'(\theta).
	\label{XiBoundFinPlace}
	\end{equation}

				\subsubsection{\textbf{At $v$ such that $T_v = 0$, $\pi_v$ ramified}}
	
	The number of such places is finite and depends only on $\pi$. $\mathcal{B}_v$ is standard. Since the local vectors concerned are classical vectors, we shall use Proposition \ref{MCDwoD} instead of Theorem \ref{MatrixCoeffDecay}, combined with (\ref{XiBoundFinPlace}) to get an inequality similar to (\ref{LbyXi}):
	\begin{equation}
		|l^{|\cdot |^s}(W_{e,v})|^2 \leq A_v(\epsilon) \lVert W_{e,v} \rVert^2 \int_{F_v^{\times}} \Xi_v(a(y))^{1-2\theta-\epsilon} d^{\times}y \leq C(\theta, \epsilon) \lVert W_{e,v} \rVert^2.
	\label{LbyXiwoD}
	\end{equation}
	\begin{lemma}
		For any $\epsilon > 0$, there is a constant $C(\theta, \epsilon) $ such that
		$$ \left| \frac{l^{|\cdot |^s}(n(T_v)W_{e,v})}{L(s+1/2,\pi_v')} \right| \leq C(\theta, \epsilon) \frac{\zeta_v(2)^{1/2}\lVert W_{e,v} \rVert}{\sqrt{L(1, \pi_v' \times \bar{\pi}_v')}}, s \in i\R. $$
	\label{Case2LocalEst}
	\end{lemma}

		\subsubsection{\textbf{At $v$ such that $T_v = 0$, $\pi_v$ unramified, $v \in \left\{ v_1,v_1',v_2,v_2' \right\}$}}
		
	The number of possible places is at most $4$ and $v(\pi')\leq 2$. Since the vectors concerned are $a(\cdot)$-translates of classical vectors and $s\in i\R$, (\ref{LbyXiwoD}) still applies and gives
	\begin{lemma}
		For any $\epsilon > 0$, there is a constant $C(\theta, \epsilon)$ such that
		$$ \left\lvert \frac{l^{|\cdot |^s}(W_{e,v})}{L(s+1/2, \pi_v')} \right\rvert \leq C(\theta,\epsilon) \frac{\zeta_v(2)^{1/2}\lVert W_{e,v} \rVert}{\sqrt{L(1, \pi_v' \times \bar{\pi}_v')}}, s \in i\R. $$
	\label{Case3LocalEst}
	\end{lemma}

		\subsubsection{\textbf{At $v$ such that $T_v = 0$, $\pi_v$ unramified, $v\notin \left\{ v_1,v_1',v_2,v_2' \right\}$}}
	
	In this case $e_v$ is spherical and we have
	\begin{lemma}
		For $v\notin \left\{ v_1,v_1',v_2,v_2' \right\}$, we have
		$$ \left\lvert \frac{l^{|\cdot |^s}(W_{e,v})}{L(s+1/2, \pi_v')} \right\rvert = \frac{\zeta_v(2)^{1/2} \lVert W_{e,v} \rVert}{\sqrt{L(1,\pi_v' \times \bar{\pi}_v')}} = |W_{e,v}(1)|, s \in i\R. $$
	\label{Case4LocalEst}
	\end{lemma}
	Note that almost all $v$ fall into this case.

	\subsection{Estimation of the Eisenstein Contribution}
	The goal of this section is to establish Lemma \ref{EisL}. We rewrite
	\begin{align}
	S_{\text{Eis}}(v_1,v_1',v_2,v_2') &= \sum_{\xi \in \widehat{F^{\times} \backslash \A^{(1)}}} \int_{-\infty}^{\infty} \sum_{\Phi \in \mathcal{B}(\pi(\xi, \xi^{-1}))} C(\varphi_0, \Phi; v_1,v_1',v_2,v_2') \cdot \nonumber \\
	&\  l^h(n(T)(E(\Phi,i\tau)-E_N(\Phi,i\tau))) \frac{d\tau}{4\pi}. \label{EisSum}
	\end{align}
with
	$$ C(\varphi_0, \Phi; v_1,v_1',v_2,v_2') = \langle a(\frac{\varpi_{v_1}}{\varpi_{v_1'}}) \varphi_0 a(\frac{\varpi_{v_2}}{\varpi_{v_2'}}) \overline{\varphi}_0, E(\Phi,i\tau) \rangle. $$
	Recall the notation $\pi_{i\tau, \xi} = \pi(\xi |\cdot|^{i \tau}, \xi^{-1} |\cdot|^{-i\tau})$. The treatment of $l^h(n(T)(E(\Phi,i\tau)-E_N(\Phi,i\tau)))$ is similar to that of $l^h(n(T)e)$ in the previous section, except that we can take $\theta = 0$. One starts with
	$$ l^h(n(T)(E(\Phi,i\tau)-E_N(\Phi,i\tau))) = \int_{\Re (s) \gg 1} \mathcal{M}(h)(-s) l^{|\cdot|^s}(n(T)(E(\Phi,i\tau)-E_N(\Phi,i\tau))) \frac{ds}{2\pi i} $$
	and
	\begin{equation}
		l^{|\cdot|^s}(n(T)(E(\Phi,i\tau)-E_N(\Phi,i\tau))) =  \Lambda(s+1/2, \pi_{i\tau, \xi}) \prod_v \frac{l^{|\cdot|^s}(n(T_v)W_{\Phi_{i\tau},v})}{L(s+1/2, \pi_{i\tau, \xi, v})},
	\label{EisLExpression}
	\end{equation}
	where
	$$ L(s, \pi_{i\tau, \xi, v}) = L(s+i\tau, \xi_v)L(s-i\tau, \xi_v^{-1}), $$
	$$ \Lambda(s, \pi_{i\tau, \xi, v}) = \Lambda(s+i\tau, \xi) \Lambda(s-i\tau, \xi^{-1}), $$
	and $\Lambda( s, \xi )$ is the complete (${\rm GL}_1$) $L$-function. (\ref{EisLExpression}) has an analytic continuation and admits simple poles at $s = 1/2 \pm i\tau$ only when $\xi = 1 $ is the trivial character and $\tau \neq 0$. We proceed by shifting the contour to $\Re s = 0$ and get
	\begin{eqnarray}
		& &  l^h(n(T)(E(\Phi,i\tau)-E_N(\Phi,i\tau))) \nonumber \\
		&=& \int_{\Re (s)=0} \mathcal{M}(h)(-s) l^{|\cdot|^s}(n(T)(E(\Phi,i\tau)-E_N(\Phi,i\tau))) \frac{ds}{2\pi i} + \label{UnitaryTerm} \\
		& &  \begin{matrix} 1_{\xi = 1} \mathcal{M}(h)(-1/2+i\tau) \Lambda_F^*(1) \Lambda(1+2i\tau, \xi) \prod_v \frac{l^{|\cdot|^{1/2+i\tau}}(n(T_v)W_{\Phi_{i\tau},v})}{L(1+i\tau, \pi_{i\tau, \xi, v})} + \\
			1_{\xi = 1} \mathcal{M}(h)(-1/2-i\tau) \Lambda_F^*(1) \Lambda(1-2i\tau, \xi) \prod_v \frac{l^{|\cdot|^{1/2-i\tau}}(n(T_v)W_{\Phi_{i\tau},v})}{L(1-i\tau, \pi_{i\tau, \xi, v})}. \end{matrix} \label{PolarTerm}
	\end{eqnarray}
	We shall need to bound the contribution of the poles (\ref{PolarTerm}) which doesn't exist in the cuspidal case.

	We first consider the contribution on the line $\Re(s) = 0$, i.e. bound of (\ref{UnitaryTerm}) and give explicit choice of basis $\mathcal{B}(\pi_{i\tau, \xi})$. Note that the operator of taking flat section from $\pi(\xi, \xi^{-1})$ to $\pi_{i\tau, \xi}$ is $K$-equivariant and preserves the inner product, so choosing $\mathcal{B}(\pi_{i\tau, \xi})$ is the same as choosing $\mathcal{B}(\pi(\xi, \xi^{-1}))$. We proceed again as in Section 6.3.1 to 6.3.4, replacing $\pi'$ there by $\pi_{i\tau, \xi}$, taking $\theta = 0$, using Remark \ref{LocalBoundRSLFactorEis} instead of Remark \ref{LocalBoundRSLFactorCusp}, to determine $\mathcal{B}(\pi_{i\tau, \xi})$. Note that the restriction of $\mathcal{B}(\pi_{i\tau, \xi})$ to $K$ doesn't depend on $\tau \in \R$. We therefore get similar bounds for
	$$ \sup_{s \in i\R} \left\lvert l^{|\cdot|^s}(n(T_v)W_{\Phi_{i\tau},v}) \right\rvert, \forall v \mid \infty; \sup_{s \in i\R} \left\lvert \frac{l^{|\cdot|^s}(n(T_v)W_{\Phi_{i\tau},v})}{L(s+1/2, \pi_{i\tau, \xi, v})} \right\rvert, \forall v < \infty $$
	as in the previous section, and deduce the following lemma,
	\begin{lemma}
		We have, for $s\in i\R, \forall \epsilon > 0$,
		$$ |l^{|\cdot |^s}(n(T)(E(\Phi,i\tau)-E_N(\Phi,i\tau)))| \ll_{\epsilon, \varphi_0} |L(s+1/2, \pi_{i\tau, \xi})| |T|^{-1/2+\epsilon} \lambda_{\Phi_{i\tau},\infty}^{1/2+\epsilon}. $$
	\label{MBoundEis}
	\end{lemma}
	\begin{remark}
		We list the differences between bounding local terms here and in the previous section but omit the details of the proof, since they are too similar to each other:
		\begin{itemize}
			\item[(1)]	In the case of Section 6.3.1, we use the bound
			$$ |l^{|\cdot |^s}(n(T_v)W_{\Phi_{i\tau},v})|^2 \leq \dim (K_v W_{\Phi_{i\tau},v}) \lVert W_{\Phi_{i\tau},v} \rVert^2 \int_{F_v^{\times}} \Xi_v(n(-T_v)a(y)n(T_v)) d^{\times}y. $$
			Since $\Xi_v$ is a matrix coefficient, one always has $\Xi_v \leq 1$, so $\Xi_v \leq \Xi_v^{1-\epsilon}$ for any $\epsilon > 0$. We get
			$$ |l^{|\cdot |^s}(n(T_v)W_{\Phi_{i\tau},v})| \ll_{\epsilon } (1+|T_v|)^{-1/2 + \epsilon} (\dim (K_v W_{\Phi_{i\tau},v}))^{1/2} \lVert W_{\Phi_{i\tau},v} \rVert. $$
			Note that we can not directly take $\theta = 0$ in the bounds of
			$$ \int_{F_v^{\times}} \Xi_v(n(-T_v)a(y)n(T_v))^{1-2\theta} d^{\times}y $$
			there, because the implicit constant depending on $\theta$ tends to infinity as $\theta \to 0$.
			\item[(2)]	Every $ \frac{\zeta_v(2)^{1/2}\lVert W_{e,v} \rVert}{\sqrt{L(1, \pi_v' \times \bar{\pi}_v')}} $ should be replaced by $\frac{\zeta_v(2)^{1/2}}{\zeta_v(1)}\lVert W_{\Phi_{i\tau},v} \rVert$ according to Lemma \ref{NormIdenEis} instead of Lemma \ref{NormIden}. Corresponding to $\lVert e \rVert_{X(F)} = 1$, the normalization here is $\lVert E(\Phi, i\tau) \rVert_{{\rm Eis}} = 1$.
			\item[(3)]	Here is the list of bounds: In the case of Section 6.3.1, we have, all for $s \in i\R$,
			$$ |l^{|\cdot|^s}(n(T_v)W_{\Phi_{i\tau},v})| \leq C(\theta, \epsilon) \lambda_{\Phi_{i\tau},v}^{1/2+\epsilon} |T_v|_v^{-1/2+\epsilon} \lVert W_{\Phi_{i\tau},v} \rVert, v \mid \infty; $$
			$$ \left| \frac{l^{|\cdot|^s}(n(T_v)W_{\Phi_{i\tau},v})}{L(s+1/2, \pi_{i\tau, \xi, v})} \right| \leq C(\theta, \epsilon) q_v^{v(\varphi_0)/2+\epsilon} |T_v|_v^{-1/2+\epsilon} \frac{\zeta_v(2)^{1/2}}{\zeta_v(1)} \lVert W_{\Phi_{i\tau},v} \rVert, v < \infty. $$
			In the case of Section 6.3.2 and 6.3.3, we have
			$$ \left| \frac{l^{|\cdot|^s}(n(T_v)W_{\Phi_{i\tau},v})}{L(s+1/2, \pi_{i\tau, \xi, v})} \right| \leq C(\theta, \epsilon) \frac{\zeta_v(2)^{1/2}}{\zeta_v(1)} \lVert W_{\Phi_{i\tau},v} \rVert. $$
			In the case of Section 6.3.4, we have
			$$ \left| \frac{l^{|\cdot|^s}(n(T_v)W_{\Phi_{i\tau},v})}{L(s+1/2, \pi_{i\tau, \xi, v})} \right| = \frac{\zeta_v(2)^{1/2}}{\zeta_v(1)} \lVert W_{\Phi_{i\tau},v} \rVert. $$
		\end{itemize}
	\end{remark}

	We then consider the contribution of (\ref{PolarTerm}). The local factors for which $T_v \neq 0$ are bounded by using (\ref{LocalZetaDecay}) and (\ref{LocalZetaDecayInfinite}). For those for which $v \in \left\{ v_1,v_1',v_2,v_2' \right\}$ and $T_v=0$, we use instead 
	$$ |l^{|\cdot|^{1/2\pm i\tau}}(W_{\Phi_{i\tau},v})| \leq \lVert W_{\Phi_{i\tau},v} \rVert \left( \int_{\text{supp} W_{\Phi_{i\tau},v} } |y| d^{\times}y \right)^{1/2}. $$
	Note that if $\Phi_{i\tau,v}$ lies in a standard basis, then $\text{supp} W_{\Phi_{i\tau},v} \subset \vO_v$; if $\Phi_{i\tau,v}$ lies in a nice basis w.r.t. $a(\varpi_v^n), n \in \N$, then $\text{supp} W_{\Phi_{i\tau},v} \subset \varpi_v^{-n} \vO_v$. We distinguish the $9$ types of positions of $\{ v_1,v_1',v_2,v_2' \}$ described in Proposition \ref{ATranslationTypes}, take into account the choice of local basis described in the beginning of Section 6.3 and (\ref{Mellin<0}), and get
	\begin{equation}
		(\ref{PolarTerm}) \ll_{F,\epsilon, h_0, \pi} \lambda_{\Phi_{i\tau}, \infty}^8 Q^{(\kappa -1)/2+\epsilon} E.
	\label{PolarBound}
	\end{equation}
	In fact, \textbf{Type 1,3,6} give the contribution $E$, other types give less.

	The final part of the argument is a little bit different from the cuspidal case. Because the amplification has ``less'' impact on the Eisenstein part than on the cuspidal part. In fact, in the typical situation (\textbf{Type 1}), for $v \in \{ v_1,v_1',v_2,v_2'\}$, amplification changes the constraint $v(\pi_{i\tau, \xi}) \leq v(\varphi_0) = 0$ into $v(\pi_{i\tau, \xi}) \leq 1$. But $v(\pi_{i\tau, \xi}) = 2 v(\xi)$, the above two constraints are both equivalent to $v(\xi) = 0$. Hence the Eisenstein series $E(\Phi, i\tau)$ giving non zero contribution remain the same with or without amplification and depend only on $\varphi_0$. We may simply insert the convex bound of $L(s+1/2, \pi_{i\tau, \xi})$ into Lemma \ref{MBoundEis}, and combine with (\ref{PolarBound}), (\ref{Mellin>0}) to get
	\begin{align*}
		S_{\text{Eis}}(v_1,v_1',v_2,v_2') &\ll_{F,\epsilon,h_0,\varphi_0} \sum_{\xi \in \widehat{F^{\times} \backslash \A^{(1)}}} \int_{-\infty}^{\infty} \sum_{\Phi \in \mathcal{B}(\pi(\xi, \xi^{-1}))} C(\varphi_0, \Phi; v_1,v_1',v_2,v_2') |T|^{-1/2+\epsilon} \lambda_{\Phi_{i\tau},\infty}^{1/2+\epsilon} \cdot \\
		&  \int_{\Re (s)=0} \mathcal{M}(h)(-s) (1+|s|)^{1/2} (1+|\tau|)^{1/2} \frac{ds}{2\pi i} \frac{d\tau}{4\pi} + \\
		&  \sum_{\xi \in \widehat{F^{\times} \backslash \A^{(1)}}} \int_{-\infty}^{\infty} \sum_{\Phi \in \mathcal{B}(\pi(\xi, \xi^{-1}))} C(\varphi_0, \Phi; v_1,v_1',v_2,v_2') \lambda_{\Phi_{i\tau}, \infty}^8 Q^{(\kappa -1)/2+\epsilon} E \frac{d\tau}{4\pi} \\
		&\ll_{h_0}  (|T|^{-1/2+\epsilon} + Q^{(\kappa -1)/2+\epsilon} E) \lVert P_{\text{Eis}}(\Delta_{\infty}^{10} a(\frac{\varpi_{v_1}}{\varpi_{v_1'}}) \varphi_0 a(\frac{\varpi_{v_2}}{\varpi_{v_2'}}) \overline{\varphi}_0)  \rVert \cdot \left( \text{Trace of } \Delta_{\infty}^{-A} \right)^{1/2},
	\end{align*}
for some $A > 1$. $\lVert P_{\text{Eis}}(\Delta_{\infty}^{10} a(\frac{\varpi_{v_1}}{\varpi_{v_1'}}) \varphi_0 a(\frac{\varpi_{v_2}}{\varpi_{v_2'}}) \overline{\varphi}_0)  \rVert$ can be bounded by some constant depending only on $\varphi_0$, while the trace of laplacian depends only on $\pi$ by Theorem \ref{WeylLaw}, we thus have
$$ S_{\text{Eis}}(v_1,v_1',v_2,v_2') \ll_{F,\epsilon,h_0,\varphi_0} Q^{(\kappa -1)/2+\epsilon} E. $$
This completes the first part of Lemma \ref{EisL}.

In the situation of \textbf{Type 3,8}, the convex bound of $L(s+1/2, \pi_{i\tau, \xi})$ contributes one more factor $E^{1/4}$ while the trace of laplacian is also increased by a factor of $E$, we get
$$ S_{\text{Eis}}(v_1,v_1',v_2,v_2') \ll_{F,\epsilon,h_0,\varphi_0} (Q^{-1/2+\epsilon}E^{1/4}+Q^{(\kappa -1)/2+\epsilon} E) E^{1/2}. $$
But in $\Sigma_3$ the increased factor is killed by the denominator by $E^{-1}$, hence this situation contributes less than the typical situation.

In the situation of \textbf{Type 6}, the convex bound of $L(s+1/2, \pi_{i\tau, \xi})$ contributes one more factor $E^{1/2}$ while the trace of laplacian is also increased by a factor of $E^2$, we get
$$ S_{\text{Eis}}(v_1,v_1',v_2,v_2') \ll_{F,\epsilon,h_0,\varphi_0} (Q^{-1/2+\epsilon}E^{1/2}+Q^{(\kappa -1)/2+\epsilon} E) E. $$
But in $\Sigma_3$ the increased factor is killed by the denominator by $E^{-2}$, hence this situation contributes less than the typical situation.

The other situations listed in Proposition \ref{ATranslationTypes} obviously contribute less. We conclude the second part of Lemma \ref{EisL}.

\address{Han WU \\ HG J 16.2 \\ ETHZ D-Math \\ R\"amistrasse 101 \\ CH-8009 Z\"urich \\ Suisse \\ han.wu@math.ethz.ch}

\end{document}